\documentclass{article}

\usepackage{amsfonts}
\usepackage{graphics}
\usepackage{amsmath}
\usepackage{amscd}
\usepackage{rlepsf}
\usepackage{amssymb}
\usepackage{euscript}
\usepackage{color}

\newtheorem{Theorem}{Theorem}[section]
\newtheorem{Lemma}[Theorem]{Lemma}

\newtheorem{Definition}[Theorem]{Definition}
\newtheorem{Example}[Theorem]{Example}
\newtheorem{Proposition}[Theorem]{Proposition}
\newtheorem{Corollary}[Theorem]{Corollary}

\newtheorem{Remark}[Theorem]{Remark}

\numberwithin{equation}{section}

\newenvironment{Proof}[1][Proof]{\textbf{#1.} }{\ \rule{0.5em}{0.5em}}

\def \Hom {{\rm Hom}}
\def \p {\psi}
\def \cm {{\bf c}}

\def \ad {\mathrm{ad}}
\def \A {\EuScript{A}}

\def \a {\alpha}

\def \GL {{\rm GL}}
\def \B {{\EuScript{B}}}

\def \I {\EuScript{I}}
\def \J {\EuScript{J}}

\def \U {\EuScript{U}}
\def \V {\EuScript{V}}

\def \pexp{\EuScript{P}\exp}

\def \sexp{\EuScript{S}\exp}

\def \Hc {\EuScript{H}}
\def \Xc  {\mathcal{X}}
\def \Yc  {\mathcal{Y}}
\def \Zc  {\mathcal{Z}}
\def \Wc  {\mathcal{W}}

\def \Der {\mathrm{Der}}

\def \ga {\gamma}

\def \e {\epsilon}

\def \w {\omega}

\def \id {\mathrm{id}}

\def \lg {\mathfrak {g}}
\def \Nc {{\EuScript{N}}}
\def \R {\mathbb{R}}
\def \K {\mathcal{K}}
\def \X {\EuScript{X}}

\def \W {\Omega}

\def \G {\Gamma}

\def \f {\phi}

\def \d {\partial}

\def \M {\EuScript{M}}

\def \N {\mathbb N}

\def \Rank {\mathrm{Rank}}

\def \a {\alpha}

\def \b {\beta}

\def \g {\gamma}

\def \W {\Omega}

\def \t {\widetilde}

\def \tn {\otimes}

\def \tr {\triangleright}

\def \Gc {\mathcal{G}}

\def \ra {\xrightarrow}

\def \le {\mathfrak{e}}

\def \Re {\EuScript{R}}

\def \S {\EuScript{S}}

\def \Dc {\EuScript{D}}
\def \D {\mathcal{D}}

\def \T { \EuScript{T} }
\def \Tc { \mathcal{T} }

\def \Q {{\EuScript{Q}}}

\def \Wc {{\EuScript{W}}}

\def \Sc {\mathcal{S}}

\def \H {\mathcal{H}}

\def \Z {\mathbb{Z}}

\def \k {\mathfrak{k}}

\def \op {{{\rm op}}}

\def \Box {{\EuScript{B}}}
\def \Sets {{\mathrm{Sets}}}

\def \h {{\mathrm{h}}}
\def \vm {{\mathrm{v}}}

\def \forall {\textrm{ {for each} } }

\begin{document}
\author{Jo\~{a}o  Faria Martins \\{ \small Departamento de Matem\'{a}tica }\\ {\small
Faculdade de Ci\^{e}ncias e Tecnologia,   Universidade Nova de Lisboa}\\
{\small Quinta da Torre,
2829-516 Caparica,
Portugal } \\ {\it \small jn.martins@fct.unl.pt}\\\\ Roger Picken  \\ {\small Departamento de Matem\'{a}tica,} \\ {\small Instituto Superior T\'{e}cnico, TU Lisbon } \\ {\small Av. Rovisco Pais, 1049-001 Lisboa, Portugal} \\ {\it \small rpicken@math.ist.utl.pt}}

\title{A Cubical {Set} Approach to 2-Bundles with Connection and Wilson Surfaces}

\date{\today}

\maketitle

\begin{abstract}
{In the context of non-abelian gerbes we define a cubical version of
categorical group 2-bundles with connection {over a smooth manifold}. We
define their two-dimensional parallel transport, study its properties,
and define {non-abelian} Wilson surface functionals.}
\end{abstract}


\noindent{\bf 2000 Mathematics Subject Classification:} {53C29 
(primary); 
18D05, 
70S15 (secondary)
}

\section{Introduction}

The aim of this paper is to address {the differential geometry of} (categorical group) 2-bundles over a smooth manifold $M$ and their two dimensional {parallel transport} with a minimal use of two dimensional category theory, {the ultimate goal being to define Wilson {surface} observables.}   The only categorical notion needed is that of an (edge symmetric, strict) double groupoid (with thin structure), which is equivalent to a  crossed module or to a categorical group; see {\cite{BH1,BHS,BH6,BL,BS,BM}}. {We also use} the concept of a cubical set \cite{BH2,J1,GM}, a cubical analogue of a simplicial set, familiar in algebraic topology; see for example \cite{Ma}. 

{Our definition of a 2-bundle with connection will be given in the framework of cubical sets.
Given a {crossed module of Lie groups} $\Gc=(\d\colon E \to G, \tr)$, where $\tr$ is a left action of $G$ on $E$ by automorphisms, the definition of a {cubical} $\Gc$-2-bundle with connection $\B$ over a manifold $M$ is an almost exact cubical analogue of the simplicial version considered in} {\cite{H,BS1,BS2,BrMe}. Following \cite{H,MP}, we will consider a coordinate neighbourhood description of 2-bundles with connection. {For a discussion of the {\it total space} of a 2-bundle see \cite{RS,Bar,Wo}.}

We also define the thin homotopy double groupoid of a smooth manifold $M${, constructed from smooth maps from the square to $M$, identified modulo thin homotopy}. An advantage of the cubical setting over the simplicial setting is that subdivision is very easy to understand. In a cubical $\Gc$-2-bundle with connection, all connection forms are in principle {only locally defined}. Therefore, given a smooth map $[0,1]^2 \to M$, to define its {holonomy (for brevity we will use the term holonomy, instead of the more accurate term, parallel transport)}, one needs to subdivide $[0,1]^2$ into smaller squares, consider all the locally defined holonomies (which we will define and analyse carefully) and patch it all together by using the 1- and  2-transition functions of the cubical $\Gc$-2-bundle, {and the transition data of the connection}. A double groupoid provides a convenient context for doing this type of calculations, {and is easier to handle than the decomposition of $[0,1]
 ^2$ into regions by means of a trivalent embedded graph of \cite{P}.}
{Citing \cite{BHS,BH1}}, {double groupoids trivially have an algebraic inverse to subdivision}. {This was the motivation for our cubical set approach to 2-bundles with connection and their holonomy.}

{We derive the local  two-dimensional holonomy (based on a crossed module), the transition functions and their properties by extending results from our previous study \cite{FMP} of holonomy and categorical holonomy in a principal fibre bundle.  Even though its definition is {apparently non-symmetric in the two variables parametrising $[0,1]^2$}, the local 2-dimensional holonomy which is associated to maps $[0,1]^2 \to M$ is covariant with respect to  the dihedral group of symmetries of the square. This important result (the {Non-Abelian Fubini's {Theorem}}) ultimately follows from the crossed module rules, and would not hold if a pre-crossed module {were} used.

Let $\Gc=(\d\colon E \to G, \tr)$ be a Lie crossed module. {We show (in the final section) that the cubical $\Gc$-2-bundle holonomy} which we define can be associated to oriented embedded 2-spheres $\Sigma\subset M$  yielding an element ${\mathcal W}(\B,\Sigma) \in \ker \d\subset E$ (the Wilson sphere observable) independent of the parametrisation of the sphere {and the chosen coordinate neighbourhoods},  up to acting by elements of $G$. This follows from the invariance of cubical $\Gc$-2-bundle holonomy under thin homotopy {(up to acting by elements of $G$)} and the fact that the mapping class group of the sphere $S^2$ is $\{\pm 1\}$.  This Wilson sphere observable depends only on the equivalence class of the cubical $\Gc$-2-bundle with connection $\B$.
For surfaces other than the sphere embedded in $M$, a holonomy can still be defined but it will a priori {(since the mapping class group is more complicated)} depend on the isotopy type of the parametrisation. We will illustrate this point with the case of Wilson tori.

{An important problem that {follows on} from this construction is the definition of a gauge invariant action in the space of all 2-bundles with connection over a smooth closed 4-dimensional manifold, analogous  to the Chern-Simons action for principal bundles with connection over a 3-dimensional closed manifold - see \cite{B}}. Given that a gauge invariant sphere holonomy was defined, this would permit a physical definition of invariants of knotted spheres in $S^4$ analogue to the Jones polynomial; see for example {\cite{W,Ko,AF,CR}.}

\tableofcontents
\section{Preliminaries}

\subsection{The Box Category and cubical sets }\label{Box}

\subsubsection{Cubical sets}\label{cs}
The box category $\Box$, see \cite{J1,BH2,BH3,BHS,GM}, is defined as the category whose set of objects is the set of standard $n$-cubes $D^n \doteq I^n$, where $I\doteq [0,1]$, and whose set of morphisms is the set of maps generated by the cellular maps $\delta^\pm_{i,n}\colon D^n \to D^{n+1}$, where $i=1,\ldots, n+1$ and $\sigma_{i,n} \colon D^{n+1}\to D^n$, $i=1,\ldots, n+1$. {We have put:}
$$ \delta^-_{i,n} (x_1,\ldots,x_{i-1}, x_{i+1}, \ldots,  
x_n)=(x_1,\ldots,x_{i-1},0, x_{i+1}, \ldots,  x_n)$$

$$ \delta^+_{i,n} (x_1,\ldots,x_{i-1}, x_{i+1}, \ldots,  x_n)=(x_1,\ldots,x_{i-1},1, x_{i+1}, \ldots,  x_n)$$

$$ \sigma_{i,n+1} (x_1,\ldots,x_{n+1})= (x_1,\ldots,x_{i-1}, x_{i+1}, \ldots,  x_n).$$
We will usually abbreviate $\delta_{i,n}=\delta_i$ and $\sigma_{i,n}=\sigma_i$.

\begin{Definition}[Cubical set]
{A cubical object $K$ in the category of sets  (abbrev. ``cubical set'') is a functor $\Box^{\rm op} \to \Sets$, the category of sets; see \cite{BH3,J1,GM,BHS}. Here $\Box^\op$ is the opposite category of the box category $\Box$. A morphism of cubical sets (a cubical map) is a natural tranformation of functors. We can analogously define cubical objects in any category, for example in  the category of smooth manifolds and their smooth maps ({defining cubical manifolds}), or more generally in the  category of smooth spaces \cite{BHo,Ch}.}
\end{Definition}
{Unpacking this definition, we can see that a cubical set $K$ is defined as being an assignment of a set $K_n$ (the set of $n$-cubes) to each $n \in \N$, together with face maps $\d_i^\pm\colon K_n \to K_{n-1}$ and degeneracy maps $\e_i\colon K_{n-1} \to K_n$, where $i\in \{1,\ldots,n\}$ satisfying the cubical  relations:}
\begin{equation}
 \begin{CD}\begin{CD}\d_i^\a\d_j^\b&=\d_{j-1}^\b\d^\a_i \quad\quad (i<j) \\
\e_i\e_j&=\e_{j+1}\e_i \quad\quad\quad (i\leq j) \end{CD}
\quad \quad \quad \begin{CD}\d_i^\a\e_j &=\left \{ \begin{CD} 
\e_{j-1}\d^\a_i  &\quad\quad\quad\quad\quad(i<j)\\ 
\e_j\d^a_{i-1}    &\quad\quad\quad\quad\quad(i>j)\\ 
\id               &\quad\quad\quad\quad\quad(i=j)
 \end{CD}\right.
           \end{CD}
\end{CD}
\end{equation}
Here $\a,\b \in \{-,+\}$. {The {description} of a cubical manifold is analogous, but each $K_n$ is to be a smooth manifold, and all faces and degeneracies are to be smooth.}
{A degenerate cube is a cube in the image of some degeneracy map.}
A cubical set $K$ for which $K_i$ consists only of degenerate cubes if $i>n$ will be called $n$-truncated.

\begin{Definition}[Dihedral cubical set]
{If a cubical set $K$ has an action of the group of symmetries of the $n$-cube  (the $n$-hyperoctahedral group) in each set $K_n$, compatible with the faces and degeneracies in the obvious way, it will be called a {\it dihedral cubical set.}  A cubical map $K\to K'$ between dihedral cubical sets that preserves the actions will be called a dihedral cubical map.}
\end{Definition}
{Dihedral cubical sets are called cubical sets with reversions and interchanges in \cite{GM}. To relate the two definitions, note that the {$n$-hyperoctahedral} group is generated by reflections and interchanges of coordinates, and is therefore isomorphic to $\Z_2^n \rtimes S_n$.}
\begin{Example}
Let $M$ be a manifold. The smooth singular cubical set $C(M)$ of $M$ is given by all smooth maps $D^n\to M$, where $D^n=[0,1]^n$ is the $n$-cube, with the obvious faces and degeneracies, \cite{BH3,BH3}. This is a dihedral cubical set in the obvious way. {We can also see   $C(M)$  as being a cubical object in the category of smooth spaces \cite{BHo}, by giving the set of $n$-cubes the smooth structure of \cite{Ch,BHo}}. 
\end{Example}

\begin{Example}
Analogously, given a smooth manifold $M$, the restricted smooth singular cubical set $C_r(M)$ of $M$ is given by all smooth maps $f\colon D^n\to M$ for which there exists an $\e>0$ such that $f(x_1,x_2,\ldots x_n)=f(0,x_2,\ldots x_n)$ if $x_1 \leq \e$, and analogously for any other face of $D^n$, of any dimension. We will abbreviate this condition by saying that $f$ has a product structure close to the {boundary of the $n$-cube.}  {This condition allows the composition of $n$-cubes to be defined, which we will be needing shortly. In the terminology of  \cite{BH3}, this example is a cubical set with connections and compositions. }
\end{Example}

\subsection{Lie crossed modules}

All Lie groups and Lie algebras are taken to be finite-dimensional. For details on (Lie) crossed modules see, for example, \cite{B1,BM,FM,FMP,B,BL}, and references therein.

\begin{Definition}[{Crossed module {and Lie crossed module}}]\label{LCM}
 A crossed module {(of groups)} ${\Gc= ( \d\colon E \to  G,\tr)}$ is given by a group morphism $\d\colon E \to G$ together with a  left action $\tr$ of $G$ on $E$ by automorphisms, {such that}:

\begin{enumerate}

 \item {$\d(g \tr e)=g \d(e)g^{-1}; \forall g \in G, \forall e \in E,$}

  \item $\d(e) \tr f=efe^{-1};\forall e,f  \in E.$

\end{enumerate}
If both $G$ and $E$ are Lie groups, $\d\colon E \to G$ is a smooth morphism,  and the left action of $G$ on $E$ is smooth then $\Gc$ will be called a Lie crossed module.

\end{Definition}

A morphism $\Gc \to \Gc'$ between the {Lie crossed modules} ${\Gc= ( \d\colon E \to  G,\tr)}$  and {$\Gc'=(\d'\colon E' \to G',\tr')$} is given by a pair of {smooth morphisms} $\f\colon G \to G'$ and  $\psi\colon E \to E'$ making the diagram:
$$ \begin{CD}
  E @>\d>> G \\
 @V \psi VV   @VV \f V \\
  E' @>\d'>> G' \\
   \end{CD}
$$
commutative. In addition we must have {$\psi(g \tr e)=\f(g) \tr' \psi(e)$ for each $ e \in E$ and each $ g \in G$.}

Given a Lie crossed module ${\Gc= ( \d\colon E \to  G,\tr)}$, then the induced Lie algebra map $\d\colon \le \to \lg$, together with the derived action of $\lg$ on $\le$ (also denoted by $\tr$) is a differential crossed module, in the sense of the following definition - see \cite{BS1,BS2,B,BC}.

\begin{Definition}[{Differential crossed module}] {A differential crossed module (or crossed module of Lie algebras) ${\mathfrak{G}=(\d \colon \le \to  \lg,\tr )}$, is given by a Lie algebra morphism $\d\colon \le \to \lg$ together with a left action of $\lg$ on the underlying vector space of  $\le$, {such that}:}
\begin{enumerate}
 \item For any $X \in \lg$ the map $v \in \le \mapsto X \tr v \in \le$ is a derivation of $\le$, in other words $$X \tr [u,v]=[X \tr u,v]+[u, X \tr v];\forall X \in \lg, \forall u ,v\in \le.$$
\item The map $\lg \to \Der(\le)$ from $\lg$ into the derivation algebra of $\le$ induced by the action of $\lg$ on $\le$ is a Lie algebra morphism. In other words:
$$[X,Y] \tr v=X \tr (Y \tr v)-Y \tr(X \tr v); \forall X,Y \in \lg, \forall v\in \le.$$
\item $\d( X \tr v)= [X,\d(v)];\forall X \in \lg, \forall v\in \le.$
\item $\d(u) \tr v=[u,v];  \forall u,v\in \le.$
\end{enumerate}
\end{Definition}
\noindent {Note that  the map $(X,v) \in \lg \times \le \mapsto X \tr v \in \le$  is necessarily bilinear.}

A very useful identity satisfied in any differential crossed module is the following:
\begin{equation}\label{vu}
{\d(u)\tr v=[u,v]=-[v,u]=-\d(v) \tr u, \forall u,v \in \le.}
\end{equation}
This will be used several times in this paper.

{Given a Lie crossed module ${\Gc= ( \d\colon E \to  G,\tr)}$, we will also denote the induced action of $G$ on $\le$ by $\tr$. Finally, }
given a differential crossed module, ${\mathfrak{G}=(\d \colon \le \to  \lg,\tr )}$ there exists a unique crossed module of simply connected Lie groups ${\Gc=(\d \colon E \to  G,\tr )}$ whose differential form is $\mathfrak{G}$, up to isomorphism. The proof of this result is standard Lie theory, together with the lift of the Lie algebra action to a Lie group action, which can be found in \cite{K}, {Theorem} 1.102.

\subsubsection{The edge symmetric double groupoid $\D(\Gc)$ where $\Gc$ is a crossed module}\label{esdg}

The definition of an edge symmetric (strict) double groupoid $\K$ (with thin structure) can be found for example in \cite{BH1,BHS,BHKP,BS}. These are 2-truncated cubical  sets for which the set of 1-cubes $\K_1$ is a groupoid, with set of objects given by the set of 0-cubes,  and also with two partial compositions, vertical and horizontal, in the set $\K_2$ of 2-cubes (squares),  each defining groupoid structures for which the set of objects is the set of 1-cubes. These horizontal and vertical compositions should verify the interchange law:
$$\begin{CD}(k_1 k_2)\\(k_3 k_4)\end{CD} = \left (\begin{CD} k_1\\ k_3\end{CD} \right) \left (\begin{CD} k_2\\ k_4\end{CD} \right),  \forall k_1,k_2,k_3,k_4 \in \K_2, $$
familiar in 2-dimensional category theory,
and be compatible with faces and degeneracies, in the obvious way. In particular, the identity maps  of the vertical and horizontal compositions are given by {degenerate squares}.

There is also an extra condition that should be verified, which is the existence of a thin structure, {meaning that} there exist, among the squares of $\K$, special elements called thin such that:
\begin{enumerate}
 \item Degenerate squares are thin.
 \item Given $a,b,c,d \in \K_1$ with $ab=cd$, there exists a unique thin square $k$ whose boundary is:
$$\begin{CD} &* @>d>> &*\\
              &@A cAA \hskip-1.3cm \scriptstyle{ }  &@AA bA\\
              &* @>> a> &*
  \end{CD}\quad \quad;$$
in other words such that $\d_d(k)=a,\d_r(k)=b,\d_u(k)=d$ and $\d_l(k)=c$, where we have put $\d_d=\d^-_2,\d_r=\d^+_1,\d_u=\d^+_2$ and $\d_l=\d^-_1$.
\item Any {composition} of thin squares is thin.
\end{enumerate}

Let ${\Gc= ( \d\colon E \to  G,\tr)}$ be a crossed module. Given that the categories of crossed modules, categorical groups and double groupoids  with a unique object $*$ are equivalent (see \cite{BH1,BH6,BHS,BS,BL}), we can construct a double groupoid $\D(\Gc)$ out of $\Gc$. 
 The 1-cubes $\D^1(\Gc)$ of $\D(\Gc)$ are given by all elements of $G$, with product as composition, and the {unique} source and target maps to the set $\{*\}$. The $2$-{cubes} $\D^2(\Gc)$ of $\D(\Gc)$, which we will {also} call squares in $\Gc$, have the form:
\begin{equation} \label{square}\begin{CD} &* @>W>> &*\\
              &@A ZAA \hskip-1.3cm \scriptstyle{e}  &@AA YA\\
              &* @>> X> &*
  \end{CD}\end{equation}
where $X,Y,Z,W \in G$ and $e \in E$ is such that $\d(e)^{-1}XY=ZW$. The horizontal and vertical compositions are:
$$ {\begin{CD} &* @>W>> &*\\
              &@A ZAA \hskip-1.3cm \scriptstyle{e}  &@AA YA\\
              &* @>> X> &*
  \end{CD}  \quad \quad \begin{CD} &* @>W'>> &*\\
              &@A YAA \hskip-1.3cm \scriptstyle{e'}  &@AA Y'A\\
              &* @>> X'> &*
  \end{CD}\quad   =  \quad  \begin{CD} &* @>WW'>> &*\\
              &@A ZAA \hskip-1.3cm \scriptstyle{(X \tr e') e} &@AA Y'A\\
              &* @>> XX'> &*
\end{CD}}$$
 and $${
\begin{CD}  \begin{CD} &* @>W'>> &*\\
              &@A Z' AA \hskip-1.3cm \scriptstyle{e'}  &@AA Y'A\\
              &* @>> W> &*
 \end{CD} \\ \begin{CD} &* @>W>> &*\\
              &@A ZAA \hskip-1.3cm \scriptstyle{e}  &@AA YA\\
              &* @>> X> &*
  \end{CD}  \end{CD}\quad \quad  =  \quad \quad  \begin{CD} &* @>W'>> &*\\
              &@A ZZ' AA \hskip-1.3cm \scriptstyle{e Z \tr e'} &@AA YY'A\\
              &* @>> X> &*
  \end{CD}}$$
{The thin structure on $\D(\Gc)$ is given by: } a square is thin if the element of $E$ assigned to it is $1_E$.

{Alternatively the 
thin structure can be given }by introducing the following special degeneracies, usually called connection maps (not to be confused with differential geometric connections) $\ulcorner,\llcorner,\urcorner,\lrcorner\colon \D^1(\Gc)\to \D^2(\Gc)$, whose images are thin:
$$ 
{\ulcorner\left( * \ra{X} *\right)= \quad\quad
\begin{CD} &* @>1_G>> &*\\
              &@A 1_G AA \hskip-1.3cm \scriptstyle{1_E}  &@AA X^{-1} A\\
              &* @>> X> &*
 \end{CD}
\quad \quad, \quad\quad  
\llcorner\left( * \ra{X} *\right)= \quad\quad
\begin{CD} &* @>X>> &*\\
              &@A1_G AA \hskip-1.3cm \scriptstyle{1_E}  &@AA X A\\
              &* @>> 1_G> &*
 \end{CD}}
$$

$${ \urcorner\left( * \ra{X} *\right)= \quad\quad
\begin{CD} &* @>1_G>> &*\\
              &@A X AA \hskip-1.3cm \scriptstyle{1_E}  &@AA 1_G A\\
              &* @>> X> &*
 \end{CD}
\quad \quad, \quad  \quad
\lrcorner\left( * \ra{X} *\right)= \quad\quad
\begin{CD} &* @>X>> &*\\
              &@A X^{-1} AA \hskip-1.3cm \scriptstyle{1_E}  &@AA 1_G A\\
              &* @>> 1_G> &*
 \end{CD}}
$$
{Here we are using results of \cite{BHS,BH1,BH2,BH3,Hi}, where it is shown that 
the existence of special degeneracies, satisfying a set of axioms, is equivalent to the existence of a thin structure. Then  an element of $\D^2(\Gc)$ is thin {if and only if} it is the composition of degenerate squares and the images of special degeneracies; see \cite{Hi,BHS}.
}

The set $\D^2(\Gc)$ is actually a $D_4$-space, where  $D_4$ is the dihedral group of symmetries of the  {square. This can be inferred from the existence of a thin structure.} Consider the following representative elements ${\rho}_{\pi/2},r_x,r_y$ and $r_{xy}$ of $D_4$, where ${\rho}_{\pi/2}$ denotes {anticlockwise} rotation by $90$ degrees, and $r_x,r_y,r_{xy}$ denote reflection in the $y=0$, $x=0$ and $x=y$ axis (recall that these last three elements are generators of $D_4\cong \Z_2^2 \rtimes S_2$). 
Under the action of these elements of $D_4$, the square (\ref{square}) is transformed into, respectively:
\begin{equation*}\begin{CD} &* @>Y^{-1}>> &*\\
              &@A WAA \hskip-1.3cm \scriptstyle{Z^{-1}\tr e}  &@AA XA\\
              &* @>> Z^{-1}> &*
  \end{CD} \quad,\quad\quad \begin{CD} &* @>X>> &*\\
              &@A Z^{-1}AA \hskip-1.3cm \scriptstyle{Z \tr e^{-1}}  &@AA Y^{-1}A\\
              &* @>> W> &*
  \end{CD}\quad\quad,\quad\quad \begin{CD} &* @>W^{-1}>> &*\\
              &@A YAA \hskip-1.3cm \scriptstyle{X\tr e^{-1}}  &@AA ZA\\
              &* @>> X^{-1}> &*
  \end{CD}\quad ,\quad\quad \begin{CD} &* @>Y>> &*\\
              &@A XAA \hskip-1.3cm \scriptstyle{e^{-1}}  &@AA WA\\
              &* @>> Z> &*
  \end{CD}.\end{equation*}
In fact each element of $D_4$ acts on $\D^2(\Gc)$ by automorphisms, though some times permuting the horizontal and vertical multiplications, {or} the order of multiplications. 

The horizontal and vertical inverses $e^{-\h}$ and $e^{-\vm}$ of an element $e \in \D^2(\Gc)$ are given by {$e^{-\h}=r_y(e)$} and {$e^{-\vm}=r_x(e)$;} we will often identify an element of $\D^2(\Gc)$ with the element of  $E$ assigned to it, whenever there is no ambiguity.

There are  two  particular maps  $\Phi,\Phi'_g\colon \D^2(\Gc) \to \D^2(\Gc)$, where $g\in G$, called folding maps, which we would like to make explicit. These are defined as:
$${\Phi\left(\quad\begin{CD} &* @>W>> &*\\
              &@A ZAA \hskip-1.3cm \scriptstyle{e}  &@AA YA\\
              &* @>> X> &*
  \end{CD}  \quad\right) = \quad  \begin{CD} &* @>ZWY^{-1}X^{-1}>> &*\\
              &@A 1_GAA \hskip-2.4cm \scriptstyle{e}  &@AA 1_G A\\
              &* @>> \quad1_G\quad  > &*
  \end{CD}}$$and $${\quad\Phi'_g\left(\quad\begin{CD} &* @>W>> &*\\
              &@A ZAA \hskip-1.3cm \scriptstyle{e}  &@AA YA\\
              &* @>> X> &*
  \end{CD}  \quad\right) =\quad  \begin{CD} &* @> ZWY^{-1}X^{-1} >> &*\\
              &@A gAA \hskip-2.2cm \scriptstyle{g \tr e}  &@AA g A\\
              &* @>>\quad 1_G\quad > &*
  \end{CD}.}$$ 
There also exists an action of $G$ on $\D^2(\Gc)$, which has the form:
$$ g\tr \left (\quad\begin{CD} &* @>W>> &*\\
              &@A ZAA \hskip-1.3cm \scriptstyle{e}  &@AA YA\\
              &* @>> X> &*
  \end{CD} \quad\right )\quad \quad = \quad \quad\quad \quad  \begin{CD} &* @>g  Wg^{-1}>> &*\\
              &@A g Z g^{-1}AA \hskip-1.3cm \scriptstyle{g \tr e}  &@AA gYg^{-1}A\\
              &* @>> gXg^{-1}> &*
  \end{CD} $$
\subsubsection{Flat $\Gc$-colourings, the edge symmetric triple groupoid $\T(\Gc)$ and the nerve $\Nc(\Gc)$ of the crossed module $\Gc$}\label{estg}

{Going one dimension up, following \cite{BHS,BH1,BH2,BH3,BH6}, we can analogously define an edge symmetric triple groupoid $\T(\Gc)$ {of thin 3-cubes in $\Gc$,} from the crossed module ${\Gc= ( \d\colon E \to  G,\tr)}$.
}

{The 1- and 2-cubes  of $\T(\Gc)$ are already defined, being $\T^1(\Gc)=\D^1(\Gc)$ and $\T^2(\Gc)=\D^2(\Gc)$, so let us define the set of {thin} 3-cubes $\T^3(\Gc)$ of $\T(\Gc)$.} Consider the set of assignments  ($\Gc$-colourings of $D^3$) of an element of $G$ to each edge of the standard cube $D^3=[0,1]^3$ in $\R^3$ and of an element of $E$ to each face of $D^3$. Each of these assignments can be mapped to the set of {$\Gc$-colourings of $D^2$, i.e. assignments} of elements of $G$ to the set of edges of the standard square $D^2$ in $\R^2$, and an element of $E$ to its unique face in several different ways, by using the maps $\delta^\pm_i,i=1,2,3$ of \ref{Box}.

Given a $\Gc$-colouring $\cm_2$ of $D^2$, we put $X^\pm_i(\cm_2)=\d^\pm_i(\cm_2)\in G$ as being $\cm_2\circ \delta^\pm_i(D^1)$ {where $i=1,2$}. We also put $e({\cm_2}) = \cm_2(D^2)$.  Analogously, if $\cm_3$ is a $\Gc$-colouring of $D^3$, we put $e^\pm_i(\cm_3)=\d^\pm_i(\cm_3)$ as being the colouring of $D^2$ given by $\cm_3 \circ \delta^\pm_i$ where $i=1,2,3$.
\begin{Definition}[Flat $\Gc$-colouring]

A $\Gc$-colouring $\cm_2$ of $D^2$ is said to be flat if it yields an element of $\D^2(\Gc)$, in the obvious way, in other words if $$\d(e(\cm_2))^{-1}X^-_2(\cm_2)X^+_1(\cm_2)=X^-_1(\cm_2)X^+_2(\cm_2).$$
Analogously,  a $\Gc$-colouring $\cm_3$ of $D^3$ is said to be flat if:
\begin{enumerate}
\item Each restriction $\d^\pm_j({\cm_3})$ of ${\cm_3}$ is a flat $\Gc$-colouring of $D^2$.
\item The following  holds:
\begin{equation} \label{ha} e^+_3({\cm_3})= {\begin{CD}
&\ulcorner(\d^+_2\d^-_1 ({{\cm_3}}))
 &e^+_2({{\cm_3}}) 
&\urcorner(\d^+_2\d^+_1({{\cm_3}}))
 \\ 
&\rho_{\pi/2}(e^-_1({{\cm_3}})) 
& \quad e^-_3({{\cm_3}})\quad  
& r_{xy}(e^+_1({{\cm_3}}))\\  
&\llcorner(\d^-_2\d^-_1 ({{\cm_3}}))  &\quad r_y (e^-_2({{\cm_3}}))\quad &\lrcorner(\d^-_2\d^+_1  ({{\cm_3}}))& \end{CD}}. \end{equation}
{We will call this the {\bf homotopy addition equation}, following the terminology adopted in \cite{BH5}. {Note that we are expressing the fact that the non-abelian composition of five faces of a cube agrees with the sixth face.}}
 
\end{enumerate}
The set $\T^3(\Gc)$  of ({\it thin}) $3$-cubes in $\Gc$  is given by the set of flat $\Gc$-colourings of the $3$-cube.
\end{Definition}

{The set $\T^3(\Gc)$ of {thin} $3$-cubes in $\Gc$ has three interchangeable associative compositions (horizontal, vertical and upwards), as well as boundary maps, $\d^\pm_i, i=1,2,3$. These compositions are induced by the horizontal and vertical composition of squares in $\Gc$ in the unique way such that the boundary maps $\d^\pm_i$ {in the transverse directions} are groupoid morphisms.}
By considering the obvious degeneracies $\e^i\colon \D^1(\Gc) \to \D^2(\Gc), i=1,2$ and 
$\e^i\colon \D^2 (\Gc) \to \T^3(\Gc), i=1,2,3$, obtained by projecting in the $i^{{\rm th}}$ direction (see \ref{cs}), we can see that we  obtain a 3-truncated cubical set $\T(\Gc)$, which is a strict {triple groupoid.}

 {By continuing this process, {one gets} a cubical set $\Nc(\Gc)$, which is called the cubical nerve of $\Gc$. The $n$-cubes of $\Nc(\Gc)$ are given by all $\Gc$-colourings of the $n$-cube $D^n$ such that for each 2- and 3-dimensional face of $D^n$ the restriction of the colouring to it is flat. This is a cubical manifold if $\Gc$ is a Lie crossed module.  The geometric  realisation of $\Nc(\Gc)$ is called the cubical classifying space of $\Gc$; see \cite{BHS,BH4} and \cite{BH5} for the simplicial version. Note that more generally we can take $\Gc$ to  be a crossed module of groupoids \cite{BHS,FMPo}, with completely analogous definitions.}

Note that the homotopy addition equation (\ref{ha}) can be expressed in several different ways by using the $D_4$-symmetry, and applying the maps {$\Phi,\Phi'_g$}. {In particular}, we get the equivalent equation:
\begin{equation}\label{ha2} 
 \Phi'_{\d^-_2\d_1^-(\cm_3)}(e^+_3(\cm_3))=\begin{CD}e^-_1(\cm_3)\quad e^+_2(\cm_3) \quad r_{x	}(e^+_1(\cm_3))\quad  r_x (e^-_2(\cm_3))\\ \Phi(e^-_3(\cm_3))
\end{CD}
\end{equation} 

\subsection{Construction of the thin {homotopy} double groupoid of a smooth manifold}\label{thin2g}

Let $M$ be a smooth manifold.  We now construct the thin {homotopy} double groupoid $\Sc_2(M)$ of $M$. For the analogous construction of the fundamental thin categorical group of a smooth manifold see \cite{FMP}.

\subsubsection{1-paths, 2-paths and 1-tracks} \label{1-Tracks}

\begin{Definition}[1-path]
A 1-path is given by a smooth map $\gamma\colon [0,1] \to M$ such that there exists an $\epsilon >0$ such that $\g$ is constant in $[0,\e] \cup [1-\e,1]$; in the terminology of \cite{CP}, this can be abbreviated by saying that each end point of $\g$ has a sitting instant. Given a 1-path $\g$, define the source and target {or initial and end point} of $\g$ as $\sigma(\g)=\g(0)$ and $\tau(\g)=\g(1)$, respectively. 

\end{Definition}

Given two 1-paths $\g$ and $\f$ with $\tau(\g)=\sigma(\f)$, their concatenation $\g\f$ is  defined in the usual way:
$$(\g\f)(t)=\left \{ \begin{CD} \g(2t), \textrm{ if }t \in [0,1/2] \\ \f(2t-1), \textrm{ if } t \in [1/2,1]\end{CD} \right.$$ 
Note that the concatenation of two 1-paths is also a 1-path, {and in particular is smooth due to the sitting instant condition}.

\begin{Definition}[2-paths]\label{2paths}

A 2-path $\G$ is given by a smooth map $\G\colon [0,1]^2 \to M$ such that 
there exists an $\epsilon >0$ for which:
\begin{enumerate}

\item $\G(t,s)=\G(0,s)$ if $0 \leq t \leq \epsilon$ and $s \in [0,1]$,

\item $\G(t,s)=\G(1,s)$ if  $1-\epsilon \leq t \leq 1$ and $s \in [0,1]$,

\item $\G(t,s)=\G(t,0)$ if $0 \leq s \leq \epsilon$ and $t \in [0,1]$,

\item $\G(t,s)=\G(t,1)$ if  $1-\epsilon \leq s \leq 1$ and $t \in [0,1]$.

\end{enumerate}
{We abbreviate this by saying that $\G$ has a product structure close to the boundary of $[0,1]^2$. }
\end{Definition}

{Given a 2-path $\Gamma$, define the following 1-paths:}
\begin{align*}
\d_l(\G)(s)&=\G(0,s), s \in [0,1],
&\d_r(\G)(s)&=\G(1,s),  s \in [0,1],\\
\d_d(\G)(t)&=\G(t,0), t \in [0,1],
&\d_u(\G)(t)&=\G(t,1), t \in [0,1].
\end{align*}

If $\G$ and $\G'$ are 2-paths such that $\d_r(\G)=\d_l(\G')$ their horizontal concatenation $\G {\circ_\h} \G'$ is defined in the obvious way, in other words:

$$  \big (\G {\circ_\h} \G' \big )(t,s)=\left \{ \begin{CD} \G(2t,s), \textrm{ if }t \in [0,1/2] \textrm{ and }  s \in [0,1] \\ \G'(2t-1,s), \textrm{ if } t \in [1/2,1] \textrm{ and }  s \in [0,1] \end{CD} \right. $$ 
Similarly, if $\d_u(\G)=\d_d(\G')$ we can define a vertical concatenation $\G {\circ_\vm}\G'$ as:
$$\big(\G {\circ_\vm}\G' \big)(t,s)=\left \{ \begin{CD} \G(t,2s), \textrm{ if }s \in [0,1/2]  \textrm{ and }  t \in [0,1]\\ \G'(t,2s-1), \textrm{ if } s \in [1/2,1]  \textrm{ and }  t \in [0,1] \end{CD} \right.$$ 
{Note that again both concatenations are smooth due to the product structure condition.}

\begin{Definition}
Two 1-paths $\f$ and $\g$ are said to be rank-1 homotopic  (and we write $\f\cong_1 \g$) if there exists a 2-path $\G$ such that:
\begin{enumerate}
\item $\d_l(\G)$ and $\d_r(\G)$ are constant.
 \item $\d_u(\G)=\g$ and $\d_d(\G)=\f$.
  \item $\Rank({{\EuScript{D}}}_{v} \G) \leq 1, \forall v \in [0,1]^2.$
\end{enumerate}
{Here ${\EuScript{D}}$ denotes the derivative.} 
\end{Definition} 
Thus if $\g$ and $\f$ are rank-1 homotopic, they have the same initial and end-points.  Note also that rank-1 homotopy is an equivalence relation. Given a 1-path $\g$, the equivalence class to which it belongs is denoted by $[\g]$. Rank-1 homotopy is one of a number of notions of ``thin'' equivalence between paths or loops, and was introduced in 
\cite{CP}, following a suggestion by A. Machado.

We denote the set of $1$-paths of $M$ by $S_1(M)$. The quotient of $S_1(M)$ by the relation of thin homotopy is denoted by $\S_1(M)$. We call the elements of $\S_1(M)$  1-tracks. 
 The {concatenation} of $1$-tracks together with the source and target maps $\sigma, \tau \colon \S_1(M) \to M$, defines a groupoid $\Sc_1(M)$ whose set of morphisms is  $\S_1(M)$ and whose set of objects is $M$.

\subsubsection{2-Tracks}\label{2-Tracks}

We {recall} the notation of \ref{cs}. 
\begin{Definition}\label{thin2}

Two 2-paths $\G$ and $\G'$ are said to be rank-2 homotopic  (and we write $\G\cong_2 \G'$) if there exists a smooth map ${J}\colon {[0,1]^3}  \to M$ such that:

\begin{enumerate}

\item {$J(t,s,0)=\G(t,s),\, J(t,s,1)=\G'(t,s)$ for $s,t\in [0,1]$.}  In other words $J\circ \delta^-_3=\G$ and $J\circ \delta^+_3=\G'.$
 
 \item $J \circ \delta^\pm_i$ is a rank-1 homotopy { from $\G\circ \delta^\pm_i$ to $\G'\circ \delta^\pm_i$, }where $i=1,2$.

\item There exists an $\epsilon >0$ such that $J(t,s,x)={J}(t,s,0)$ if $x\leq \epsilon$ and $s,t \in [0,1]$, and analogously for all the other faces of $[0,1]^3$.  We will describe this condition by saying that $J$ has a product structure close to the boundary of $[0,1]^3$.

 \item  $\Rank ( {{\EuScript{D}}}_{v} {J}) \leq 2$ for any $v \in {[0,1]}^3$. 
\end{enumerate}
\end{Definition} 

Note that {rank-2 homotopy} is an equivalence relation. To prove {transitivity} we need to use the penultimate condition of the previous definition. 
We denote by $S_2(M)$ the set of all 2-paths of $M$. The quotient of $S_2(M)$ by the relation of rank-2 homotopy is denoted by $\S_2(M)$. We call the elements of $\S_2(M)$ 2-tracks.
 If $\G \in S_2(M)$, we denote the equivalence class in $\S_2(M)$ to which $\G$ belongs by $[\G]$.

\subsubsection{Horizontal and vertical compositions of 2-tracks}\label{HVC}

  Suppose that $\G$ and $\G'$ are 2-paths with {$\d_u(\G)\cong_1 \d_d(\G')$.} Choose a rank-1 homotopy ${J}$ connecting $\d_u(\G)$ and $\d_d(\G')$. Then $[\G]\circ_\vm[\G']$ is defined as {$[(\G {\circ_\vm}J) {\circ_\vm}\G'].$} The fact that this composition is well defined in $\S_2(M)$ is not tautological (and was left as an  open problem in \cite{MP}). However this follows immediately from the following lemma proved in \cite{FMP}.

\begin{Lemma}\label{MAIN}
Let $f \colon \d ({{D}}^3) \to M$ be a smooth map such that $\Rank ({{\EuScript{D}}}_v f)\leq 1, \forall \\ v \in \d ({{D}}^3)$. Here ${{D}}^3=[0,1]^3$. Suppose that $f$ is constant in a neighbourhood of each vertex of $\d({{D}}^3) $. In addition, suppose also that in a neighbourhood $I \times [-\epsilon,\epsilon]$ of each edge $I$ of $\d ({{D}}^3)$, {$f(x,t)=\f(x)$, where $(x,t) \in I \times [-\epsilon,\epsilon]$ and $\f\colon I \to M$ is smooth.}
 Then $f$ can be extended to a smooth map $F \colon {{D}}^3 \to M$ such that  $\Rank ({{\EuScript{D}}}_w F)\leq 2, \forall w \in {{D}}^3$. Moreover we can choose  $F$  so that it has a product structure close to the boundary of ${{D}}^3$.
\end{Lemma}

\begin{Remark}
This basically says that any smooth map $f\colon S^2 \to M$ for which the rank of the derivative is {less} than or equal to 1, for each point in $S^2$, can be extended to all of the unit 3-ball, in such a way  that the rank of the derivative of the resulting map at each point is {less than or equal to}  2.
\end{Remark}

Analogously the  horizontal composition of 2-paths descends to $\S_2(M)$.  These compositions are obviously associative, and admit units and inverses. Note that the interchange law is also verified.

Finally, a 2-track $[\G]$ is thin if it admits a representative which is a thin map, in other words for which $\Rank (\Dc_x \G) \leq 1, \forall x \in [0,1]^2$. Lemma \ref{MAIN} implies that if $a,b,c,d\colon [0,1] \to M$ are 1-paths with $[ab]=[cd]$ then there exists a unique 2-track $[\G]$ for which $\d_d([\G])=[a]$, $\d_r([\G])=[b]$, $\d_l([\G])=[c]$ and $\d_u([\G])=[d]$.

 Therefore the following theorem holds:
\begin{Theorem}Let $M$ be a smooth manifold.
The horizontal and vertical {compositions}  in $\S_2(M)$ together with the boundary maps $\d_u,\d_d,\d_l,\d_r\colon \S_2(M) \to \S_1(M)$ define a double groupoid $\Sc_2(M)$, {called the thin homotopy double groupoid of $M$,}  whose set of objects is given by all points of $M$, set of 1-morphisms by the set $\S_1(M)$ of 1-tracks on $M$, and set of 2-morphisms  by all 2-tracks in $\S_2(M).$ In addition, { $\S_2(M)$ admits a thin structure given by:} a 2-track is thin if it admits a representative whose  derivative has rank {less than or equal to}  1 (in other words if it is thin as a smooth map).
\end{Theorem}

\begin{Remark}
Another possible argument to prove that the compositions of 2-tracks are well defined is to adapt the arguments in \cite{BH1,BHS,BH2,BH3,BHKP}, which lead to the construction of the fundamental double groupoid of a triple of spaces and of a Hausdorff space (and can be continued to define the homotopy $\w$-groupoid of a filtered space).  The same technique therefore leads to the construction of the fundamental $\w$-groupoid of a smooth manifold. Details will appear elsewhere.
\end{Remark}

This construction should be compared with \cite{HKK,BHKP}, where the thin strict 2-groupoid of a Hausdorff space was defined, using a different notion of thin equivalence (factoring through a graph).  For analogous non-strict constructions see \cite{M,BS1,MP}.

\subsection{Connections and categorical connections in principal fibre bundles}

To approach non-abelian integral calculus based on a crossed module, it is convenient (since the proofs are slightly {easier}) to consider categorical connections in principal fibre bundles. For details of this approach see \cite{FMP}. {For a treatment of non-abelian integral calculus based on a crossed module, using forms on the base space of the principal bundle, see \cite{SW1,SW2,SW3,FMP2}. }

\subsubsection{Differential crossed module valued forms}
Let $M$ be a  smooth manifold with its Lie algebra of vector fields denoted by $\X(M)$. Consider also  a differential crossed module ${\mathfrak{G}=(\d\colon \le \to \lg,\tr)}$. In particular the map $(X,e) \in \lg \times \le  \mapsto X \tr e \in \le$ is bilinear. 

Let $a \in \A^n(M,\lg )$  and $b \in \A^m(M,\le)$ be $\lg$- and $\le$-valued (respectively) differential forms on  $M$. We define $a \tn^\tr b  $ as being the  $\le$-valued covariant tensor field on $M$ such that
$$(a \tn^\tr b)(A_1,\ldots, A_n,B_1,\ldots B_m)= a (A_1,\ldots, A_n)\tr b(B_1,\ldots, B_m);A_i,B_j \in \X(M).$$
{We also define an alternating tensor field  $a \wedge^\tr b \in \A^{n+m}(M,\le)$, being given by }
$${ a \wedge^\tr b=\frac{(n+m)!}{n! m!}{\rm Alt}(a  \tn^\tr b).}$$
Here ${\rm Alt}$ denotes the natural projection from the vector space of $\le$-valued covariant tensor fields on $M$ onto the vector space of $\le$-valued differential forms on $M$.
For example, if  $a \in \A^1(M,\lg)$ and $b \in \A^2(M,\le)$, then $a \wedge^\tr b$ satisfies:
\begin{equation}\label{BB}
(a \wedge^\tr b)(X,Y,Z)=a(X) \tr b(Y,Z)+a(Y) \tr b(Z,X)+a(Z) \tr b(X,Y),
\end{equation}
where $X,Y,Z \in \X(M).$  

\subsubsection{Categorical connections in principal fibre bundles}
{In \cite{FMP} we defined categorical connections in terms of differential forms on the total space of a principal fibre bundle. }
Let $M$ be a smooth manifold and $G$ a Lie group with Lie algebra $\lg$. Let also $\pi\colon P \to M$ be a smooth principal $G$-bundle over $M$. Denote the fibre at each point $x \in M$ as $P_x\doteq \pi^{-1}(x)$.

\begin{Definition}
{Let  ${\mathcal G}=(\d\colon E\to  G,\triangleright)$  
be a Lie crossed module, where $\tr$ is a Lie group left action of $G$ on $E$ by automorphisms. Let also  $\mathfrak{G}=(\d\colon \le\to \lg,\tr)$ be the associated differential crossed module. A $\Gc$-categorical connection on $P$ is a pair $(\w,m)$, where $\w$ is a connection 1-form on $P$, {i.e.} $\w\in \A^1(P,\lg)$ is a $1$-form on $P$ with values in $\lg$  such that: 
\begin{itemize}
\item $R_g^*(\w)=g^{-1}\w g , \forall g \in G,$ (i.e. $\w$ is $G$-equivariant)
\item $\w(A^\#)=A,\forall A \in \lg$;
\end{itemize}
{where $A^\#$ denotes the vertical vector field associated to $A$ coming from the $G$-action on $P$, }
and $m\in \A^2(P,\le)$ is a 2-form on $P$ with values in $\le$, the Lie algebra of $E$, such that:\begin{itemize}
\item {$m$ is $G$-equivariant, in the sense that $R_g^*(m)=g^{-1} \tr m$ for each $g \in G$.} 
\item {$m$ is horizontal,} in other words:
 $$m(X,Y)=m(X^H,Y^H), \textrm{ for each } X,Y \in \X(P).$$  
In particular $m(X_u,Y_u)=0$ if either of the vectors $ X_u,Y_u \in T_u P$ is  vertical, where $u \in P$. Here the map $X \in \X(P) \mapsto X^H \in \X(P)$ denotes the horizontal projection of vector fields on $P$ with respect to the connection 1-form $\w$.
\end{itemize}
{Finally $(\w,m)$ satisfies} the ``vanishing of the fake curvature condition'' \cite{BS1,BS2,BrMe}:
\begin{equation}
\d(m)=\W,
\end{equation} 
where $\W=d\w +\frac{1}{2} \w \wedge^\ad \w \in \A^2(P,\lg)$ is the curvature 2-form of $\w$.}
\end{Definition}

\subsubsection{The  categorical curvature 3-form of a $\Gc$-categorical connection}\label{curv}

Let $P$ be a principal $G$-bundle over $M$. Let $\w \in \A^1(P,\lg)$ be a connection 1-form on $P$. Given an $n$-form $a$ 
on $P$, the exterior covariant derivative of $a$ is given by $$Da=d a \circ ({H \times H \ldots \times  H} ).$$
Let $\W\in \A^2(P,\lg)$ be the ($G$-equivariant) curvature 2-form of the connection $\w$. It can be defined as the exterior covariant derivative $D\w$ of the connection 1-form $\w$ and also by the Cartan structure equation $\W=d\w +\frac{1}{2} \w \wedge^\ad \w$.
It is therefore natural to define:
\begin{Definition}[{Categorical curvature}]
Let ${\Gc=(\d\colon E \to G,\tr)}$  {be a crossed module of Lie groups}, and let $P \to  M$ be  a smooth principal $G$-bundle. 
The {categorical curvature} 3-form {or 2-curvature 3-form} of a $\Gc$-categorical connection $(\w,m)$  on $P$ is defined as 
$\M=Dm$, where the exterior covariant derivative $D$ is taken with respect to $\w$.
\end{Definition}

The following equation is an analogue of Cartan's structure equation.
\begin{Proposition}[{Categorical structure equation}]\label{cse}
We have:
$\M=dm {+}\w\wedge^\tr m$. {In particular {the} 2-curvature 3-form $\M$ is $G$-{equivariant}, in other words: $R_g^*(\M)=g^{-1} \tr \M,$
for each $g \in G$.}
\end{Proposition}
{This categorical structure equation follows directly from the following natural lemma, easy to prove; see \cite{FMP}:}
\begin{Lemma}\label{ext}
Let $a$ be a $G$-{equivariant} horizontal $n$-form in $P$.  Then $Da=da+\w\wedge^\tr a$.
\end{Lemma}

 Recall that the  {usual} Bianchi identity can be written as $D\W=0$, which is the same as saying that $d\W+\w\wedge^\ad\W=0$.

\begin{Corollary}
The 2-curvature 3-form of a categorical connection is $\mathfrak{k}$-valued, where $\k$ is the Lie algebra of $K=\ker(\d)$.
\end{Corollary}
\begin{Proof}
We have $\d(\M)=\d(d m +\w\wedge^\tr m)=d\W+\w\wedge^\ad \W=0,$ by the Bianchi identity.
\end{Proof}

The 2-curvature $3$-form of a categorical connection satisfies the following.
\begin{Proposition}[{2-Bianchi 	identity}]
Let $\M\in \A^3(P,\le)$ be the {2-curvature} 3-form of $(\w,m)$. Then the exterior covariant derivative $D \M$ of $M$ vanishes,
which by {Lemma} \ref{ext} is the same as:
$d \M +\w \wedge^{\tr} \M=0.$
\end{Proposition}

\subsubsection{Local form}\label{lf}

Let $P\to M$ be a {principal $G$-bundle} with a categorical connection $(\w,m)$.  Let $\{U_i\}$ be an open cover of $M$, with local sections $\sigma_i \colon U_i \to P$ of $P$. The local form of $(\w,m)$ is given by the forms $(\w_i,m_i)$, where $\w_i =\sigma_i^*(\w)$ and $m_i=\sigma_i^*(m)$, and we have $\d(m_i)=d\w_i+\frac{1}{2}\w_i \wedge^\ad \w_i=\W_i=\sigma_i^*(\W)$, and also $\w_j=g_{ij}^{-1} \w_i g_{ij}+g_{ij}^{-1} dg_{ij}$ and $m_j=g_{ij}^{-1} \tr m_i$. Here $\sigma_ig_{ij}=\sigma_j$. Conversely, given forms $\{(\w_i,m_i)\}$ satisfying these conditions then there exists a unique categorical connection $(\w,m)$ in $P$ whose local form ({with respect to} the given sections $\sigma_i $) is $(\w_i,m_i)$.  

Note that locally the {2-curvature} $3$-form of a categorical connection reads $\M_i=dm_i+\w_i\wedge^\tr m_i$, {with $\M_j=g_{ij}^{-1} \tr \M_i$ }and the 2-Bianchi identity is $d\M_i+\w_i\wedge^\tr \M_i=0$.

\subsection{Holonomy and categorical holonomy in a principal fibre bundle}\label{hln}

Let $P$ be a principal $G$-bundle over the manifold $M$. Let $\w \in \A^1(P,\lg)$ be a connection on $P$.  Recall that $\w$ determines a parallel {transport} along smooth curves. Specifically, given $x \in M$ and a smooth curve $\gamma\colon [0,1] \to M$, with $\gamma(0)=x$, then there exists a smooth map:
 $$ (t,u) \in [0,1] \times P_x \mapsto {\mathcal{H}}_\w(\ga,t,u) \in  P,$$ 
uniquely defined by  the conditions:

\begin{enumerate} 
\item $\frac{d}{d t} {\mathcal{H}}_\w(\ga,t,u)=\left(\t{\frac{d}{dt} \ga(t)}\right)_{{\mathcal{H}}_\w(\ga,t,u) };\forall t \in [0,1], \forall u \in P_x,$ {where $\,\t{}\,$ denotes the horizontal lift,} 

\item ${\mathcal{H}}_\w(\ga,0,u)=u; \forall u \in P_x. $

\end{enumerate}
In particular this implies that ${\mathcal{H}}_\w(\ga,t)$, {given by $u\mapsto \mathcal{H}_\w(\ga,t,u)$}, maps $P_x$ bijectively into  $P_{\ga(t)}$, for any $t \in [0,1]$.  
We will also use the notation 
${\mathcal{H}}_\w(\ga,1,u)\doteq u \g. $
Therefore if $\g$ and $\g'$ are such that $\g(1)=\g'(0)$ we have:
$ (u\g)\g'=u(\g\g')$.
Recall that the parallel transport is $G$-equivariant, in other words:
$$\H_\w(\ga,t,ug)=\H_\w(\ga,t,u)g, \forall g \in G, \forall u \in P_x.$$

\subsubsection{A form of the Ambrose-Singer Theorem}\label{AST}

Let $M$ be a smooth manifold. Let $D^n\doteq [0,1]^n$ be the $n$-cube, where $n \in \N$. A map $f\colon  D^n \to M$ is said to be smooth if {its partial derivatives of any order} exist and are continuous as maps $D^n \to M$.

The well known relation between curvature and parallel transport can be summarised in the {following lemma,} proved for instance in \cite{FMP,SW2}.
\begin{Lemma}\label{Main1}
Let $G$ be a Lie group with Lie algebra $\lg$. Let $P$ be a smooth principal $G$-bundle over the manifold $M$. Consider a smooth map $\G\colon [0,1]^2 \to M$.  For each $s,t \in [0,1]$, define the curves $\g_s,\g^t \colon [0,1] \to M$ as $\g_s(t)=\g^t (s)=\G(t,s)$.
Consider a connection $\w \in \A^1(P,\lg)$. 
Choose $u \in P_{\g^0(0)}$, and let  $u_s={\mathcal{H}}_\w(\g^0,s,u)$, and analogously $u^t={\mathcal{H}}_\w(\g_0,t,u)$ where $s,t \in  [0,1]$.  The following holds for each $s,t \in [0,1]$:
\begin{equation}
\w\left (\frac{\d}{\d s} { {\mathcal{H}}_\w(\g_s,t,u_s)}\right)=\int_{0}^{t} \W \left (\t{\frac{\d}{\d t'}\g_s(t')},\t{\frac{\d}{\d s}\g_s(t')} \right)_{{\mathcal{H}}_\w(\g_s,t',u_s) }  d t',
\end{equation}
and by reversing the roles of $s$ and $t$ we also have:
\begin{equation}
\w\left (\frac{\d}{\d t} { {\mathcal{H}}_\w(\g^t ,s,u^t)}\right)=-\int_{0}^{s} \W \left (\t{\frac{\d}{\d t}\g_{s'}(t)},\t{\frac{\d}{\d s'}\g_{s'}(t)} \right)_{{\mathcal{H}}_\w(\g^t ,s',u^t) }  d s'.
\end{equation}
\end{Lemma}

Continuing the notation of the previous lemma, define the elements ${\stackrel{\w}{g}_\G}{(u,t,s)}$ by the rule:
$$\H_\w(\g^t ,s,u^t){\stackrel{\w}{g}_\G}{(u,t,s)}=\H_\w(\g_s,t,u_s).$$
Therefore
$$u{\stackrel{\w}{g}_\G}(u,t,s)=\H_\w(\hat{\g},1,u)$$
where $\hat{\g}$ is the curve $\hat{\g}=\d \G'$, starting in $\G(0,0)$ and oriented clockwise, and $\G'$ is the  truncation of $\G$ such that $\G'(t',s')=\G(t't,s's)$, for $0\leq s',t' \leq 1$.

By using the fact that $\frac{\d}{\d t} \H_\w(\g_s,t,u_s)$ is horizontal it follows that:
$$\w\left(\frac{\d}{\d t} \left (\H_\w(\g^t ,s,u^t){\stackrel{\w}{g}_\G}{(u,t,s)} \right) \right) =0.$$
{Thus,} by using the Leibniz rule together with the fact that $\w$ is a connection $1$-form,
$$ ({\stackrel{\w}{g}_\G}{(u,t,s)})^{-1}\w\left(\frac{\d}{\d t}  \H_\w(\g^t ,s,u^t)\right){\stackrel{\w}{g}_\G}{(u,t,s)}+({\stackrel{\w}{g}_\G}{(u,t,s)})^{-1} \frac{\d}{\d t}{\stackrel{\w}{g}_\G}{(u,t,s)}=0.$$

{Therefore:
\begin{align}\label{gA}
\frac{\d}{\d t}{\stackrel{\w}{g}_\G}{(u,t,s)}&=\left(\int_{0}^{s} \W \left (\t{\frac{\d}{\d t}\g_{s'}(t)},\t{\frac{\d}{\d s'}\g_{s'}(t)} \right)_{{\mathcal{H}}_\w(\g^t ,s',u^t) }  d s'\right){\stackrel{\w}{g}_\G}{(u,t,s)}.
\end{align}
}
Analogously we have (since $\frac{\d}{\d s} \H_\w(\g^t ,s,u^t)$ is horizontal):
\begin{align}\label{gB}
\frac{\d}{\d s} {\stackrel{\w}{g}_\G}{(u,t,s)} &={\stackrel{\w}{g}_\G}{(u,t,s)}\int_0^{t}\W \left (\t{\frac{\d}{\d t'}\g_s(t')},\t{\frac{\d}{\d s}\g_s(t')} \right)_{{\mathcal{H}}_\w(\g_s,t',u_s) }  d t'.
\end{align}

\subsubsection{Categorical holonomy in a principal fibre bundle}\label{chpfb}

Let $P$ be a principal fibre bundle with a $\Gc$-categorical connection $(\w,m)$. Here  ${\mathcal G}=(E \ra{\d} G,\triangleright)$  
is a Lie crossed module, where $\tr$ is a Lie group left action of $G$ on $E$ by automorphisms. Let also  $\mathfrak{G}=(\d\colon \le\to \lg,\tr)$ be the associated differential crossed module. 

As before,  for each smooth map $\G \colon [0,1]^2 \to M$,  let $\g_s(t)=\g^t (s)=\G(t,s)$. Let $a=\G(0,0)$. {Let also $u\in P_a$, $u_s=\H(\g^0,s,u)$ and $u^t=\H(\g_0,t,u).$} Define the function ${\stackrel{(\w,m)}{e_\G}}\colon P_a \times  [0,1]^2\to E$ as being the solution of the differential equation:

\begin{equation}\label{E}
{\frac{\d}{\d s} {\stackrel{(\w,m)}{e_\G}}{(u,t,s)}={\stackrel{(\w,m)}{e_\G}}{(u,t,s)} \int_{0}^{t} m \left (\t{\frac{\d}{\d t'}\g_s(t')},\t{\frac{\d}{\d s}\g_s(t')} \right)  _{{\mathcal{H}}_\w(\g_s,t',u_s) } dt',}
\end{equation}
with initial condition ${\stackrel{(\w,m)}{e_\G}}(u,t,0)=1_E,$ for each $t \in [0,1]$. Let ${\stackrel{(\w,m)}{e_\G}}(u) \doteq {\stackrel{(\w,m)}{e_\G}}(u,1,1)$. Compare with equations (\ref{gA}) and (\ref{gB}).
The apparently non-symmetric way the horizontal and vertical directions are treated will be dealt with later.

Given a smooth map $\G\colon [0,1]^2 \to M$, define:
 $$\Xc_\G=\g_0, \quad \Yc_\G=\g^1, \quad\Zc_\G=\g^0 \quad \textrm{ and }\quad \Wc_\G=\g_1.$$

\begin{Theorem}[{Non-Abelian Green's Theorem, bundle form}]\label{gtb}
For any $u \in P_a$ we have:
$$\H_\w(\Xc_\G \Yc_\G,1,u)\d{\Big(\stackrel{(\w,m)}{e_\G}(u)\Big)}=\H_\w(\Zc_\G \Wc_\G,1,u),$$
or, {in the other notation of section \ref{hln}},
$$u\Xc_\G \Yc_\G\d{\Big(\stackrel{(\w,m)}{e_\G}{(u)}\Big)}=u\Zc_\G \Wc_\G.$$
\end{Theorem}
\begin{Proof}
Let $k_x=\H_\w(\g^1,x,u^1)$ and $l_x=\H_\w(\g_x,1,u_x)$. Let $x\mapsto g_x \in G$ be defined as $k_xg_x=l_x$. We have, since {$(\frac{d}{d x} k_x) g_x$} is horizontal:
\begin{align*}
\w\left ( \frac{d}{d x} (k_x g_x)\right)=\w\left (k_x\frac{d}{d x}  g_x\right)=\w\left (k_x g_x g_x^{-1}\frac{d}{d x}  g_x\right)= g_x^{-1}\frac{d}{d x}  g_x.
\end{align*}
On the other hand:
\begin{align*}{
\w\left ( \frac{d}{d x} (k_x g_x)\right) =\w\left (\frac{d}{d x} l_x\right) =\int_{0}^{1} \W \left (\t{\frac{\d}{\d t}\g_x(t)},\t{\frac{\d}{\d x}\g_x(t)} \right)_{{\mathcal{H}}_\w(\g_x,t,u_x) }  d t.}
\end{align*}
Therefore
\begin{equation}\label{ddxg}
\frac{d}{d x}  g_x=g_x\int_{0}^{1} \W \left (\t{\frac{\d}{\d t}\g_x(t)},\t{\frac{\d}{\d x}\g_x(t)} \right)_{{\mathcal{H}}_\w(\g_x,t,u_x) }  d t. 
\end{equation}
This is a differential equation satisfied {also} by {$x\mapsto \d({\stackrel{(\w,m)}{e_\G}}(u,x,1))$, by the vanishing of the fake curvature condition $\d(m)=\W$}, and both have the same initial conditions.
\end{Proof}

{
{Note that it follows from the {Non-Abelian Green's Theorem} that:}
\begin{equation}\label{gtts}
{{\H_\w(\g^t,s,u^t)}\d \Big ({\stackrel{(\w,m)}{e_\G}(u,t,s)\Big)= {{\mathcal{H}}_\w(\g_s,t,u_s) }},\textrm{ for each } t,s \in  [0,1].}
\end{equation}
}

\begin{Lemma}[Vertical multiplication]\label{vm}
We have: $$\stackrel{(\w,m)}{e_{\G {\circ_\vm}\G'}}{(u)}={\stackrel{(\w,m)}{e_\G}}(u) \stackrel{(\w,m)}{e_{\G'}}({u\Zc_\G}).$$
{Here $\G,\G'\colon [0,1]^2 \to M$ are smooth maps such that $\d_u(\G)=\d_d(\G')$ and moreover $\G {\circ_\vm}\G'$ is smooth.}
\end{Lemma}
\begin{Proof}
Obvious from the definition.
\end{Proof}

\begin{Lemma}[Vertical inversion]\label{vi}
We have: $${\stackrel{(\w,m)}{e_\G}}(u) \stackrel{(\w,m)}{e_{\G^{{-\vm}}}}({u\Zc_\G})=1_E.$$
Here $\G^{-\vm}$ denotes the obvious vertical reversion of $\G\colon [0,1]^2 \to M$.
\end{Lemma}
\begin{Proof}
Obvious from the definition.
\end{Proof}

\begin{Lemma}[Horizontal multiplication]\label{hm}
We have: $$\stackrel{(\w,m)}{e_{\Phi {\circ_\h} \Psi}}(u)=\stackrel{(\w,m)}{e_{\Psi}}(u\Xc_{\Phi})\stackrel{(\w,m)}{e_{\Phi}}(u).$$
{Here $\Phi,\Psi'\colon [0,1]^2 \to M$ are smooth maps such that $\d_r(\Phi)=\d_l(\Psi)$ and moreover $\Phi {\circ_\h} \Psi$ is smooth.}
\end{Lemma}
\begin{Proof}
Let $\G=\Phi{\circ_\h} \Psi$. As before put $\f_s(t)=\f^t(s)=\Phi(t,s)$ and $\psi_s(t)=\psi^t(s)=\Psi(t,s)$. We have:
\begin{align*}
\frac{\d}{\d s} &\left( { \stackrel{(\w,m)}{e_{\Psi}}(u\Xc_\Phi,1,s)}{\stackrel{(\w,m)}{e_{\Phi}}{(u,1,s)}}\right)\\&={ \stackrel{(\w,m)}{e_{\Psi}}(u\Xc_\Phi,1,s)}{\stackrel{(\w,m)}{e_{\Phi}}{(u,1,s)}}\left (\int_{0}^{1} m \left (\t{\frac{\d}{\d t}\phi_s(t)},\t{\frac{\d}{\d s}\phi_s(t)} \right)_{{\mathcal{H}}_\w(\phi_s,t,u_s) } dt\right)\\&+{ \stackrel{(\w,m)}{e_{\Psi}}(u\Xc_\Phi,1,s)}\left(\int_{0}^{1} m \left (\t{\frac{\d}{\d t}\psi_s(t)},\t{\frac{\d}{\d s}\psi_s(t)} \right)  _{{\mathcal{H}}_\w(\psi_s,t,(u\Xc_\Phi)_s) } dt\right){\stackrel{(\w,m)}{e_{\Phi}}{(u,1,s)}}\\&=Q+W.
\end{align*}
Here $(u\Xc_\Phi)_s=\H_\w(\Zc_\Psi,s,u\Xc_\Phi)$.
Let us analyse each term separately. We have:

\begin{align*}
Q=\stackrel{(\w,m)}{e_{\Psi}}(u\Xc_\Phi,1,s){\stackrel{(\w,m)}{e_{\Phi}}{(u,1,s)}}\left (\int_{0}^{\frac{1}{2}} m \left (\t{\frac{\d}{\d t}\g_s(t)},\t{\frac{\d}{\d s}\g_s(t)} \right)_{{\mathcal{H}}_\w(\g_s,t,u_s) } dt\right)
\end{align*}
{where $\g_s(t) = \Phi{\circ_\h} \Psi(t,s)$.}
{On the other hand:}
\begin{multline}
W={ \stackrel{(\w,m)}{e_{\Psi}}(u\Xc_\Phi,1,s)}{\stackrel{(\w,m)}{e_{\Phi}}{(u,1,s)}}\\  \left (\d({\stackrel{(\w,m)}{e_{\Phi}}{(u,1,s)}})^{-1}\tr \left(\int_{0}^{1} m \left (\t{\frac{\d}{\d t}\psi_s(t)},\t{\frac{\d}{\d s}\psi_s(t)} \right)  _{{\mathcal{H}}_\w(\psi_s,t,(u\Xc_\Phi)_s) } dt\right)\right),
\end{multline}
and therefore
\begin{align*}
W&={ \stackrel{(\w,m)}{e_{\Psi}}(u\Xc_\Phi,1,s)}{\stackrel{(\w,m)}{e_{\Phi}}{(u,1,s)}}\\& \quad\quad\quad\quad\left(\int_{0}^{1} m \left (\t{\frac{\d}{\d t}\psi_s(t)},\t{\frac{\d}{\d s}\psi_s(t)} \right)  _{{\mathcal{H}}_\w(\psi_s,t,   (u\Xc_\Phi)_s \d({\stackrel{(\w,m)}{e_{\Phi}}{(u,1,s)}})) } dt\right)\\
&={ \stackrel{(\w,m)}{e_{\Psi}}(u\Xc_\Phi,1,s)}{\stackrel{(\w,m)}{e_{\Phi}}{(u,1,s)}}\left(\int_{0}^{1} m \left (\t{\frac{\d}{\d t}\psi_s(t)},\t{\frac{\d}{\d s}\psi_s(t)} \right)  _{{\mathcal{H}}_\w(\psi_s,t,u_s\phi_s) } dt\right)\\
&={ \stackrel{(\w,m)}{e_{\Psi}}(u\Xc_\Phi,1,s)}{\stackrel{(\w,m)}{e_{\Phi}}{(u,1,s)}}\left(\int_{\frac{1}{2}}^{1} m \left (\t{\frac{\d}{\d t}\g_s(t)},\t{\frac{\d}{\d s}\g_s(t)} \right)  _{{\mathcal{H}}_\w(\g_s,t,u_s) } dt\right).
\end{align*}
Therefore both sides of the equation of the lemma satisfy the same differential equation, and they have the same initial condition.
\end{Proof}

\begin{Lemma}[Horizontal inversion]\label{hi}
We have: 
$$\stackrel{(\w,m)}{e_{\G^{{-\h}}}}({u\Xc_\G}){{\stackrel{(\w,m)}{e_\G}}}(u)=1_E,$$
{where $\G^{-\h}$ denotes the obvious horizontal {reversion} of $\G\colon [0,1]^2 \to M$.}\end{Lemma}
\begin{Proof}
Analogous to the proof of the previous result.
\end{Proof}

\begin{Lemma}[Gauge transformations]\label{gt}
We have: 
$${{\stackrel{(\w,m)}{e_\G}}}(ug)=g^{-1} \tr {{\stackrel{(\w,m)}{e_\G}}}(u).$$
\end{Lemma}
\begin{Proof}
Analogous to the proof of the previous result.
\end{Proof}

\subsubsection{The Non-Abelian Fubini's Theorem}\label{naft}
{We continue with the notation of \ref{chpfb}.} Again  let $\G\colon [0,1]^2\to M$ be a smooth map, $a=\G(0,0)$ and $u \in P_a$.
Define  $\stackrel{(\w,m)}{f_\G}(u,t,s)$ by the differential equation: 
\begin{equation}\label{F}
{\frac{\d}{\d t} \stackrel{(\w,m)}{f_\G}(u,t,s)=\stackrel{(\w,m)}{f_\G}(u,t,s)\int_{0}^{s} m \left (\t{\frac{\d}{\d s'}\g^t(s')},\t{\frac{\d}{\d t}\g^t(s')} \right)_{\H_\w(\g^t,s',u^t)} ds',}
\end{equation}
{with initial condition $\stackrel{(\w,m)}{f_\G}(u,0,s)=1_E,$ for each $s\in [0,1]$.
Note that the differential equation for $\stackrel{(\w,m)}{f_\G}$ is obtained from the differential equation for $\stackrel{(\w,m)}{e_\G}$, equation (\ref{E}), by reversing the roles of $s$ and $t$. Let $\stackrel{(\w,m)}{f_\G}(u,1,1) \doteq {\stackrel{(\w,m)}{f_\G}}(u)$. The following holds.}

\begin{Theorem}[Non-abelian Fubini's Theorem, bundle form]\label{naftbf}
$${\stackrel{(\w,m)}{e_\G}}(u)\stackrel{(\w,m)}{f_\G}(u)=1.$$
\end{Theorem}
\begin{Proof}
In fact we show for every $t,s\in [0,1]$:
\begin{equation}\label{Fubts}
{\stackrel{(\w,m)}{e_\G}}(u,t,s)\stackrel{(\w,m)}{f_\G}(u,t,s)=1.
\end{equation}
In the following put ${\stackrel{(\w,m)}{e_\G}}(u,t,s)=e(t,s)$. 
Let $\theta$ be the canonical left invariant 1-form in $E$ (the Maurer-Cartan 1-form); see \ref{wkl}.
Taking the $t$ derivative of (\ref{E}), we obtain:
$$
\frac{\d}{\d t}\theta\left(\frac{\d}{\d s} e(t,s)\right)= m \left (\t{\frac{\d}{\d t}\g_s(t)},\t{\frac{\d}{\d s}\g_s(t)} \right) _{{\mathcal{H}}_\w(\g_s,t,u_s) }.
$$
By  (\ref{gtts}) and the $G$-equivariance of $m$:
$$
\d(e(t,s)) \tr \frac{\d}{\d t}\theta\left(\frac{\d}{\d s} e(t,s)\right) =
 m \left (\t{\frac{\d}{\d t}\g_s(t)},\t{\frac{\d}{\d s}\g_s(t)} \right) _{{\mathcal{H}}_\w(\g^t ,s,u^t)}.
$$
{We also have:}
\begin{align*}
&\frac{\d}{\d s} \left( \d(e(t,s)) \tr 
\theta\left(\frac{\d}{\d t} e(t,s)\right) \right) \\
&= \d(e(t,s)) \tr  \left(\d\left(\theta\left(\frac{\d}{\d s} e(t,s) \right) \right)
\tr \theta\left(\frac{\d}{\d t} e(t,s) \right)\right)
+ \d(e(t,s)) \tr \frac{\d}{\d s}\left(\theta\left(\frac{\d}{\d t} e(t,s)\right)\right)  \\
&= \d(e(t,s)) \tr \left( \left[ \theta\left(\frac{\d}{\d s} e(t,s) \right), \theta\left(\frac{\d}{\d t} e(t,s)\right) \right]
+ \frac{\d}{\d s}\theta\left(\frac{\d}{\d t} e(t,s)\right) \right) \\
&= \d(e(t,s)) \tr \frac{\d}{\d t}\theta\left(\frac{\d}{\d s} e(t,s)\right).
\end{align*}
The second equation follows from the definition of a differential crossed module, and the third from the fact $d\theta(X,Y)=-[X,Y]$ for each $X,Y \in \le$.
Combining the two equations and integrating in $s$, with ${\stackrel{(\w,m)}{e_\G}}(u,t,0)=1_E$, we obtain:
\begin{align*}
\frac{\d}{\d t} {\stackrel{(\w,m)}{e_\G}}(u,t,s)&=\left( \int_0^s m \left (\t{\frac{\d}{\d t}\g_{s'}(t)},\t{\frac{\d}{\d s'}\g_{s'}(t)} \right) _{{\mathcal{H}}_\w(\g^t ,s',u^t)} ds'\right){\stackrel{(\w,m)}{e_\G}}(u,t,s),
 \end{align*}
with initial condition ${\stackrel{(\w,m)}{e_\G}}(u,0,s)=1_E,$ (set $t=0$ in (\ref{E})), from which (\ref{Fubts}) follows as an immediate consequence.
\end{Proof}
Note that by using 
the {Non-Abelian Fubini's Theorem},  lemmas \ref{hm} and \ref{hi} follow directly from lemmas \ref{vm} and \ref{vi}. 

{From the {Non-Abelian Fubini's Theorem} and \ref{chpfb} it follows that the two-dimensional holonomy of a categorical connection is covariant with respect to the action of the dihedral group {$D_4\cong \Z_2^2 \rtimes \Z_2 $} of symmetries of the square;  see \ref{ddgm}.}

\subsection{Dependence of the categorical holonomy on a smooth family of squares}\label{dch}
{
In this subsection we prove a fundamental result giving the variation of the 2-holonomy of a smooth family of 2-paths in terms of the 2-curvature, analogous to equation (\ref{ddxg}) for the variation of the 1-holonomy of a smooth family of 1-paths in terms of the curvature. }
Let $P\to M$ be a principal $G$-bundle over the smooth manifold $M$ with a $\Gc$-categorical connection $(\w,m)$. Here  ${\mathcal G}=(E \ra{\d} G,\triangleright)$  
is a Lie crossed module, where $\tr$ is a Lie group left action of $G$ on $E$ by automorphisms. Let  $\mathfrak{G}=(\d\colon \le\to \lg,\tr)$ be the associated differential crossed module. 

{Consider a smooth map $J\colon [0,1]^3 \to M$. {Put $J(t,s,x)=\G^x(t,s),$ }where $x,t,s \in [0,1]$}. Define $q(x)=J(0,0,x)$, for each $x\in [0,1]$. Choose $u\in P_{q(0)}$   and let $u(x)=\H_\w(q,x,u)$. We want to analyse the dependence on $x$ of the categorical holonomy ${\stackrel{(\w,m)}{e_{\G^x}}}{(u(x),t,s)}$,  see equation (\ref{E}). To this end, we now prove the following {well known} technical lemma, also appearing in \cite{FMP}.

\subsubsection{A well-known lemma}\label{wkl}
Let $G$ be a Lie group. Consider a $\lg$-valued smooth function $V(s,x)$ defined on $[0,1]^2$.
Consider  the following differential equation in $G$:
 $$\frac{\d}{\d s} a(s,x)=a(s,x)V(s,x),$$ 
with initial condition $a(0,x)=1_G, \forall x \in [0,1]$.
We want to know $\frac{\d}{\d x} a(s,x)$.

Let $\theta$ be the canonical $\lg$-valued $1$-form on $G$. Thus $\theta$ is left {invariant} and satisfies $\theta(A)=A, \forall A \in \lg$, being defined uniquely by these properties. Also $
d \theta(A,B)=- \theta ([A,B]),$
where $A,B \in \lg$. 
We have:
$$\frac{\d }{ \d x} \theta\left (\frac{\d}{\d s} a(s,x) \right )= \frac{\d }{ \d x}\theta\big (  a(s,x)V(s,x) \big )= \frac{\d }{ \d x}V(s,x).
$$
On the other hand: 
\begin{align*}\frac{\d }{ \d x} \theta\left (\frac{\d}{\d s} a(s,x) \right ) 
&= d a^*(\theta)\left(  \frac{\d }{ \d x}, \frac{\d }{ \d s}\right) {+}\frac{\d }{ \d s} a^*(\theta) \left ( \frac{\d }{ \d x} \right){+} a^*(\theta)\left(  \left [\frac{\d }{ \d x}, \frac{\d }{ \d s} \right]\right) \\
&=d \theta\left(  \frac{\d }{ \d x}  a(s,x) , \frac{\d }{ \d s} a(s,x) \right)  {+}\frac{\d }{ \d s} \theta \left ( \frac{\d }{ \d x}a(s,x) \right).
\end{align*}  
Therefore:
$${\theta \left ( \frac{\d }{ \d x} a(s,x) \right)+ \int_0^s \big(d \theta\left(  \frac{\d }{ \d x}  a(s',x) , \frac{\d }{ \d s'} a(s',x) \right) {+}\frac{\d }{ \d x} V(s',x)\big)ds'=\theta\left (\frac{\d }{ \d x}a(0,x)\right).}$$
Since  $\frac{\d }{ \d x}a(0,x)=0$ (due to the initial conditions) we have the following:
\begin{Lemma}\label{WKL}
$$ \frac{\d }{ \d x} a(s,x)= a(s,x)  \int_0^s \left ({-}d \theta\left(  \frac{\d }{ \d x}  a(s',x) , \frac{\d }{ \d s'} a(s',x) \right) {+} \frac{\d }{ \d x} V(s',x)\right )ds',$$
{for each $x,s \in [0,1]$.}
\end{Lemma}
\subsubsection{The relation between 2-curvature and categorical holonomy}

The following main theorem is more general than the analogous result in \cite{FMP,SW2} since it is valid for any smooth homotopy $J$ connecting two 2-paths $\G$ and $\G'$, and in particular the  basepoints of the 2-paths may vary with the parameter $x$.  For this reason the proof is considerably longer, forcing several integrations by parts.

\begin{Theorem}\label{Main2}
Let $M$ be a smooth manifold. Let ${\Gc=(\d\colon E \to G,\tr)}$ be a Lie crossed module. Let $P \to M$ be a principal $G$-bundle over $M$. Consider a $\Gc$-categorical connection $(\w,m)$ on $P$.
Let $J\colon [0,1]^3 \to M$ be a smooth map. Let $J(t,s,x)=\G^x(t,s)=\g^x_s(t)=\g^{x,t}(s);\forall t,s,x \in [0,1]$. {Define $q(x)=\G^x(0,0)$. Choose $u \in P_{q(0)}$, the fibre of $P$ at $q(0)$. Let  $u(x)=\H_\w(q,{x},u)$ and $u(x,s)=\H_\w(\g^{x,0},s,u(x))$, where $s,x \in [0,1]$.}  

Consider the map $(s,x)\in [0,1]^2 \mapsto e_{\G^x}(s)\in E$ defined by:
\begin{equation}\label{E1}
{\frac{\d}{\d s} e_{\G^x}(s) =e_{\G^x}(s)\int_{0}^{1} m \left (\t{\frac{\d}{\d t}\g^x_s(t)},\t{\frac{\d}{\d s}\g^x_s(t)} \right)  _{{\mathcal{H}}_\w(\g_s^x,t,u(x,s)) } dt ,}
\end{equation}
with initial condition:
\begin{equation}\label{E2}
e_{\G^x}(0)=1_E, \forall x \in [0,1],
\end{equation}
{Let $e_{\G^x}=e_{\G^x}(1)$}. For each $x\in [0,1]$, we have:
\begin{align*}\frac{d }{ d x}e_{\G^x} &=e_{\G^x} \int_0^1 \int_{0}^{1} \M\left ( \t{\frac{\d}{\d x} \g_s^x(t)},\t{ \frac{\d}{\d t} \g_s^x(t)},\t{ \frac{\d}{\d s} \g_s^x(t)}\right)_{{\mathcal{H}}_\w(\g_s^x,t,u(x,s)) }  dtds\\
&\quad \quad 
+e_{\G^x}\int_0^1 m\left(\t{\frac{\d}{\d n} \hat{\g}^x(n)},\t{ \frac{\d}{\d x} \hat{\g}^x(n)}\right)_{{\mathcal{H}}_\w(\hat{\g}^x,n,u(x)) } dn,
\end{align*}
where $\hat{\g}^x=\d \G^x$, starting at $\G^x(0,0)$ and  oriented clockwise. {Here $\M\in \A^3(P,\le)$ is the categorical curvature $3$-form of $(\w,m)$; see \ref{curv}.}
\end{Theorem}
\begin{Proof}
{Consider the smooth} map $f\colon [0,1]^3 \to P$ such {that} $f(x,s,t)={{\mathcal{H}}_\w(\g_s^x,t,u(x,s)) }$, for each $ x,s,t \in [0,1].$ By definition we have:
$$\frac{\d}{\d t } f(x,s,t) = \t{ \frac{\d}{\d t}\g^x_s(t)}_{{\mathcal{H}}_\w(\g_s^x,t,u(x,s)) }$$ and therefore 
$\w(\frac{\d}{\d t } f(x,s,t))=0$. 
We also have: $$
\left (\frac{\d}{\d s } f(x,s,t)\right) ^H = \t{ \frac{\d}{\d s}\g^x_s(t)}_{{\mathcal{H}}_\w(\g_{s}^x,t,u(x,s)) }$$ and $$
\left ( \frac{\d}{\d x } f(x,s,t)\right)^H = \t{ \frac{\d}{\d x}\g^x_s(t)}_{{\mathcal{H}}_\w(\g_{s}^x,t,u(x,s)) }.$$
Note also that $m(X,Y)$, $\W(X,Y)$ and $\M(X,Y,Z)$  vanish if either $X$,  $Y$ or $Z$ is vertical. 

By the 2-structure equation, see Proposition \ref{cse}, and equation (\ref{BB}) it follows that (since $\M$ is horizontal):
\begin{align*}
\int_0^1 \int_{0}^{1} & \M\left ( \t{\frac{\d}{\d x} \g_s^x(t)},\t{ \frac{\d}{\d t} \g_s^x(t)},\t{ \frac{\d}{\d s} \g_s^x(t)}\right)_{{\mathcal{H}}_\w(\g_s^x,t,u(x,s)) }  dtds\\
&=\int_0^1 \int_{0}^{1} \M\left ( \frac{\d}{\d x} f{(x,s,t)}, \frac{\d}{\d t}f{(x,s,t)}, \frac{\d}{\d s} f{(x,s,t)}\right) dtds\\
&=\int_0^1 \int_{0}^{1} dm\left ( \frac{\d}{\d x} f{(x,s,t)}, \frac{\d}{\d t}f{(x,s,t)}, \frac{\d}{\d s} f{(x,s,t)}\right) dtds\\
& \quad \quad +\int_0^1 \int_{0}^{1} \w\left ( \frac{\d}{\d x} f{(x,s,t)}\right) \tr m\left (\frac{\d}{\d t}f{(x,s,t)}, \frac{\d}{\d s} f{(x,s,t)}\right) dtds\\
& \quad \quad -\int_0^1 \int_{0}^{1} \w\left ( \frac{\d}{\d s} f{(x,s,t)}\right) \tr m\left (\frac{\d}{\d t}f{(x,s,t)}, \frac{\d}{\d x} f{(x,s,t)}\right) dtds.
\end{align*}
{Using {{Lemma}} \ref{Main1} and integration by parts, we rewrite the integral in the last term:}
\begin{align*}
&\int_{0}^{1} \w\left ( \frac{\d}{\d s} f{(x,s,t)}\right) \tr m\left (\frac{\d}{\d t}f{(x,s,t)}, \frac{\d}{\d x} f{(x,s,t)}\right) dt \\
   &= \int_{0}^{1} \int_0^{t} \W  \left (\frac{\d}{\d t'}f{(x,s,t')}, \frac{\d}{\d s} f{(x,s,t')}\right){dt'} \tr m\left (\frac{\d}{\d t}f{(x,s,t)}, \frac{\d}{\d x} f{(x,s,t)}\right) dt \end{align*}
or
\begin{multline*}
\int_{0}^{1}   \W  \left (\frac{\d}{\d t'}f{(x,s,t')}, \frac{\d}{\d s} f{(x,s,t')}\right)dt' \tr \int_0^1 m\left (\frac{\d}{\d t'}f{(x,s,t')}, \frac{\d}{\d x} f{(x,s,t')}\right) dt' \\   - \int_0^1\W\left (\frac{\d}{\d t}f{(x,s,t)}, \frac{\d}{\d s} f{(x,s,t)}\right) \tr \left (\int_0^t m\left (\frac{\d}{\d t'}f{(x,s,t')}, \frac{\d}{\d x} f{(x,s,t')}\right) dt' \right )dt.
\end{multline*}
{Using equation (\ref{vu}), we have for the final term:}
\begin{align*}
&\int_0^1\W\left (\frac{\d}{\d t}f{(x,s,t)}, \frac{\d}{\d s} f{(x,s,t)}\right) \tr \left (\int_0^t m\left (\frac{\d}{\d t'}f{(x,s,t')}, \frac{\d}{\d x} f{(x,s,t')}\right) dt' \right )dt\\
 &=-\int_0^1 \int_0^t \W\left (\frac{\d}{\d t'}f{(x,s,t')}, \frac{\d}{\d x} f{(x,s,t')}\right) dt' \tr m\left (\frac{\d}{\d t}f{(x,s,t)}, \frac{\d}{\d s} f{(x,s,t)}\right)   dt\\
&=- \int_0^1 \w\left (\frac{\d}{\d x}f{(x,s,t)}\right)  \tr m\left (\frac{\d}{\d t}f{(x,s,t)}, \frac{\d}{\d s} f{(x,s,t)}\right)   dt \\
&\quad \quad + \int_0^1 \int_0^s \W\left (\frac{\d}{\d s'}f{(x,s',0)}, \frac{\d}{\d x} f{(x,s',0)}\right) ds' \tr m\left (\frac{\d}{\d t}f{(x,s,t)}, \frac{\d}{\d s} f{(x,s,t)}\right)   dt \\
&=- \int_0^1 \w\left (\frac{\d}{\d x}f{(x,s,t)}\right)  \tr m\left (\frac{\d}{\d t}f{(x,s,t)}, \frac{\d}{\d s} f{(x,s,t)}\right)   dt \\
&\quad \quad + \int_0^1 \w\left (\frac{\d}{\d x}f{(x,s,0)} \right) \tr m\left (\frac{\d}{\d t}f{(x,s,t)}, \frac{\d}{\d s} f{(x,s,t)}\right)   dt.
\end{align*}
{where we have used Lemma \ref{Main1} twice. Combining the previous equations, yields}
\begin{multline}\label{intM}
\int_0^1 \int_{0}^{1}  \M\left ( \t{\frac{\d}{\d x} \g_s^x(t)},\t{ \frac{\d}{\d t} \g_s^x(t)},\t{ \frac{\d}{\d s} \g_s^x(t)}\right)_{{\mathcal{H}}_\w(\g_s^x,t,u(x,s)) }  dtds\\
=\int_0^1 \int_{0}^{1} dm\left ( \frac{\d}{\d x} f{(x,s,t)}, \frac{\d}{\d t}f{(x,s,t)}, \frac{\d}{\d s} f{(x,s,t)}\right) dtds\\
 -\int_0^1 \int_{0}^{1}   \W  \left (\frac{\d}{\d t'}f{(x,s,t')}, \frac{\d}{\d s} f{(x,s,t')}\right)dt' \tr \int_0^1 m\left (\frac{\d}{\d t'}f{(x,s,t')}, \frac{\d}{\d x} f{(x,s,t')}\right) dt' ds\\
- \int_0^1 \int_0^1 \w\left (\frac{\d}{\d x}f{(x,s,0)} \right) \tr m\left (\frac{\d}{\d t}f{(x,s,t)}, \frac{\d}{\d s} f{(x,s,t)}\right)   dt ds.
\end{multline}

For the second term on the right hand side in the theorem, we obtain:
\begin{multline}\label{intdn}
\int_0^1  m\left(\t{\frac{\d}{\d n} \hat{\g}^x(n)},\t{ \frac{\d}{\d x}  \hat{\g}^x(n)}\right)_{{\mathcal{H}}_\w(\hat{\g}^x,n,u(x))} dn\\
=  \int_0^1 m\left ( \frac{\d}{\d s}f(x,s,0), \frac{\d}{\d x}f(x,s,0)\right )ds + \int_0^1 m\left ( \frac{\d}{\d t}f(x,1,t), \frac{\d}{\d x}f(x,1,t)\right )dt \\
-g^{-1} \tr \big(\int_0^1 m\left ( \frac{\d}{\d s}f'(x,s,1), \frac{\d}{\d x}f'(x,s,1)\right )ds-  \int_0^1 m\left ( \frac{\d}{\d t}f(x,0,t), \frac{\d}{\d x}f(x,0,t)\right )dt\big),
\end{multline}
{where we have put
$g(x,s)=\d(e_{\G^x}(s))$ and $f'(x,s,t)=\H_\w\left(\g^{x,t},s,\H_\w(\g^x_0,t,u(x))\right)$; also $g=g(x,1)$.} Therefore $f(x,s,1)=f'(x,s,1) \d(e_{\G^x}(s))$ by the {Non-Abelian Green's Theorem.} Note that $f'(x,0,t)=f(x,0,t)$.  We will be using the function $f'$ again shortly. 

\medskip
Thus it remains to prove that $e_{\G^x}^{-1} \frac{d }{ d x} e_{\G^x}$ is equal to the sum of the right hand sides of (\ref{intM}) and (\ref{intdn}).

\medskip
By {Lemma} \ref{WKL}, we have
\begin{equation}\label{torefer}
 \frac{d }{ d x} e_{\G^x} =e_{\G^x} (A_x-B_x),
\end{equation}
where
\begin{align*}
A_x&=\int_0^1 \int_{0}^{1} \frac{\d }{ \d x} \left ( m \left (\t{\frac{\d}{\d t}\g^x_s(t)},\t{\frac{\d}{\d s}\g^x_s(t)} \right)  _{{\mathcal{H}}_\w(\g_{s} ^x,t,u(x,s)) } \right) dt ds\\
B_x&=\int_0^1 d \theta\left(  \frac{\d }{ \d x}  e_{\G^x}(s) , \frac{\d }{ \d s} e_{\G^x}(s) \right)   ds.
\end{align*}
Let us analyse $A_x$ and $B_x$ separately. Using the well known equation:
$$
d\a(X,Y,Z)=X\a(Y,Z)+Y\a(Z,X)+Z\a(X,Y)+\a(X,[Y,Z])+ \a(Y,[Z,X])+\a(Z,[X,Y]),
$$
valid for any smooth 2-form {$\a$} in a manifold, and any three vector fields $X,Y,Z$ in $M$, we obtain for $A_x$:
\begin{align*}
A_x&= \int_{0}^{1}\int_{0}^{1} \frac{\d }{ \d x}  m \left (\frac{\d}{\d t} f(x,s,t),\frac{\d}{\d s}f(x,s,t) \right) dt ds\\ 
  &  =    \int_0^1 \int_{0}^{1} dm\left ( \frac{\d}{\d x} f(x,s,t), \frac{\d}{\d t}f(x,s,t), \frac{\d}{\d s}f(x,s,t)\right)dtds
\nonumber\\ 
&  {{-} \int_0^1 \int_{0}^{1} \frac{\d}{\d t} m\left ( \frac{\d}{\d s}f(x,s,t), \frac{\d}{\d x}f(x,s,t)\right )+ \frac{\d}{\d s} m\left ( \frac{\d}{\d x}f(x,s,t), \frac{\d}{\d t}f(x,s,t)\right )dtds}\nonumber\end{align*}
or
\begin{align}\label{Ax}
  A_x&= \int_{0}^{1}\int_{0}^{1} dm\left ( \frac{\d}{\d x} f(x,s,t), \frac{\d}{\d t}f(x,s,t), \frac{\d}{\d s}f(x,s,t)\right)dtds\\ 
&   +\int_0^1m\left ( \frac{\d}{\d s}f(x,s,0), \frac{\d}{\d x}f(x,s,0)\right )ds+\int_0^1 m\left ( \frac{\d}{\d t}f(x,1,t), \frac{\d}{\d x}f(x,1,t)\right )dt \nonumber\\
&   -\int_0^1 m\left ( \frac{\d}{\d t}f(x,0,t), \frac{\d}{\d x}f(x,0,t)\right )dt-\int_0^1m\left ( \frac{\d}{\d s}f(x,s,1), \frac{\d}{\d x}f(x,s,1)\right )ds\nonumber.
\end{align}
{Recall that $g(x,s)=\d(e_{\G^x}(s))$  and $f'(x,s,t)=\H_\w\left(\g^{x,t},s,\H_\w(\g^x_0,t,u(x))\right)$, and the relation $f(x,s,1)=f'(x,s,1) \d(e_{\G^x}(s))$. 
We thus have:}
$$
\w\left (\frac{\d}{\d s} f'(x,s,1)\right)=\w\left (\frac{\d}{\d s}\left( f(x,s,1)g^{-1}(x,s)\right)\right),
$$
which since $\frac{\d}{\d s} f'(x,s,1)$ is horizontal implies, by using the Leibniz rule and the fact that $\w$ is a connection $1$-form, that:
$$g(x,s)\w\left (\frac{\d}{\d s} f(x,s,1)\right)g^{-1}(x,s)+g(x,s) \frac{\d}{\d s }g^{-1}(x,s)=0. $$
Analogously (this will be used later):
\begin{multline*}g^{-1}(x,s)\w\left(\frac{\d}{\d x}\left(\H_\w (\g^{x,1},s, u(x)\g^x_0)\right)\right)g(x,s)\\=-g^{-1}(x,s)\frac{\d}{\d x}g(x,s)+\w\left(\frac{\d}{\d x} \big ( u(x,s)\g_s^x \big)\right),\end{multline*}
which is the same as:
$$\frac{\d}{\d x} g(x,s)=g(x,s) \w\left(\frac{\d}{\d x} f(x,s,1)\right)-\w\left(\frac{\d}{\d x} f'(x,s,1)\right)g(x,s).
$$
The very last term {$R$} of (\ref{Ax}) can be simplified as follows (since $m$ is horizontal and $G$-equivariant):
\begin{align*}
R&=-\int_0^1m\left ( \frac{\d}{\d s}f(x,s,1), \frac{\d}{\d x}f(x,s,1)\right )ds\nonumber\\
&=-\int_0^1 g^{-1}(x,s)\tr m\left ( \frac{\d}{\d s}f'(x,s,1), \frac{\d}{\d x}f'(x,s,1)\right )ds\nonumber\\
&=-g^{-1}(x,1) \tr \int_0^1  m\left ( \frac{\d}{\d s}f'(x,s,1), \frac{\d}{\d x}f'(x,s,1)\right )ds \nonumber\\ 
& \quad \quad +\int_0^1 \int_0^s\frac{\d}{\d s}g^{-1}(x,s)\tr m\left ( \frac{\d}{\d s'}f'(x,s',1), \frac{\d}{\d x}f'(x,s',1)\right )ds'ds;\nonumber
\end{align*}
{(the penultimate equation follows from integrating by parts). Therefore:} 
\begin{multline}\label{Ax2}
R=-g^{-1}(x,1) \tr \int_0^1  m\left ( \frac{\d}{\d s}f'(x,s,1), \frac{\d}{\d x}f'(x,s,1)\right )ds \\
  \quad \quad -\int_0^1 \int_0^s \w \left (\frac{\d}{\d s}f{(x,s,1)}  \right) g^{-1}(x,s)  \tr m \left ({\frac{\d}{\d s'}f'(x,s',1)},{\frac{\d}{\d x}f(x,s',1)} \right)ds'ds.
\end{multline}

{We now analyse $B_x$, for each $x \in [0,1]$. We have:}
\begin{align*}
B_x &=d \theta\Big( e_{\G^x}^{-1}(s) \frac{\d }{ \d x}  e_{\G^x}(s) ,e_{\G^x}^{-1}(s) \frac{\d }{ \d s} e_{\G^x}(s) \Big) \\
&=- \Big [ e_{\G^x}^{-1}(s) \frac{\d }{ \d x}  e_{\G^x}(s) ,e_{\G^x}^{-1}(s) \frac{\d }{ \d s} e_{\G^x}(s) \Big] \\
&= -\left ( g^{-1}(x,s)\frac{\d }{ \d x} g(x,s)\right) \tr \left (e_{\G^x}^{-1}(s) \frac{\d }{ \d s} e_{\G^x}(s)\right)\\
&=-\w\left(\frac{\d}{\d x} f(x,s,1)\right)  \tr \int_0^1 m \left ({\frac{\d}{\d t}f(x,s,t)},{\frac{\d}{\d s}f(x,s,t)}  \right)    dt\\&\quad+\left(g^{-1}(x,s)\w\left(\frac{\d}{\d x} f'(x,s,1)\right) g(x,s)\right) \tr \int_0^1 m \left ({\frac{\d}{\d t}f(x,s,t)},{\frac{\d}{\d s}f(x,s,t)}  \right)dt.  
\end{align*}
By using {Lemma} \ref{Main1}, this may be rewritten as $B_x=C_x+C_x'$, where
\begin{multline*}
C_x =-\int_0^1 \W \left ({\frac{\d}{\d t}f(x,s,t)},{\frac{\d}{\d x}f(x,s,t)}  \right)    dt\tr \int_0^1 m \left ({\frac{\d}{\d t}f(x,s,t)},{\frac{\d}{\d s}f(x,s,t)}  \right)    dt\\
 -\int_0^s \W \left ({\frac{\d}{\d s'}f(x,s',0)},{\frac{\d}{\d x}f(x,s',0}  \right)    ds'\tr \int_0^1 m \left ({\frac{\d}{\d t}f(x,s,t)},{\frac{\d}{\d s}f(x,s,t)}  \right)    dt 
\end{multline*}
and
\begin{multline*}
C_x'=\left(g^{-1}(x,s) \left (\int_0^s \W \left ({\frac{\d}{\d s'}f'(x,s',1)},{\frac{\d}{\d x}f'(x,s',1)}\right) ds'\right)
 g(x,s)\right) \tr\\ \int_0^1 m \left ({\frac{\d}{\d t}f(x,s,t)},{\frac{\d}{\d s}f(x,s,t)}  \right)  dt
\\+\left(g^{-1}(x,s)\left( \int_0^1 \W \left ({\frac{\d}{\d t}f(x,0,t)},{\frac{\d}{\d x}f(x,0,t)}\right)dt\right) 
 g(x,s)\right) \tr \\\int_0^1 m \left ({\frac{\d}{\d t}f(x,s,t)},{\frac{\d}{\d s}f(x,s,t)}  \right)  dt,
\end{multline*}
{Again using $\d(m)=\W$ and $\d(u)\tr v=[u,v]= -[v,u]= -\d(v) \tr u; \forall u,v \in \le$,} together with $\d(m)=\W$ and {Lemma} \ref{Main1}, for all but the second term of the right hand side of the previous equation, we obtain:
\begin{multline}\label{tcompare}
{B_x} =\\
\int_0^1 \left (\int_{0}^{1} \W \left (\frac{\d}{\d t}f{(x,s,t)},\frac{\d}{\d s}f{(x,s,t)}  \right ) dt  \tr  \int_0^1 m \left (\frac{\d}{\d t}f{(x,s,t)},\frac{\d}{\d x}f{(x,s,t)}  \right)   dt\right)ds. \\
-\int_0^1 \int_0^1 \w \left ({\frac{\d}{\d x}f(x,s,0)}   \right) \tr m \left ({\frac{\d}{\d t}f(x,s,t)},{\frac{\d}{\d s}f(x,s,t)}  \right)    dtds
\\ -\int_0^1 \int_0^s \w \left (\frac{\d}{\d s}f{(x,s,1)}  \right) g^{-1}(x,s)  \tr m \left ({\frac{\d}{\d s'}f'(x,s',1)},{\frac{\d}{\d x}f'(x,s',1)} \right)ds'ds    
\\-\int_0^1 \int_0^1 \w \left (\frac{\d}{\d s}f{(x,s,1)}  \right) g^{-1}(x,s)  \tr m \left ({\frac{\d}{\d t}f(x,0,t)},{\frac{\d}{\d x}f(x,0,t)} \right)dtds  
\end{multline}
 Finally, since (given that $\w$ is a connection 1-form):
\begin{align*}
\w \left (\frac{\d}{\d s}f{(x,s,1)}  \right)g^{-1}(x,s)&=g^{-1}(x,s) \w \left (\frac{\d}{\d s}f{(x,s,1)} g^{-1}(x,s)  \right)\\
&=g^{-1}(x,s) \w \left (\frac{\d}{\d s}f'{(x,s,1)} - f(x,s,1) \frac{\d }{\d s}g^{-1}(x,s)  \right)\\
&=-g^{-1}(x,s) \w \left (  f(x,s,1) \frac{\d }{\d s}g^{-1}(x,s)  \right)\\
&=-\frac{\d }{\d s}g^{-1}(x,s),
\end{align*}
 the last term {$R'$} of the previous expression is rewritten as follows:
\begin{align*}
R'&=-\int_0^1 \int_0^1 \w \left (\frac{\d}{\d s}f{(x,s,1)}  \right) g^{-1}(x,s)  \tr m \left ({\frac{\d}{\d t}f(x,0,t)},{\frac{\d}{\d x}f(x,0,t)} \right)dsdt\nonumber\\
&=\int_0^1 \int_0^1 \frac{\d}{\d s} g^{-1}(x,s)  \tr m \left ({\frac{\d}{\d t}f(x,0,t)},{\frac{\d}{\d x}f(x,0,t)} \right)dsdt,
\end{align*} or
\begin{equation}\label{Bx2}
R'= g^{-1}(x) \tr \int_0^1    m \left ({\frac{\d}{\d t}f(x,0,t)},{\frac{\d}{\d x}f(x,0,t)} \right)dt-  \int_0^1    m \left ({\frac{\d}{\d t}f(x,0,t)},{\frac{\d}{\d x}f(x,0,t)} \right)dt,
\end{equation}
{where we have put $g(x)=g(x,1)$.}
Combining $A_x-B_x$ from equations (\ref{Ax}), (\ref{Ax2}), (\ref{tcompare}), (\ref{Bx2}), four terms cancel and the remaining terms are equal to the sum of the right hand sides of (\ref{intM}) and (\ref{intdn}). {This finishes the proof of Theorem \ref{Main2}.}
\end{Proof}

\subsubsection{Invariance under thin homotopy}
From  {Theorem} \ref{Main2} and the fact that the horizontal lift $X \mapsto \t{X}$ of vector fields on $M$ defines a linear map $\X(M) \to \X(P)$ we obtain the following:
\begin{Corollary}\label{cormain2}
{Let $M$ be a smooth manifold. Let also ${\Gc=(\d\colon E \to G,\tr)}$ be a Lie crossed module. Let $P \to M$ be a principal $G$-bundle over $M$, and consider a $\Gc$-categorical connection $(\w,m)$ on $P$. If $\G$ and $\G'$ are { rank-2 homotopic (see {Definition} \ref{thin2}) {2-paths} } $[0,1]^2 \to M$ then $\stackrel{(\w,m)}{e_\G}(u,t,s)=\stackrel{(\w,m)}{e_{\G'}}(u,t,s)$, whenever $u\in P_{\G(0,0)}$, the fibre of $P$ at $\G(0,0)=\G'(0,0)$, and for each $t,s \in [0,1]$.}
\end{Corollary}

\subsubsection{A (dihedral) double groupoid map}\label{ddgm}
Let $P$ be a principal $G$ bundle over $M$. We define a double groupoid $\D^2(P)$ whose set of objects is $M$, and whose set of morphisms $x \to y$ is given by all {right} $G$-equivariant maps $a\colon P_x \to P_y$. A 2-morphism is given by a square of the form:
\begin{equation} \begin{CD} &P_z @>d>> &P_w\\
              &@A cAA \hskip-1.3cm \scriptstyle{f}  &@AA bA\\
              &P_x @>> a> &P_y
  \end{CD}\end{equation}
where $x,y,z,w \in M$ and $a,b,c,d$ are right $G$-equivariant maps. Finally $f\colon P_x \to E$ is a smooth map such that $f(ug)=g^{-1} \tr f(u)$ for each $u\in P_x$ and $g \in G$, satisfying $(b\circ a)(u)\d(f(u))=(d\circ c)(u)$, for each $u \in P_x$. The horizontal  and vertical compositions are as in \ref{chpfb}.  We also have an action of the dihedral group {$D_4 \cong \Z_2^2 \rtimes \Z_2$ of the {2-cube}} given by the horizontal and vertical reversions, and such that the interchange of coordinates is accomplished by the move $f\mapsto f^{-1}$.
As a corollary of the discussion in the last two subsections it follows:
\begin{Theorem}\label{fubd}
Whenever the principal $G$-bundle $P\to M$ is equipped with a categorical connection $(\w,m)$, the holonomy and categorical holonomy maps ${\mathcal H}_\w$ and $\stackrel{(\w,m)}{e}$ define a double groupoid morphism $\stackrel{(\w,m)}{\mathcal H}\colon \S_2(M) \to \D^2(P)$, where $\S_2(M)$ is the thin fundamental double groupoid of $M$. Given a dihedral group element $r\in D_4$ we have 
$$\stackrel{(\w,m)}{\mathcal H}(\G\circ r^{-1})=r\left (\stackrel{(\w,m)}{\mathcal H}(\G)\right). $$
\end{Theorem}

\section{Cubical $\Gc$-2-bundles with connection}\label{simport}

\subsection{Definition of a cubical $\Gc$-2-bundle}
{Recall the conventions introduced in \ref{cs} and \ref{estg}.}

{Let $M$ be a smooth manifold.}
Let $\U=\{U_i\}_{i \in \I}$ be an open cover of $M$. {From this we  can define a cubical set $C(M,\U)$.}  For each positive integer $n$ the set $C^n(M,\U)$ of $n$-cubes of $C(M,\U)$ is given by all pairs $(x,R)$, where $R$  is an assignment of an element $U_v^R \in \U$ to each vertex of $v$ of $D^n$, such that the intersection $$U^R=\bigcap_{\textrm{vertices }v \textrm{ of } D^n} U_v^R $$
is non-empty, and $x \in U^R.$
The face maps $\d^{\pm}_i\colon C^n(M,\U) \to  C^{n-1}(M,\U)$ where $i \in \{1,\ldots, n\}$ and {$n=1,2,\ldots$,} are defined by 
$$\d^\pm_i(x,R)=(x, {R \circ \delta^\pm_i}).$$ Analogously, the degeneracies are given by:
$$\e_i(x,R)=(x, {R \circ \sigma_i}).$$

{The cubical set $C(M,\U)$ is clearly a cubical object in the category of manifolds, in other words a cubical manifold.}
{Given an $x \in M$, the cubical set $C(M,\U,x)$ is given by all the cubes of $C(M,\U)$ whose associated element of $M$  is $x$.}

\begin{Definition}[Cubical $\Gc$-2-bundle]\label{cgerbe}
{Let $\Gc=(\d\colon E \to G, \tr)$ be a Lie crossed module. Let $\Nc(\Gc)$ be the cubical nerve of $\Gc$; see \cite{BHS} and \ref{estg}, which is a cubical manifold. Let $M$ be a smooth manifold and  $\U=\{U_i\}_{i \in \I}$ be an open cover of $M$.
A cubical $\Gc$-2-bundle over $(M,\U)$ is given by a  map $C(M,\U) \to \Nc(\Gc)$ of cubical manifolds.}

Unpacking this definition, we see that a cubical $\Gc$-2-bundle is specified by smooth maps $\f_{ij} \colon U_i \cap U_j \to G$, where $U_i,U_j \in \U$ have a non-empty intersection, and also by smooth maps $\p_{ijkl}\colon  U_i \cap U_j  \cap U_k \cap U_l \to  E$, where $U_i,U_j,U_k,U_l \in \U$ have a non-empty intersection, such that:

 \begin{enumerate}
 \item {We have $\d(\p_{ijkl})^{-1}\f_{ij}\f_{jl}=\f_{ik}\f_{kl}$ in $U_{ijkl} \doteq U_i \cap U_j  \cap U_k \cap U_l$. In other words, putting $\f_{ij}=X^-_2({\cm_2}),  \quad \f_{ik}=X^-_1({\cm_2}),\quad \f_{kl}=X^+_2({\cm_2}),\quad  \f_{jl}=X^+_1({\cm_2})$ and $e({\cm_2})=\p_{ijkl}$ yields a flat $\Gc$-colouring $\cm_2=(\psi,\phi)_{ijkl}$ of $D^2$, for each $x \in U_{ijkl}$.}

\item {Given $i^\pm,j^\pm,k^\pm,l^\pm\in \I$ with $U_{i^-j^-k^-l^-} \cap U_{i^+j^+k^+l^+} \neq \emptyset$, and putting $$e^\pm_3({\cm_3})=(\psi,\phi)_{i^\pm j^\pm k^\pm l^\pm}, \quad e^-_1({\cm_3})=(\psi,\phi)_{i^-k^-i^+k^+}, \quad e^+_1({\cm_3})=(\psi,\phi)_{j^-l^-j^+l^+},$$ $$e^-_2({\cm_3})=(\psi,\phi)_{i^-j^-i^+j^+} \textrm{ and } e^+_2({\cm_3})=(\psi,\phi)_{k^-l^-k^+l^+}$$ yields a flat $\Gc$-colouring $\cm_3$ of $D^3$ in $U_{i^-j^-k^-l^-} \cap U_{i^+j^+k^+l^+}$.}
\item $\f_{ii}=1_G$ in $U_i$ for all $i \in \I$.
\item $\p_{iijj}=\p_{ijij}=1_E$ in $U_{ij}$
\end{enumerate}
See {Figure \ref{labels}} for our conventions in labelling the vertices of $D^2$ and $D^3$.
\end{Definition}
\begin{figure}
\centerline{\relabelbox 
\epsfxsize 12cm
\epsfbox{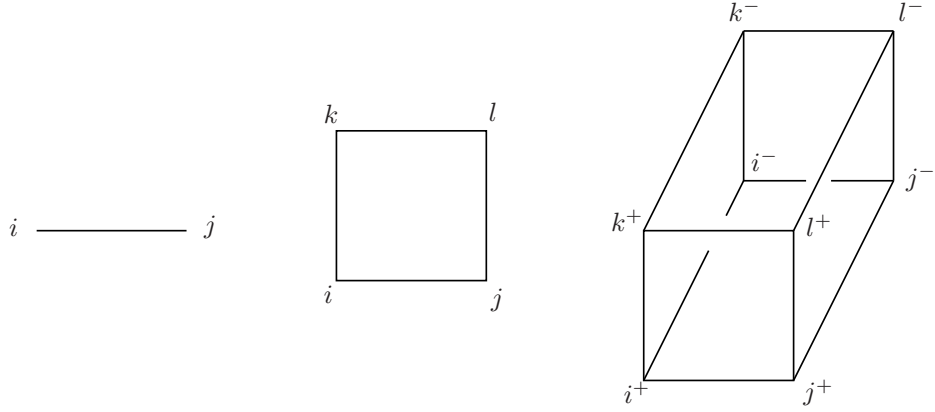}
\relabel{a}{${i}$}
\relabel{b}{${j}$}
\relabel{c}{${i}$}
\relabel{d}{${j}$}
\relabel{e}{${k}$}
\relabel{f}{${l}$}
\relabel{g}{${i^+}$}
\relabel{h}{${j^+}$}
\relabel{i}{${k^+}$}
\relabel{j}{${l^+}$}
\relabel{k}{${i^-}$}
\relabel{l}{${j^-}$}
\relabel{m}{${k^-}$}
\relabel{n}{${l^-}$}
\endrelabelbox}
\caption{\label{labels} Label conventions in {Definition} \ref{cgerbe}.}
\end{figure}
The previous definition is therefore a cubical {counterpart of the  simplicial definition} of a $\Gc$-2-bundle (and non-abelian gerbe) appearing for example in \cite{BrMe,ACG,BS1,BS2,SW3}. 
\begin{Remark}Note that in {Definition} \ref{cgerbe} the word bundle is used in the same sense as when one defines a principal bundle in terms of its transition functions, without reference to a total space; we are following {\cite{H,MP,P}}. For a discussion of the concept of total space of a non-abelian gerbe, see \cite{RS,Bar,Wo}.
\end{Remark}

\begin{Definition}[\bf  Dihedral cubical $\Gc$-2-bundles]\label{Dih}
{Recall that the cubical sets $C(M,\U)$ and $\Nc(\Gc)$  are dihedral; see \ref{cs}.
Therefore we can restrict our definition of a cubical $\Gc$-2-bundle and only allow dihedral cubical maps  $C(M,\U) \to \Nc(\Gc)$  which gives the definition of a  dihedral cubical $\Gc$-2-bundle.}
Explicitly, a cubical $\Gc$-2-bundle is said to be dihedral if the maps $\f_{ij}\colon U_{ij} \to G$ and $\psi_{ijkl}\colon U_{ijkl} \to E$ satisfy the following extra conditions:
\begin{enumerate}
 \item We have $\f_{ji}=\f_{ij}^{-1}$ in $U_{ij}$ for all $i,j \in \I$.
\item We have $ \psi_{ikjl}=\psi_{ijkl}^{-1}$, {$\psi_{jilk}=\f_{ij}\tr \psi_{ijkl}^{-1}$ and $\psi_{klij}=\f_{ik}\tr \psi_{ijkl}^{-1}$} in $U_{ijkl}$.
\end{enumerate}

\end{Definition}

\subsection{Connections in cubical $\Gc$-2-bundles}
Let  ${\mathcal G}=(\d\colon E\to  G,\triangleright)$  
be a Lie crossed module, where $\tr$ is a Lie group left action of $G$ on $E$ by automorphisms. Let also  $\mathfrak{G}=(\d\colon \le\to \lg,\tr)$ be the associated differential crossed module.

\begin{Definition}[Connection in a cubical $\Gc$-2-bundle]\label{ccg2b}
Let $M$ be a smooth manifold with an open cover $\U=\{U_i\}_{i \in \I}$. A connection in a cubical $\Gc$-2-bundle over $(M,\U)$ is given by:
\begin{itemize}
\item For any $i\in \I$ a local connection pair $(A_i,B_i)$ {defined in $U_i$}; in other words $A_i \in \A^1(U_i,\lg)$, $B_i \in \A^2(U_i,\le)$ and $\d(B_i)=d A_i +\frac{1}{2}A_i \wedge^\ad A_i=\W_{A_i}.$
\item For any ordered pair $(i,j)$ an $\le$-valued 1-form $\eta_{ij}$ in $U_{ij}$.
\end{itemize}
The conditions that should hold are:

\begin{enumerate}
 
\item For any $i \in \I$ we have $\eta_{ii}=0$.
\item For any $i,j \in \I$ we have:
$$A_j=\phi_{ij}^{-1} \left(A_i+\d (\eta_{ij})\right)\phi_{ij} +\phi_{ij}^{-1} d\phi_{ij},$$
$$B_j=\phi_{ij}^{-1} \tr \left(B_i +d\eta_{ij}+\frac{1}{2} \eta_{ij} \wedge^\ad \eta_{ij}+ A_i \wedge^\tr  \eta_{ij}\right).$$

\item {For any $i,j,k,l \in  \I$ we} have:\begin{multline*}\eta_{ik}+{\phi}_{ik} \tr \eta_{kl}  -{\phi}_{ik} {\phi}_{kl} {\phi}_{jl}^{-1}\tr \eta_{jl}-{\phi}_{ik} {\phi}_{kl} {\phi}_{jl}^{-1}{\phi}_{ij}^{-1} \tr \eta_{ij}\\
=\psi_{ijkl}^{-1}d\psi_{ijkl} +\psi_{ijkl}^{-1}\left(A_i\wedge^\tr \psi_{ijkl}\right).\end{multline*}
\end{enumerate}

\end{Definition}

\noindent
{The equivalence of cubical $\Gc$-2-bundles with connection will be dealt with in subsection \ref{Subdiv}.}

\begin{Definition}[Dihedral connection]\label{ccdg2b}
If a cubical $\Gc$-2-bundle is dihedral, then a connection in it is said to be dihedral if the following extra condition holds:
$$\eta_{ji}=-\f_{ij}^{-1} \tr \eta_{ij}, \textrm{ for each } i,j \in \I;$$
therefore, condition $3$ of the previous definition can be written as:
$$\eta_{ik}+{\phi}_{ik} \tr \eta_{kl}  +{\phi}_{ik} {\phi}_{kl} \tr \eta_{lj}+{\phi}_{ik} {\phi}_{kl} {\phi}_{lj} \tr \eta_{ji}
{=\psi_{ijkl}^{-1}d\psi_{ijkl} +\psi_{ijkl}^{-1}\left(A_i\wedge^\tr \psi_{ijkl}\right).}
 $$

\end{Definition}

\section{Non-abelian integral calculus based on a crossed module}

\subsection{{Path-ordered exponential and surface-ordered exponential}}\label{pese}
We {continue with} the notation and results of subsections \ref{hln} and \ref{dch}. Alternative direct derivations of some of the following results appear in \cite{BS1,SW1,SW2,SW3}.

Let $M$ be a manifold, and let $G$ be a Lie group with Lie algebra $\lg$. Let $\g\colon [0,1] \to M$ be a piecewise smooth map. Let $A \in \A^1(M,\lg)$ be a $\lg$-valued 1-form in $M$. We define, {as is usual}, the path ordered exponential $ {\stackrel{A}{g}}\g(t)=\pexp\left(\int_0^t A\left( \frac{d}{d t'} \g(t')\right) dt'\right)$ {to be} the solution of the differential equation:
$$ \frac{d}{dt} {\stackrel{A}{g}}\g(t)={\stackrel{A}{g}}\g(t) A\left( \frac{d}{d t} \g(t)\right) ,$$
with initial condition  ${\stackrel{A}{g}}\g(0)=1_G$; see \cite{Ch}. Put ${\stackrel{A}{g}}\g\doteq {\stackrel{A}{g}}\g(1)=\pexp\left(\int_0^1 A\left( \frac{d}{d t} \g(t)\right) dt\right)$.
 We immediately get that ${\stackrel{A}{g}}_{\g\g'}={\stackrel{A}{g}}_\g{\stackrel{A}{g}}_{\g'}$, and also ${\stackrel{A}{g}}_{\g^{-1}}=({\stackrel{A}{g}}_\g)^{-1}$. Here $\g$ and $\g'$ are piecewise smooth maps with $\g(1)=\g'(0)$.

Consider the trivial bundle $P=M \times G$ over $M$. Given $A \in \A^1(M,\lg)$ there exists a unique connection $1$-form $\w_A$ in the trivial bundle $P$ for which $A=\zeta^*(\w_A)$, where $\zeta(x)=(x,1_G)$ for each $x \in M$. We then  have that:
$$\zeta(\g(t))=\H_{\w_A}(\g,t,\zeta(\g(0)))\pexp\left(\int_0^t A\left( \frac{d}{d t'} \g(t')\right) dt'\right).$$

Let  ${\mathcal G}=(\d\colon E\to  G,\triangleright)$  
be a Lie crossed module and let  $\mathfrak{G}=(\d\colon \le\to \lg,\tr)$ be the associated differential crossed module.
As before, if we have $B \in  \A^2(M,\le)$ with $\d(B)=\W_A=d A+\frac{1}{2} A \wedge^\ad A$ we define 
$$\stackrel{(A,B)}{e_\G}(t,s)=\sexp \left(\int_0^s \int_0^t B\left( \frac{\d}{\d t'} \g_{s'}(t'), \frac{\d}{\d s'} \g_{s'}(t')\right) dt' ds'\right)$$
as being the solution of the differential equation:
$$\frac{\d}{\d s} \stackrel{(A,B)}{e_\G}(t,s)=\stackrel{(A,B)}{e_\G}(t,s)  \int_0^t \left(\stackrel{A}{g}_{\g^0(s)}\stackrel{A}{g}_{\g_{s}(t')}\right) \tr B \left(\frac{\d}{\d t'}\g_{s}(t'),\frac{\d}{\d 
s}\g_{s}(t')\right)dt' $$ 
with initial conditions $$ \stackrel{(A,B)}{e_\G}(t,0)=1_E,\forall t \in [0,1].$$
{Put $\stackrel{(A,B)}{e_\G}=\stackrel{(A,B)}{e_\G}(1,1)$.}
We can equivalently define the surface ordered exponential by the differential equation:
$$\frac{\d}{\d t} \stackrel{(A,B)}{e_\G}(t,s)= \left( \int_0^s \left(\stackrel{A}{g}_{\g_0(t)}\stackrel{A}{g}_{\g^t(s')}\right) \tr B \left(\frac{\d}{\d t}\g_{s'}(t),\frac{\d}{\d 
s'}\g_{s'}(t)\right)ds'\right) \stackrel{(A,B)}{e_\G}(t,s) $$ 
with initial conditions $$  {\stackrel{(A,B)}{e_\G}(0,s)}=1_E,\forall s \in [0,1];$$
{see the proof of Theorem \ref{naftbf} and below.}

As before, there exists a unique categorical connection  $(\w_A,m_{A,B})$ in the trivial bundle $P=M \times G$ for which $A=\zeta^*(\w_A)$ and $B=\zeta^*(m_{A,B}))$. We have that {$\stackrel{(A,B)}{e_\G}(t,s)={\stackrel{(\w_A,m_{A,B})}{e_\G}}({\zeta(\G(0,0))},t,s)$,} see \ref{naft}. The following follows immediately from the {Non-Abelian Green's  Theorem} \ref{gtb}.

\begin{Theorem}[Non-abelian Green's Theorem, elementary form]
{Consider a 2-square $\G\colon [0,1]^2 \to M$.}
{Put $\stackrel{A}{X_\G}=\stackrel{A}{g}_{\Xc_\G}$, $\stackrel{A}{Y_\G}=\stackrel{A}{g}_{\Yc_\G}$,  $\stackrel{A}{Z_\G}=\stackrel{A}{g}_{\Zc_\G}$ and $\stackrel{A}{W_\G}=\stackrel{A}{g}_{ \Wc_\G}$;} see \ref{chpfb} for this notation. We have that:
$$\d\Big(\stackrel{(A,B)}{e_\G}\Big)^{-1}\stackrel{A}{X_\G} \stackrel{A}{Y_\G}=\stackrel{A}{Z_\G}  \stackrel{A}{W_\G}.$$
\end{Theorem}
The following follows from theorems \ref{vm} and \ref{hm}. See \ref{cs} and subsection \ref{thin2g}.
\begin{Theorem}
Consider the map $\stackrel{(A,B)}{\Hc}\colon  C^2(M) \to \D^2(\Gc)$ such that:
$$\stackrel{(A,B)}{\Hc}(\G)=\quad \quad \begin{CD} &* @>\stackrel{A}{W_\G}>> &*\\
              &@A \stackrel{A}{Z_\G} AA \hskip-1.3cm \scriptstyle{\stackrel{(A,B)}{e_\G}}  &@AA \stackrel{A}{Y_\G} A\\
              &* @>> \stackrel{A}{X_\G}> &*
  \end{CD} $$
{Then $\stackrel{(A,B)}{\Hc}(\G \circ_\h \G')=\stackrel{(A,B)}{\Hc}(\G)\circ_\h \stackrel{(A,B)}{\Hc}(\G')$ and $\stackrel{(A,B)}{\Hc}(\G \circ_\vm \G')=\stackrel{(A,B)}{\Hc}(\G)\circ_\vm  \stackrel{(A,B)}{\Hc}(\G')$, whenever the compositions of $\G,\G' \colon [0,1]^2 \to M$ are well defined. }
\label{localH}
\end{Theorem}

{Passing to the quotient $\S_2(M)$ of $C^2_r(M)$ under thin homotopy it follows, by using {Theorem} \ref{Main2} and {Corollary} \ref{cormain2},} that:
\begin{Theorem}\label{imp}
{The map  $\stackrel{(A,B)}{\Hc}\colon  {\Sc_2(M)} \to \D^2(\Gc)$ defined in the previous theorem is a morphism of double groupoids with thin structure.}
\end{Theorem}

The following {result is a consequence of} {Theorem} \ref{fubd}.
\begin{Theorem}
[Non-abelian Fubini's Theorem]\label{NAFT}{The   map $\stackrel{(A,B)}{\Hc}\colon C^2(M) \to \D^2(\Gc)$} preserves the action of the dihedral group $D_4$ of the square. Concretely for any element $r$ of $D_4$ we have
$$ \stackrel{(A,B)}{\Hc}(\G \circ r^{-1})=r(\stackrel{(A,B)}{\Hc}(\G )),$$
for each smooth map $\G\colon [0,1]^2 \to M$.
 \end{Theorem}
This {follows from} the fact that $\stackrel{(A,B)}{\Hc}$ preserves horizontal and vertical reversions and moreover interchanges of coordinates, which generate the dihedral group $D_4\cong \Z_2^2\rtimes S_2$ of the square.

We finish this subsection with the following important theorem:
\begin{Theorem}\label{Main3}
Let $(A,B)$ be a local connection pair	 in $M$, by which as usual we mean $A \in \A^1(M,\lg)$, $B \in \A^2(M,\le)$ and $\d(B)=\W_A=d A+\frac{1}{2}A \wedge^\ad A.$ Let $C=d B+A \wedge^\tr B$ be the 2-curvature 3-form of $(A,B)$ as in \ref{curv} and \ref{lf}. Let $J \colon [0,1]^3 \to M$ be a smooth map such that $J^*(C)=0$. Then the colouring  $T$ of $D^3$ such that:
$$
T \circ \delta^\pm_i=\stackrel{(A,B)}{\Hc}(\d^\pm_i J),\quad {i=1,2,3}
$$
is flat;  see \ref{estg} and \ref{cs}.
\end{Theorem}
\begin{Proof}
This follows from the construction in this subsection and {Theorem} \ref{Main2}. Note the form (\ref{ha2}) for the homotopy addition equation (\ref{ha}).
\end{Proof}

\subsection{1-Gauge transformations}\label{Trans}

Let $M$ be a smooth manifold.  Let $(A,B)$ and $(A',B')$ be local connection pairs defined in $M$. {For the time being we will drop the index $i$ for the open cover and take $A$ and $B$ to be globally defined on $M$. We will return to the general case in the next section.} In other words $A,A' \in \A^1(M,\lg)$ and $B,B'\in \A^2(M,\le)$ are such that $\d(B)=\W_A=dA+\frac{1}{2}A \wedge^\ad A$ and $\d(B')=\W_{A'}$. Let $\eta\in \A^1(M,\le)$ be such that:
$$A'=A +\d(\eta) $$
and
$$B'=B+d\eta+\frac{1}{2}\eta  \wedge^\ad \eta+A \wedge^\tr \eta .$$

Given a smooth path $\g\colon [0,1] \to M$, define the following 2-square in $\Gc$:
$${\tau_{A}^{(1_G,\eta)}(\g)\quad=\quad \quad \begin{CD} &* @>\stackrel{A'}{g_\g} >> &*\\
              &@A 1_G AA \hskip-1.5cm \scriptstyle{\stackrel{(A,\eta)} {f_{\g }}} &@AA 1_GA\\
&* @>>\stackrel{A}{g_\g} > &*
 \end{CD}\quad \quad \doteq\quad \quad \begin{CD} &* @>\stackrel{A'}{g_\g} >> &*\\
             &@A 1_G AA \hskip-1.3cm \scriptstyle{\stackrel{(A_\eta,B_\eta)} {e_{\g \times I}}} &@AA 1_GA\\
&* @>>\stackrel{A}{g_\g} > &*
 \end{CD}}$$
Here $A_\eta=A+z \d(\eta) \in \A^1(M \times I, \lg)$ and $$B_\eta=B +zd\eta+\frac{1}{2}z^2 \eta\wedge^\ad \eta+zA\wedge^\tr \eta+dz \wedge \eta \in  \A^2(M \times I, \le),$$
where { $I=[0,1]$, with coordinate $z$.}
It is an easy calculation to prove that $\d(B_\eta)=\W_{A_\eta}$. In addition, $\g\times I\colon [0,1]^2\to M \times I$ is the map $(\g\times I)(t,s)=(\g(t),s)$, where $s,t \in [0,1]$. {We will see below (Remark \ref{NN}) that ${\stackrel{(A,\eta)}{f_{\g}}}=\stackrel{(A_\eta,B_\eta)} {e_{\g \times I}}$  depends only on $A,\gamma$ and $\eta$. }

Let $h\colon M \to G$ be a smooth map. It is well known (and easy to prove) that if $A''=h^{-1}A'h+h^{-1}dh$ then
$$\tau^h_{A'}(\g)= \quad \quad \quad \begin{CD} &* @>\stackrel{A''}{g_\g} >> &*\\
              &@A h(\g(0)) AA \hskip-1.3cm \scriptstyle{1_E} &@AA h(\g(1))A\\
&* @>>\stackrel{A'}{g_\g} > &*
 \end{CD}$$
{is  a 2-square in $\Gc$.}
{This leads us to the following:
\begin{Definition} We  say that $(A'',B'')$ and $(A,B)$ are related by the 1-gauge transformation $(h,\eta)$, when 
$$ 
A'' = h^{-1}(A+\d (\eta))h + h^{-1}dh
$$ 
and
$$B''=h^{-1} \tr (B+d\eta+A \wedge^\tr \eta +\frac{1}{2}\eta  \wedge^\ad \eta).$$
\label{gtdef}
\end{Definition}
We also define 2-squares relating the holonomies along $\g$ with respect to $A$ and $A''$:
\begin{equation}
 {\tau_{A}^{(h,\eta)}(\g)\doteq \begin{CD}  \tau^h_{A'}(\g)\\{\tau}_{A}^{(1_G,\eta)}(\g) \end{CD}=\quad\quad\quad   \begin{CD} &* @>\stackrel{A''}{g_\g} >> &*\\
              &@A h(\g(0)) AA \hskip-1.5cm \scriptstyle{{\stackrel{(A,\eta)}{f_{\g}}}} &@AA h(\g(1))A\\
&* @>>\stackrel{A}{g_\g} > &*
 \end{CD}} 
\label{tau}
\end{equation}
and 
\begin{equation}
{\hat{\tau}_{A}^{(h,\eta)}(\g)=r_{xy}\left({\tau}_{A}^{(h,\eta)}(\g)\right)=\quad\quad \quad \begin{CD} &* @> \quad  h(\g(1))\quad >> &*\\
              &@A \stackrel{A}{g_\g} AA \hskip-2cm 
\scriptstyle
{
\left(
\stackrel{(A,\eta)} {f_{\g} } 
\right)^{-1}
} 
&@AA  \stackrel{A''}{g_\g}  A\\
&* @>> \quad h(\g(0))  \quad > &*
 \end{CD}} 
\label{tauhat}
\end{equation}
}
see \ref{esdg}.

\begin{Remark}\label{NN}
By the {Non-Abelian Fubini's Theorem}, $\stackrel{(A_\eta,B_\eta)}{e_{\g\times I}}=\stackrel{(A_\eta,B_\eta)}{e_{\g\times I}}(1,1)$, where  $\stackrel{(A_\eta,B_\eta)}{e_{\g\times I}}(t,z)$ can be defined by either of the following differential equations:
$$\frac{\d}{\d z}\stackrel{(A_\eta,B_\eta)}{e_{\g\times I}}(t,z)=-\stackrel{(A_\eta,B_\eta)}{e_{\g\times I}}(t,z)\int_0^t \stackrel{A_z}{g_\g}(t') \tr \eta\left(\frac{\d}{\d t'}\g(t')\right)dt', $$
where {$A_z=A + z \d(\eta) \in \A^1(M,\lg)$,}
or
$$\frac{\d}{\d t}\stackrel{(A_\eta,B_\eta)}{e_{\g\times I}}(t,z)=\left (-z\stackrel{A}{g_\g}(t) \tr \eta\left(\frac{\d}{\d t}\g(t)\right)\right)\stackrel{(A_\eta,B_\eta)}{e_{\g\times I}}(t,z) $$
with initial conditions:
$$ \stackrel{(A_\eta,B_\eta)}{e_{\g\times I}}(\xi,0) =1_E \textrm{ or }  \stackrel{(A_\eta,B_\eta)}{e_{\g\times I}}(0,\xi)=1_E, \textrm{ where } \xi \in [0,1],$$
in the first and second case, respectively. {Therefore it follows that $\stackrel{(A_\eta,B_\eta)}{e_{\g\times I}}$ depends only on $A,\eta$ and $\g$, thus it can be written simply as} {$\stackrel{(A,\eta)}{f_{\g}}$.}
\end{Remark}

{There is another setting for the 2-cubes $\tau$ and $\hat{\tau}$ introduced here, which will be needed when we return to considering local connection pairs $(A_i,B_i)$  ({Definition} \ref{ccg2b}), namely}
$$
\tau_{A_i}^{(\f_{ij},\eta_{ij})}(\g), {\quad \quad \hat{\tau}_{A_i}^{(\f_{ij},\eta_{ij})}(\g) }
$$
where $\g$ is a 1-path whose image is contained in $U_{ij}$. We will refer to these {2-cubes as a transition 2-cubes for the 1-path $\g$.} Note that the relation between $A_i$ and $A_j$ is identical to that between $A$ and $A''$, replacing $h$ by 
$\f_{ij}$ and $\eta$ by $\eta_{ij}$.

\subsubsection{The group of 1-gauge transformations}\label{1gt}

Let $M$ be a smooth manifold. Let also  ${\mathcal G}=(\d\colon E\to  G,\triangleright)$   be a Lie crossed module with associated differential crossed module  $\mathfrak{G}=(\d\colon \le\to \lg,\tr)$. 
The group of 1-gauge transformations {in $M$} is the group of pairs $(h,\eta)$, where $h\colon M \to G$ is smooth, and $\eta$ is an $\le$-valued $1$-form in $M$. The product law will be given by the semidirect product: $(h,\eta)(h',\eta')=(hh',h \tr \eta'+\eta)$.
Recall that a local connection pair in $M$ is given by a pair of forms $A\in \A^1(M,\lg)$ and $B \in \A^2(M,\le)$ with $\d(B)=\W_A=d A +\frac{1}{2}A \wedge^\ad A.$
{Then defining:}
 $$(A,B)\triangleleft (h,\eta)=\left(h^{-1}Ah +\d(h^{-1}\tr \eta)+h^{-1}dh, h^{-1} \tr (B+d\eta+A \wedge^\tr \eta +\frac{1}{2}\eta  \wedge^\ad \eta) \right)$$
{which is equivalent to saying}
$$
{(A'',B'') = (A,B)  \triangleleft (h,\eta)}
$$
{in terms of Definition \ref{gtdef},}
 defines a right action of the group of $1$-gauge transformations {on} the set of local connection pairs.

\subsubsection{The coherence law for 1-gauge transformations}\label{Coher1S}

The following theorem expresses how the holonomy of a local connection pair changes under the group of 1-gauge transformations. We recall the notation of \ref{cs}, \ref{estg} and \ref{1gt}. The notion of a flat $\Gc$-colouring appears in \ref{estg}.

\begin{Theorem}[Coherence law for 1-gauge transformations]\label{Coher1}
{Let $M$ be a smooth manifold with a local connection pair $(A,B)$. Let also $(h,\eta)$ be a 1-gauge transformation, and let $(A'',B'')=(A,B) \triangleleft (h,\eta)$.}
Let $\G\colon [0,1]^2\to M$ be a smooth map. Define
$T_{(A,B)}^{(h,\eta)}(\G)=T_{(A,B)}^{(h,\eta)}$ as being the $\Gc$-colouring of the 3-cube $D^3$ such that:
$$T_{(A,B) }^{(h,\eta)}\circ\delta^-_3 =\stackrel{(A,B)}{\H}(\G), \quad T_{(A,B) }^{(h,\eta)}\circ \delta^+_3=\stackrel{(A'',B'')}{\H}(\G)  $$
and
$$T_{(A,B) }^{(h,\eta)}\circ \delta^\pm_i={\tau_{A }^{(h,\eta)}}(\d^\pm_i \G), \quad i=1,2.$$
(Note that the colourings of the edges of $D^3$ are determined from the colourings of the faces of it, given that they coincide in their intersections.)
Then  $T_{(A,B)}^{(h,\eta)}$ is flat.
\end{Theorem}
\begin{Proof}
The colouring $T_{(A',B')}^{(h,0)}(\G) $ is flat by {Lemma} \ref{gt}{; here $(A',B')=(A,B)\triangleleft (1_G, \eta)$}. Let us prove that the colouring $T_{(A,B)}^{(1_G,\eta)}(\G)$  is flat. This follows  from theorems \ref{Main2} or \ref{Main3} and the fact that if $\M_\eta=dB_\eta+A_\eta\wedge^\tr B_\eta \in \A^3(M \times \{z, z \in \R\},\le)$ is the 2-curvature 3-form of $(A_\eta,B_\eta)$ then the contraction of $\M_\eta$ with the vector field $\frac{\d}{\d z}$ vanishes.  
 A more intricate calculation of this type appears in the proof of Theorem \ref{coher2}. {The theorem follows from the fact that $\Tc(\Gc)$, the set of flat $\Gc$-colourings of the 3-cube $D^3$, is a (strict) triple groupoid (see \ref{estg}) and $T_{(A,B) }^{(h,\eta)}=T_{(A,B) }^{(1_G,\eta)}\circ_3 T_{(A',B')}^{(h,0)} $, where $\circ_3$ denotes upwards composition.}
\end{Proof}

{From remark \ref{NN} it follows:}
\begin{Corollary}\label{VV}
Suppose $\G\colon [0,1]^2 \to M$ is such that $\G(\d [0,1]^2))=x$, where $x\in M$. Given a local connection pair $(A,B)$ in $M$ and a 1-gauge transformation $(h,\eta)$ we then have:
$$\stackrel{(A,B) \triangleleft (h,\eta)} {e_\G}=h^{-1}(x)\tr \stackrel{(A,B)}{e_\G}.$$
\end{Corollary}

{By construction we have:}
\begin{Corollary}
 Given a local connection pair $(A,B)$ in $M$ and a 1-gauge transformation $(h,0)$ we then have for any smooth map $\G\colon [0,1]^2 \to M$:
$$\stackrel{(A,B) \triangleleft (h,0)} {e_\G}=h^{-1}(\G(0,0))\tr \stackrel{(A,B)}{e_\G}.$$
\end{Corollary}

{Theorem \ref{Coher1} may also be interpreted in a different way to give a relation between the holonomies for a 2-path $\G$ with image contained in $U_{ij}$, using local connection pairs $(A_i,B_i)$ and $(A_j,B_j)${; Definition \ref{ccg2b}}. Note that  $(A_j,B_j)=(A_i,B_i)\triangleleft (\f_{ij},\eta_{ij})$.}

\begin{Theorem}[{Transition 3-cube for a {2-path}}]\label{Coher1a}
{Given a connection on a cubical $\Gc$-2-bundle over a pair $(M,\U)$, let $\G\colon [0,1]^2\to M$ be a smooth 2-path with image contained in $U_{ij}$. 
 Define
$T_{(A_i,B_i)}^{(\f_{ij},\eta_{ij})}(\G)=T_{(A_i,B_i)}^{(\f_{ij},\eta_{ij})}$ as being the $\Gc$-colouring of the 3-cube $D^3$ such that:
$$T_{(A_i,B_i)}^{(\f_{ij},\eta_{ij})}\circ\delta^-_3 =\stackrel{(A_i,B_i)}{\H}(\G), \quad T_{(A_i,B_i)}^{(\f_{ij},\eta_{ij})}\circ \delta^+_3=\stackrel{(A_j,B_j)}{\H}(\G)  $$
and
$$T_{(A_i,B_i)}^{(\f_{ij},\eta_{ij})}\circ \delta^\pm_k={\tau_{A_i }^{(\f_{ij},\eta_{ij})}}(\d^\pm_i \G), \quad k=1,2.$$
Then  $T_{(A_i,B_i)}^{(\f_{ij},\eta_{ij})}$ is flat. }
\end{Theorem}

\subsubsection{Dihedral symmetry for 1-gauge transformations}
Let $M$ be a manifold with a local connection pair $(A,B)$ and a 1-gauge transformation $(h,\eta)$. Let $\g\colon [0,1] \to M$ be a smooth map.

\begin{Theorem}\label{dsgt}
We have:
\begin{enumerate}
 \item ${\tau_{A}^{(h,\eta)}(\g^{-1})=\tau_{A}^{(h,\eta)}(\g)^{-\h}}$
 \item If $(A'',B'')=(A,B)\triangleleft (h,\eta)$ then 
${\tau_{A''}^{(h,\eta)^{-1}}(\g)=\left(\tau_{A}^{(h,\eta)}(\g)\right)^{-\vm}.}$
\end{enumerate}
Recall $e^{-\h}=r_x(e)$ and $e^{-\vm}=r_y(e)$, where $e \in \D^2(\Gc)$, denote the  horizontal and vertical inversions of squares in $\Gc$.
\end{Theorem}
\begin{Proof}
The first statement is immediate. Let $h_0=h(\g(0))$,  $h_1=h(\g(1))$ and $\eta'=-h^{-1} \tr \eta$. Let also $(A',B')=(A,B) \triangleleft (0,{\eta})$. The second statement follows from:
$${\begin{CD}\tau_{A''}^{(h,\eta)^{-1}}(\g)\\ \tau_{A}^{(h,\eta)}(\g) \end{CD}\quad=\quad \quad\begin{CD} &\begin{CD} &* @>\quad\stackrel{A}{g_\g}\quad >> &*\\
              &@A h_0^{-1} AA \hskip-1.7cm \scriptstyle{\stackrel{(A'',\eta' )}{f_{\g}}} &@AA h_1^{-1}A\\
&* @>>\quad\stackrel{A''}{g_\g}\quad > &*
 \end{CD}   \\ & \begin{CD} &* @>\quad\stackrel{A''}{g_\g}\quad >> &*\\
              &@A h_0 AA \hskip-1.7cm \scriptstyle{\stackrel{(A,\eta)}{f_{\g}}} &@AA h_1A\\
&* @>>\quad\stackrel{A}{g_\g}\quad > &*
\end{CD}  \end{CD}  \quad=\quad \quad\begin{CD} &\begin{CD} &* @>\quad\quad\stackrel{A}{g_\g}\quad\quad >> &*\\
              &@A 1_GAA \hskip-2.3cm \scriptstyle{h_0^\tr\stackrel{(A'',{\eta'})}{f_{\g}}} &@AA 1_GA\\
&* @>>\quad\quad\stackrel{A'}{g_\g} \quad\quad> &*
 \end{CD}   \\ & \begin{CD} &* @>\quad\quad\stackrel{A'}{g_\g}\quad\quad >> &*\\
              &@A 1_G AA \hskip-2.3cm \scriptstyle{\stackrel{(A,\eta)}{f_{\g}}} &@AA 1_GA\\
&* @>>\quad\quad\stackrel{A}{g_\g}\quad\quad > &*
\end{CD}  \end{CD}  } $$
Now note $${{h_0\tr}\stackrel{(A'',{\eta'})}{f_{\g}}=\stackrel{(A',-\eta)}{f_{\g }}=\left(\stackrel{(A,{\eta})}{f_{\g}}\right)^{-1};}$$
the last equation can be inferred for example from the first equation of Remark \ref{NN}.
\end{Proof}

\subsection{{Equivalence of cubical $\Gc$-2-bundles with connection}}\label{Subdiv}

Let $M$ be a smooth manifold. Let  ${\mathcal G}=(\d\colon E\to  G,\triangleright)$  
be a Lie crossed module and let  $\mathfrak{G}=(\d\colon \le\to \lg,\tr)$ be the associated differential crossed module.  We freely use the material of section \ref{simport}.
\subsubsection{{A crossed module of groupoids of gauge transformations}}

We define a groupoid $M_\Gc^1$, whose set of objects $M_\Gc^0$ is given by the set of local connection pairs $(A,B)$ {in $M$}, in other words $A \in \A^1(M,\lg)$ and $B \in \A^2(M,\le)$ are smooth forms such that $\d(B)=\W_A=dA+\frac{1}{2}A \wedge^\ad A$. The set of morphisms of $M_\Gc^1$ is given by all quadruples of the form $(A,B,\f,\eta)$ where $A$ and $B$ are as above, $\f\colon M \to G$ is a  smooth map and $\eta \in \A^1(M,\le)$ is an $\le$-valued smooth 1-form in $M$. The source of $(A,B,\f,\eta)$ is $(A,B)$ and its target is $(A,B) \triangleleft (\f,\eta)$. The composition is given by the product of 1-gauge transformations; see \ref{1gt}.
We also define a totally intransitive groupoid $M_\Gc^2$, consisting of all triples of the form $(A,B,\psi)$, where $(A,B)$ is a local connection pair in $M$ and $\psi$ is a smooth map $M \to E$. The source and target of $(A,B,\psi)$ each are given by $(A,B)$, and we define $(A,B,\psi)(A,B,\psi')=(A,B,\psi\psi')$. 

{The following lemma states that this gives rise to a crossed {module} of groupoids, a notion defined in \cite{BH1,BHS,B1}, for example. We follow the conventions of \cite{FMPo}.}
\begin{Lemma}
The map $\d\colon M_\Gc^2 \to M_\Gc^1$ such that $$(A,B,\psi) \mapsto \left(A,B,\d \psi,  \psi (d \psi^{-1} )+\psi (A \tr \psi^{-1})\right)$$ is a groupoid morphism, and together with the left action: 
$$(A,B,\f,\eta) \tr (A',B',\psi)=(A,B,\f \tr \psi) ,$$
where $(A',B')=(A,B)\triangleleft (\f,\eta)$, 
of the groupoid $M_\Gc^1$ on the totally intransitive groupoid $M_\Gc^2$ defines a crossed module of groupoids $M_\Gc$.
\end{Lemma}
\begin{Proof}
{Much of this is straightforward calculations. One  complicated bit is to prove that:}
\begin{equation}
{ (A,B)  \triangleleft} \big(\d \psi,  \psi (d \psi^{-1} )+\psi (A \tr \psi^{-1})\big)=(A,B)
\label{abinv}
\end{equation}
{It is easy to see that this is true at the level of 1-forms.} {At the level of the 2-forms we need to prove:}
\begin{multline}\label{formula} B=(\d \psi)^{-1}\tr \Big(B+d(\psi (d \psi^{-1} ))+ d(\psi (A \tr \psi^{-1}))+A\wedge^\tr (\psi (d \psi^{-1} ))+ A\wedge^\tr (\psi (A \tr \psi^{-1})) \\+\frac{(\psi (d \psi^{-1} ) \wedge^\ad (\psi (d \psi^{-1})}{2}+\frac{(\psi (A \tr \psi^{-1})) \wedge^\ad (\psi (A \tr \psi^{-1}))}{2}\\+(\psi (d \psi^{-1} )\wedge^\ad (\psi (A \tr \psi^{-1}))\Big). \end{multline}
{We can eliminate two terms by using:}
$$d(\psi (d \psi^{-1} ))+ \frac{(\psi (d \psi^{-1} ) \wedge^\ad (\psi (d \psi^{-1})}{2}=0,$$
which follows from the fact $d \theta=\frac{1}{2} {\theta \wedge^\ad  \theta}$, where $\theta$ is the Maurer-Cartan form. By using the Leibnitz rule it follows that:
$$  A\wedge^\tr (\psi (A \tr \psi^{-1}))+\frac{(\psi (A \tr \psi^{-1})) \wedge^\ad (\psi (A \tr \psi^{-1}))}{2}=\psi \left (\left( \frac{A\wedge^\ad A}{2}\right ) \tr \psi^{-1}\right).$$
Also we have
$$d(\psi (A \tr \psi^{-1}))+A\wedge^\tr (\psi d \psi^{-1}))+\psi (d \psi^{-1} )\wedge^\ad (\psi (A \tr \psi^{-1}))=\psi (d A \tr\psi^{-1} ) ,$$
{using $\psi(A \tr \psi^{-1})=-(A \tr \psi) \psi^{-1}$ and $(d \psi) \psi^{-1}=-\psi d\psi^{-1}$.}

Putting everything together, formula (\ref{formula}) reduces to: 
\begin{align*}
\f^{-1}\tr &\left (B+\psi \left (\left( \frac{A\wedge^\ad A}{2}\right ) \tr \psi^{-1}\right)+\psi (d A \tr\psi^{-1} )\right)\\&=\f^{-1}\tr \left (B+\psi \left ( \d(B) ) \tr \psi^{-1}\right)\right)\\
&=\f^{-1}\tr \left (B+\psi B \psi^{-1}-B\right)\\
&=B.
\end{align*}
We have used the identity $\d(V) \tr e=Ve-eV$ for each $V \in \le$ and for each $e \in E$. This follows from the definition of a Lie crossed module.

We now prove the other difficult condition, namely:
$$ \d((A,B,\f,\eta) \tr (A',B',\p))=(A,B,\f,\eta)\d((A',B',\p)(A',B',\f^{-1},-\f^{-1} \tr \eta)$$ 
or 
\begin{multline}
(A,B,\d(\f \tr \p),(\f \tr \psi )d(\f \tr \psi )^{-1}+(\f \tr \psi ) A \tr (\f \tr \psi^{-1} ))\\
=(A,B,\f\p\f^{-1}, \eta+{(\f \tr \p) (\f\tr d \p^{-1})+(\f\tr \psi)(\f A' \tr \psi^{-1}) }-\f\d(\p)\f^{-1} \tr \eta)
\end{multline}
{Now use the fact that $A'=\f^{-1}A \f+ \f^{-1} d \f+ \d(\f^{-1} \tr \eta)$, and the terms involving $\eta$ on the right hand side cancel.}
\end{Proof}

\begin{Definition}
 The crossed module of groupoids $M_\Gc$ of the previous lemma will be called the crossed module of gauge transformations in $M$.
\end{Definition}
A very similar construction appears in \cite{SW2}. 
{Note that the collection of crossed modules $U_\Gc$, one  for each open set $U\subset M$, can naturally be assembled into a crossed module sheaf $\overline{M_\Gc}$ over $M$.}
\subsubsection{{Equivalence of cubical $\Gc$-2-bundles with connection over a pair $(M,\U)$}}\label{2bce}
\begin{Definition}
We continue to fix a smooth manifold $M$. Given a point $x \in M$, the crossed module $M_\Gc(x)$ of germs of gauge transformations is constructed in the following obvious way {from the crossed module sheaf $\overline{M_\Gc}$ over $M$}. The set of objects $M_\Gc^0(x)$ of $M_\Gc(x)$ is given by  the set of all triples $(A,B,U)$, with $(A,B) \in U_\Gc^0$, where $U$ is open and $x \in U$, with the equivalence relation $(A,B,U)\cong (A',B',U')$ if $A =A'$ and $B=B'$ in some open neighbourhood of $x$. One proceeds analogously to define the morphisms $M_\Gc^1(x)$ and the 2-morphisms $M_\Gc^2(x)$ of $M_\Gc(x)$. 
\end{Definition}
{Note that the evaluation at $x \in M$ gives  maps} 
$$
{M_\Gc^0(x) \to \Hom(T_x(M),\lg) \times \Hom(\wedge^2(T_x),\le), }
$$ 
$$
{M_\Gc^1(x) \to G \times \Hom(T_x(M),\le) \textrm{ and } M_\Gc^2(x) \to E.}
$$ Therefore the set $\Nc(M_\Gc(x))^n$  of $n$-cubes of the cubical nerve $\Nc(M_\Gc(x))$ of $M_\Gc(x)$ (see \cite{BHS,BHS} and \ref{estg}), comes with a a naturally defined map} 
\begin{multline*}
{t_x\colon\Nc(M_\Gc(x))^n \to \left(\Hom(T_x(M),\lg) \times \Hom(\wedge^2(T_x),\le)\right)^{a_n}} \\ {\times \left (G \times \Hom(T_x(M),\le)\right)^{b_n} \times E^{c_n},}
\end{multline*}
where $a_n,b_n$ and $c_n$ denote the number of vertices, edges and two dimensional faces  of the $n$-cube $[0,1]^n$.

Consider the  bundle $\cup_{x \in M} \Nc(M_\Gc(x))$, of cubical sets,  which is a itself a cubical 
set, where the set of $n$-cubes is given by  $\cup_{x \in M} \Nc(M_\Gc(x))^n$, with the obvious 
faces and degeneracies.  The set of $n$-cubes of $\cup_{x \in M} \Nc(M_\Gc(x))$ can be turned into 
a smooth space \cite{BHo,Ch} by saying that a map $f\colon V \to \cup_{x \in M} \Nc(M_\Gc(x))^n$ is 
smooth if $\left(\cup_{x \in M} t_x \right)\circ f$ is smooth, where $V$ is some open set in some 
$\R^i$. This upgrades the cubical set $\cup_{x \in U} \Nc(M_\Gc(x))$ to a cubical object in the 
category of smooth spaces, a cubical smooth space.

\begin{Theorem}
Let $\U$ be an open cover of $M$. A cubical {$\Gc$-2-bundle} with connection {over $(M,\U)$} is given by a cubical map $C(M,\U,x) \ra{f_x} \Nc(M_\Gc(x))$, the cubical nerve of the crossed module of groupoids $M_\Gc(x)$, for each $x \in M$. This is to verify the following smoothness condition: The collection 
$$
{\bigcup_{x \in M} f_x \colon C(M,\U,x) \to \bigcup_{x \in M} \Nc(M_\Gc(x)) }
$$
is a map of cubical smooth spaces {(recall that $C(M,\U)$ is a cubical manifold)}.
\end{Theorem}
\begin{Proof}
 Easy calculations.
\end{Proof}

\begin{Definition}
{We say that two cubical $\Gc$-2-bundles with connection  $\B$ and $\B'$ over a }pair $(M,\U)$, say $(\f_{ij},\psi_{ijkl}, A_i,B_i,\eta_{ij})$ and $(\f_{ij}',\psi_{ijkl}', A_i',B_i',\eta_{ij}')$, are equivalent (and we write $\B\cong_\U \B'$) if the associated cubical maps $C(M,\U,x) \to {\Nc(M_\Gc(x))}$, where $x \in M$, are homotopic, {through} a smooth homotopy (in the sense above).
\end{Definition}
{The fact that  the {cubical nerve} of a crossed module of groupoids is a Kan cubical set \cite{BH5,BHS} can be used to prove that this is an equivalence relation.}

 Explicitly, $\B\cong_\U \B'$  if there exist smooth maps $\Phi_i\colon U_i \to G$ and $\Psi_{ij}\colon U_{ij} \to E$, as well as smooth forms ${\mathcal E_i}\in  \A^1(U_i,\le)$ such that:
\begin{enumerate}
 \item We have $$\d(A_i,B_i,{\Psi_{ij}^{-1}})(A_i,B_i,\Phi_i,{\mathcal E_i})(A_{i}',B_{i}',\f_{ij}',\eta_{ij}')=(A_i,B_i,\f_{ij},\eta_{ij})(A_j,B_j,\Phi_j,{\mathcal E_j}),$$
{where we suppose $(A'_i,B'_i)=(A_i,B_i)\triangleleft (\Phi_i,{\mathcal E_i}) $ and $(A_j,B_j)=(A_i,B_i)\triangleleft (\f_{ij},\eta_{ij})$.}
\item The colouring $T$ of $D^3$ such that $\d^-_3(T)=(\f,\psi)_{ijkl}$,\quad $\d^+_3(T)=(\f',\psi')_{ijkl}$ (see subsection \ref{Trans2}),  and
$$T^-_1=\quad \begin{CD} &* @>\f'_{ij}>> &*\\
              &@A \Phi_{i}AA \hskip-1.3cm \scriptstyle{\Psi_{ij} }  &@AA \Phi_j A\\
              &* @>> \f_{ij}> &*
  \end{CD} \quad, \quad T^+_1=\quad \begin{CD} &* @>\f'_{kl}>> &*\\
              &@A \Phi_{k}AA \hskip-1.3cm \scriptstyle{\Psi_{kl} }  &@AA \Phi_l A\\
              &* @>> \f_{kl}> &*
  \end{CD}$$ $$T^-_2=\quad \begin{CD} &* @>\f'_{ik}>> &*\\
              &@A \Phi_{i}AA \hskip-1.3cm \scriptstyle{\Psi_{ik} }  &@AA \Phi_k A\\
              &* @>> \f_{ik}> &*
  \end{CD}\quad, \quad T^+_2=\quad \begin{CD} &* @>\f'_{jl}>> &*\\
              &@A \Phi_{j}AA \hskip-1.3cm \scriptstyle{\Psi_{jl} }  &@AA \Phi_l A\\
              &* @>> \f_{jl}> &*
  \end{CD}$$
{is flat for each $x \in U_{ij}$ and any $i,j$; see \ref{estg}. We have put $T^\pm_i=T \circ \delta^\pm_i=\d^\pm_i(T)$.}
\end{enumerate}
We can easily see that this defines an equivalence relation on the set of cubical $\Gc$-2-bundles over $(M,\U)$.

\subsubsection{Subdivisions of covers and the equivalence of cubical $\Gc$-2-bundles over a manifold}

Let $\U=\{U_i\}_{i \in \I}$ be an open cover of $M$. A subdivision $\V$ of $\U$ is a map $i\in \I \mapsto S_i$, where $S_i$ is a set, together with open sets {$V_{a} \subset U_i$,} for each $a\in S_i$ such that $U_i=\cup_{a \in S_i} V_{a}$. If we are given a cubical $\Gc$-2-bundle with connection $\B$ over $C(M,\U)$, we immediately have another one, $\B_\V$ over 
$\V=\{V_{a}\}_{a\in S_i,\, i\in \I}$, 
provided by the obvious cubical map $C(M,\V) \to C(M,\U)$. Its structure maps are such that e.g.  $\f_{ab}=\f_{ij}|_{V_a\cap V_b}$, where $a\in S_i$ and $b\in S_j$, and analogously for all  the remaining information needed to specify a cubical $\Gc$-2-bundle with connection. For the same reason, it is easy to see that if $\B\cong_\U \B'$ then $\B_\V\cong_\V \B_\V'$ for any subdivision $\V$ of $\U$.

If $\U=\{U_i\}_{i \in \I}$ and $\Wc=\{W_j\}_{j \in \J}$ are open covers of $M$, then $\U \cap \Wc$ is the open cover {$\{U_i \cap W_j\}_{(i,j) \in \I \times \J}$.} It is a subdivision of both $\U$ and $\Wc$ in the obvious way.

\begin{Definition}[Equivalence of cubical $\Gc$-2-bundles with connection]
Two cubical $\Gc$-2-bundles with connection $\B$ and $\B'$ over the open covers $\U=\{U_i\}_{i \in \I}$ and $\Wc=\{W_j\}_{j \in \J}$ of $M$, respectively, are called equivalent if $$\B_{\U \cap \Wc}\cong_{\U \cap \Wc} \B'_{\U \cap \Wc}$$
\end{Definition}

The following follows from the previous discussion.
\begin{Theorem}
{
Equivalence of cubical $\Gc$-2-bundles with connection is an equivalence relation.}
\end{Theorem}

\subsection{{Coherence law for transition 2-cubes}\label{Trans2}}

Let $\B$ be a cubical $\Gc$-2-bundle with connection over $(M,\U)$
({Definition} \ref{ccg2b}). Suppose $\g$ is a 1-path whose image is contained in the overlap $U_{ijkl}$. Recall the notation in \ref{esdg}, \ref{estg} and subsection \ref{Trans}, in particular the notion of transition 2-cube for the path $\g$.
{Recall from Definition \ref{cgerbe} the  2-{cube} (for each $x\in M$):}
\begin{equation} 
{(\psi,\phi)}_{ijkl}=\quad \quad \begin{CD}&*@>\f_{kl}>>&*\\
                        &@A \f_{ik}AA \hskip-1.3cm \p_{ijkl} &@AA\f_{jl}A\\          
&*@>> \f_{ij}>&* \end{CD}
\label{psi}
\end{equation}

\begin{Theorem}[Coherence law for {transition 2-cubes}]\label{coher2}
Let $\g\colon [0,1] \to {U_{ijkl}\subset M}$ be a smooth map. We have:
\begin{multline}\label{lhs}\begin{CD}\hat{\tau}^{({\phi}_{ik},\eta_{ik})}_{A_i}(\g) \quad  \hat{\tau}^{({\phi}_{kl},\eta_{kl})}_{A_k}(\gamma)  \quad  {\left ({\hat{\tau}^{({\phi}_{jl},\eta_{jl})}_{A_j}}\right)}^{-\h}(\gamma) \quad {\left({\hat{\tau}^{({\phi}_{ij},\eta_{ij})}_{A_i}}\right)}^{-\h}(\gamma)
\\ \Phi({(\psi,\phi)}_{ijkl}(\g(0)))
\end{CD} \\=\Phi'_{\stackrel{A_i}{g_{\g}}}({(\psi,\phi)}_{ijkl}(\g(1))),\end{multline}
and therefore the $\Gc$-colouring $T$ of $D^3$ such that:
$$T\circ \delta^-_2={(\psi,\phi)}_{ijkl}(\g(0)),  \quad   T\circ \delta^+_2={(\psi,\phi)}_{ijkl}(\g(1)) $$
and	
{\begin{align*}
T\circ\delta^-_1&={\tau}^{({\phi}_{ik},\eta_{ik})}_{A_i}(\g), &T\circ \delta^+_3&= \hat{\tau}^{({\phi}_{kl},\eta_{kl})}_{A_k}(\gamma),\\
 T\circ \delta^+_1 &={\tau}^{({\phi}_{jl},\eta_{jl})}_{A_j}(\gamma),  &T\circ \delta^-_3&= \hat{\tau}^{({\phi}_{ij},\eta_{ij})}_{A_i}(\gamma),
\end{align*}}
is flat.
\end{Theorem}
\begin{Proof}
By {Theorem} \ref{dsgt}, the left hand side $F(\g)$ of (\ref{lhs}) is (we omit the $\g$):
{\begin{align*}
\begin{CD}\hat{\tau}^{({\phi}_{ik},\eta_{ik})}_{A_i} \quad  \hat{\tau}^{({\phi}_{kl},\eta_{kl})}_{A_k}  \quad  \hat{\tau}^{{({\phi}_{jl},\eta_{jl})^{-1}}}
_{A_l} \quad {\hat{\tau}^{({\phi}_{ij},\eta_{ij})^{-1}}_{A_j}}
\\ \Phi({(\psi,\phi)}_{ijkl}(\g(0))
\end{CD}\quad,
\end{align*}}
which can also be written as:
\begin{multline*}
\left [\begin{CD}\hat{\tau}^{(1,\eta_{ik})}_{A_i} \quad  \hat{\tau}^{(1,{\phi}_{ik} \tr \eta_{kl})}_{{\phi}_{ik}\tr A_k}  \quad  \hat{\tau}^{{(1,-{\phi}_{ik} {\phi}_{kl} {\phi}_{jl}^{-1}\tr \eta_{jl})}}
_{{\phi}_{ik} {\phi}_{kl}\tr A_l} \quad 
{\hat{\tau}^{(1,-{\phi}_{ik} {\phi}_{kl} {\phi}_{jl}^{-1}{\phi}_{ij}^{-1} \tr \eta_{ij})}_{{\phi}_{ik} {\phi}_{kl} {\phi}_{jl}^{-1}\tr A_j }}
\\ \id
\end{CD}\right]\quad \quad \circ_\h\\ \left [\begin{CD}{\hat{\tau}^{{\phi}_{ik}}_{{\phi}_{ik} {\phi}_{kl} {\phi}_{jl}^{-1}{\phi}_{ij}^{-1}\tr A_i }}\quad {
{\hat{\tau}^{{\phi}_{kl}}_{ {\phi}_{kl} {\phi}_{jl}^{-1}{\phi}_{ij}^{-1}\tr A_i }}\quad {\hat{\tau}^{{\phi}_{jl}^{-1}}_{ {\phi}_{jl}^{-1}{\phi}_{ij}^{-1}\tr A_i }}\quad 
{\hat{\tau}^{{\phi}_{ij}^{-1}}_{{\phi}_{ij}^{-1}\tr A_i }}\quad }
 \\\Phi({(\psi,\phi)}_{ijkl}(\g(0))\end{CD} \right ]
\end{multline*}
{Here we have put $\f\triangleright A=A\triangleleft \f^{-1}=\f A \f^{-1}+\f d \f^{-1}$.}
Let $\g_t\colon [0,1] \to M$ be the path $\g_t(t')=\g(t't)$, where $t,t' \in [0,1]$. Let also $F'(\g_t)\in E$ be the element assigned to the square $F(\g_t)$. We then have (by using Remark \ref{NN}):

\begin{align*}\frac{d}{d t} F'(\g_t)&=F'(\g_t)\stackrel{A_i}{g_{\g_t}}\tr\left(\eta_{ik}+{\phi}_{ik} \tr \eta_{kl}  -{\phi}_{ik} {\phi}_{kl} {\phi}_{jl}^{-1}\tr \eta_{jl}-{\phi}_{ik} {\phi}_{kl} {\phi}_{jl}^{-1}{\phi}_{ij}^{-1} \tr \eta_{ij} \right)_{\frac{d}{d t} \g(t)}
\\&=F'(\g_t)\stackrel{A_i}{g_{\g_t}}\tr\left(\psi_{ijkl}^{-1}d\psi_{ijkl} +\psi_{ijkl}^{-1}\left(A_i\tr \psi_{ijkl}\right) \right)_{\frac{d}{d t} \g(t)}
 \end{align*}
On the other hand:
\begin{align*}\frac{d}{d t} \left({\stackrel{A_i}{g_{\g_t}}}\tr\psi_{ijkl}(\g(t))\right)
&=\left({\stackrel{A_i}{g_{\g_t}}}  A_i\tr\psi_{ijkl}+{\stackrel{A_i}{g_{\g_t}}}\tr d\psi_{ijkl}\right)_{\frac{d}{d t} \g(t)}\\
&=\left({\stackrel{A_i}{g_{\g_t}}}\tr\psi_{ijkl}\right)\left({\stackrel{A_i}{g_{\g_t}}}\tr\psi_{ijkl}^{-1}\right)\left({\stackrel{A_i}{g_{\g_t}}}  A_i\tr\psi_{ijkl}+{\stackrel{A_i}{g_{\g_t}}}\tr d\psi_{ijkl}\right)_{\frac{d}{d t} \g(t)}\\ 
&=\left({\stackrel{A_i}{g_{\g_t}}}\tr\psi_{ijkl}\right) {\stackrel{A_i}{g_{\g_t}}}\tr\left(\psi_{ijkl}^{-1}d\psi_{ijkl} +\psi_{ijkl}^{-1}\left(A_i\tr \psi_{ijkl}\right) \right)_{\frac{d}{d t} \g(t)}.
\end{align*}
This proves that $F'(\g_t)={\stackrel{A_i}{g_{\g_t}}}\tr\psi_{ijkl}(\g(t))$, which by taking $t=1$ finishes the proof.
\end{Proof}

\section{Wilson spheres and tori}

\subsection{Holonomy for an arbitrary 2-path in a smooth manifold}\label{gencase}
{We recall the notation of subsections \ref{pese}, \ref{Trans} and \ref{Trans2}.}
\subsubsection{Patching together local holonomies and transition functions}\label{Patch}
Let $M$ be a smooth manifold. Let also  ${\mathcal G}=(\d\colon E\to  G,\triangleright)$  be a Lie crossed module with associated differential crossed module $\mathfrak{G}=(\d\colon \le\to \lg,\tr)$.
Let $\U=\{U_i\}_{i \in \I}$ be an open cover of $M$. Let $\B$  be a cubical {$\Gc$-2-}bundle over $(M,\U)$ with connection, given by $\{\f_{ij}, \psi_{ijkl}\}_{i,j,k,l \in \I}$ (Definition \ref{cgerbe}) and $\{A_i,B_i,\eta_{ij}\}_{i,j \in \I}$ (Definition \ref{ccg2b}).  

Let $\G\colon [0,1]^2 \to M$ be a {2-path}. Let $\Q$ denote a subdivision of $[0,1]^2$ into rectangles $\{Q_R\}_{R\in \Re}$, where $\Re$ is some index set, by means of partitions of each $[0,1]$ factor, together with an assignment, to each $R\in \Re$, of $i_R\in \I$, such that $\G(Q_R) \subset  U_{i_R}$. Such subdivisions {with open set assignments (partitions $\Q$ of $\Gamma$) do} exist because of the {Lebesgue Covering Lemma.}

For each $R\in \Re$, let $\G_{R}\colon[0,1]^2\to M$ denote the restriction of $\G$ to $Q_{R}$, rescaled and reparametrized to be a 2-path $[0,1]^2\to M$. We reparametrize again to introduce additional 2-paths, which are thickened 1-paths, constant horizontally (e.g. ${\hat{\g}_{ij}}$ in {Figure {\ref{fig1}} or constant vertically   (e.g. $\g_{ik}$ in Figure \ref{fig1})}, or thickened points, constant both horizontally and vertically ($p_{ijkl}$ in {Figure \ref{fig1}}). To each 2-path in this array we assign a 2-cube of {the double groupoid $\Dc^2(\Gc)$, see \ref{esdg},} as follows:
$$ \G_i \mapsto \H(\G_i) \doteq \stackrel{(A_i,B_i)}{\H}(\G_i)
$$
$$ 
{\hat{\g}_{ij} \mapsto \hat{\tau}({\hat{\g}}_{ij}) }\doteq \hat{\tau}_{A_i}^{(\f_{ij}, \eta_{ij})}(\G_i|_{\{1\}\times [0,1]})
\,\,\, {\rm or}\,\, \,\,
\g_{ik} \mapsto \tau(\g_{ik}) \doteq \tau_{A_i}^{(\f_{ik}, \eta_{ik})}(\G_i|_{[0,1]\times \{1\}}) \,\,\,\,
$$
$$
p_{ijkl} \mapsto \psi(x)_{ijkl} \doteq (\psi, \f)_{ijkl}(x)
$$
where $x\in M$ is the image of the constant 2-path $p_{ijkl}$.
{See Theorem \ref{localH} and equations (\ref{tau}), (\ref{tauhat}) and (\ref{psi}) for the definitions. }

\begin{figure}
\centerline{\relabelbox 
\epsfxsize 9cm
\epsfbox{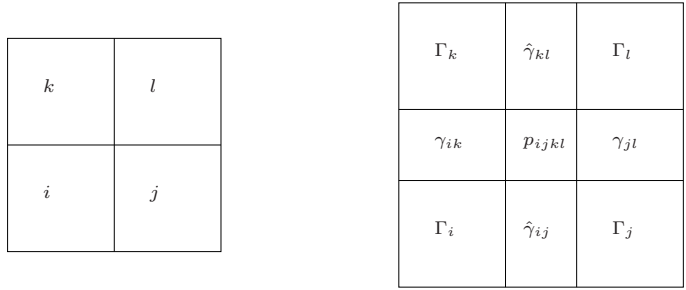}
\relabel{i}{$\scriptstyle{i}$}
\relabel{j}{$\scriptstyle{j}$}
\relabel{k}{$\scriptstyle{k}$}
\relabel{l}{$\scriptstyle{l}$}
\relabel{1}{$\scriptstyle{\G_k}$}
\relabel{2}{$\scriptstyle{\hat{\g}_{kl}}$}
\relabel{3}{$\scriptstyle{\G_l}$}
\relabel{4}{$\scriptstyle{\g_{ik}}$}
\relabel{5}{$\scriptstyle{p_{ijkl} }$}
\relabel{6}{$\scriptstyle{\g_{jl}}$}
\relabel{7}{$\scriptstyle{\G_i}$}
\relabel{8}{$\scriptstyle{\hat{\g}_{ij}}$}
\relabel{9}{$\scriptstyle{\G_j}$}
\endrelabelbox}
\caption{\label{fig1} {Decomposition of $\G$ for the definition of the holonomy of $(\G,\Q)$}}
\end{figure}

\begin{Definition}
\label{lochol}
{Given a 2-path $\G\colon [0,1]^2 \to M$ and a partition $\Q$ of $\G$, the} holonomy of $(\G,\Q)$  for the cubical $\Gc$-2-bundle with connection $\B$, written 
$$  \stackrel{\B}{\H}(\G,\Q),$$
 or simply
${\H}(\G,\Q)$ if the cubical $\Gc$-2-bundle with connection is clear from the context, is the composition of the 2-cubes of $\Dc^2(\Gc)$ obtained from the above assignments. {This is well defined due to the associativity and interchange law for the composition of squares in $\Gc$, which make up a double groupoid; see {\ref{esdg}}.}
\end{Definition}

{In the remainder of this chapter we will see that the 2-dimensional holonomy of {Definition} \ref{lochol} does not depend (up to rather simple transformations) on the chosen partition of $\G$, the chosen coordinate neighbourhoods, the choice of cubical $\Gc$-2-bundle with connection within the same equivalence class, or the choice of $\G$ within the same thin homotopy equivalence class.  Furthermore, in the final section we will see how it can be associated to oriented embedded 2-spheres in a manifold, therefore defining Wilson {2-Sphere} observables. 
}

\subsubsection{Independence under {\it subdividing} partitions	}\label{subdivide}
\begin{Proposition}\label{iusp}
 Suppose we introduce an extra point in one of the partitions underlying $\Q$, so as to subdivide one of the rows or columns of the partition of $[0,1]^2$. For this new subdivision, suppose we assign each of its rectangles to the same open set as that assigned by $\Q$ to the rectangle in which it is contained, and call this new subdivision and assignment $\Q'$. Then
$$
{\H}(\G,\Q')= {\H}(\G,\Q)
$$
\end{Proposition}
\begin{Proof}
(For the case of subdividing a row). The only change in the holonomy for $\Q'$ is in the contributions along the subdivided row, where the open set assignments look like {Figure  {\ref{fig1}}} with $i=k$ and $j=l$. Since $\tau(\g_{ii})$ and $\psi_{ijij}(p)$ are thin elements of  $\Dc^2(\Gc)$ (from  Definition \ref{ccg2b} and Section \ref{Trans}, and from Definition \ref{cgerbe} respectively), the composition of the three rows of rectangles after subdividing equals the composition of the original row of rectangles before subdividing. 
\end{Proof}

\subsubsection{The case of paths}\label{paths}
Let $\g: [0,1]\rightarrow M$ be a path. Let $\Q$ denote a subdivision of $[0,1]$ into subintervals $\{q_r\}_{r=1,\dots. s}$, together with an assignment, for each $r$, of $i_r\in \I$, such that $\g(q_r)\subset 
U_{i_r}$. For each $r$, let $\g_r:[0,1]\rightarrow M$ denote the restriction of $\g$ to $q_r$, rescaled and reparametrized to be a 1-path $[0,1]\rightarrow M$. As for the case of 2-paths, we reparametrize again to introduce constant 1-paths $p_{r,r+1}$ with image $x_r= \g_r(1)=\g_{r+1}(0)$ between $\g_r$ and $\g_{r+1}$. To each of these 1-paths we assign an element of $G$ as follows:
$$ \g_r \mapsto \stackrel{A_{i_r}}{g_{\g_r}}
$$
$$ p_{r,r+1} \mapsto \f_{i_r i_{r+1}}(x_r)
$$

\begin{Definition}
\label{locholpath}
The holonomy of $(\g,\Q)$ for the cubical $\Gc$-2-bundle with connection $\B$, written 
$$  \stackrel{\B}{\H}(\g,\Q),$$
 or simply
${\H}(\g,\Q)$ if the cubical $\Gc$-2-bundle with connection is clear from the context, is the composition of the 1-cubes of $\Dc^1(\Gc)$ obtained from the above assignments. Concretely, we have the formula:
$$
\stackrel{\B}{\H}(\g,\Q)=\stackrel{A_{i_1}}{g_{\g_1}} \f_{i_1i_2}(x_1)\stackrel{A_{i_2}}{g_{\g_2}} \f_{i_2i_3}(x_2)\ldots \stackrel{A_{i_s}}{g_{\g_s}}.
$$
\end{Definition}

\vskip 20 pt

Let $\g$ be a 1-path, and let $\Q$, $\Q'$ be based on the same subdivision of $[0,1]$ into subintervals $\{q_r\}_{r=1,\dots. s}$, but with different assignments $i_r$ and $i'_r$ to each $q_r$. As in {Definition} \ref{lochol}, we replace $\g$ by a product of 2-paths which are constant vertically, corresponding to $\g_r$, or constant horizontally and vertically, corresponding to $x_r$. We introduce the notation:
\begin{equation}
\label{taunotation}
 \stackrel{\B}{\tau}(\g,\Q,\Q')  \doteq \tau((\g_1)_{i_1 i'_1}) \psi(x_1)_{i_1 i_2 i'_1 i'_2} \tau((\g_2)_{i_2 i'_2}) \ldots
 \tau((\g_s)_{i_s i'_s}).
\end{equation}
When $\B$ is understood we will drop it from the notation. In particular, this denotes the evaluation of a row of the holonomy formula of {Definition} \ref{lochol}, with $\g$ being the restriction of $\G$ to one of the horizontal lines in the partition of $[0,1]^2$. We have:
$$
\d_d\tau(\g,\Q,\Q')=\H(\g,\Q)\, \textrm{ and }\, \d_u\tau(\g,\Q,\Q')=\H(\g,\Q'),
$$
with $\B$ understood everywhere.

\vskip 10 pt

\subsubsection{The dependence of the holonomy on {the partition} $\Q$}

We want to study the effect on the holonomy of substituting the subdivision with open set assignments $\Q$ by $\Q'$. Since by the previous proposition, the holonomy is unaffected by subdividing the partition of $[0,1]^2$, we can assume {that the underlying subdivision of $[0,1]^2$ is the same for $\Q$ and $\Q'$, thus that $\Q$ and $\Q'$  differ { only with respect to} the open set assignments.}

\begin{Theorem}[Coherence law for 2-holonomy]\label{gc2hg} Let $\G\colon [0,1]^2\to M$ be a smooth map.
Suppose $\Q$ and $\Q'$ are given by the same subdivision of $[0,1]^2$ into rectangles $\{Q_R\}_{R\in \Re}$, and assignments $i_R$ and $i'_R$ respectively to each rectangle $Q_R$ such that $\G(Q_R) \subset  U_{i_R} \cap U_{i'_R}$. Then the respective holonomies of $\G$  are related by the homotopy addition equation {(\ref{ha})} for $T\in D^3$, where $T$ is given by:
$$T\circ \delta^-_3 =\H(\G,\Q) \textrm{ and } T\circ \delta^+_3 =\H(\G,\Q') $$
and 
$$T \circ \delta^\pm_i=\tau(\d^\pm_i(\G), \d^\pm_i \Q,  \d^\pm_i \Q'), \quad i=1,2;$$
where $\d^\pm_i \Q$ and $\d^\pm_i \Q'$ are the restrictions of $\Q$ and $\Q'$ to the corresponding faces. {In other words the colouring $T$ of $D^3$ is flat.}
\end{Theorem}

\begin{Proof}
Analogously to the procedure in {Definition} \ref{lochol}, but now in three dimensions, we take the 3-path $\G \times {\rm id}_{[0,1]}$, with its domain $[0,1]^3$ partitioned into rectangular solids by the partition of the domain of $\G$ underlying $\Q$ and $\Q'$. We then reparametrize to replace the vertical surfaces and lines of the partition by 3-paths that are constant horizontally or vertically, or both horizontally and vertically. The flat cube {$T\in \T^3(\Gc)$} is the composition of elementary flat cubes of the following types.

 To each 2-path $\G_R$, we assign (see Theorem \ref{Coher1a})
$$
T(\G_R, \Q, \Q') = T_{(A_{i_R}, B_{i_R})}^{(\f_{i_R {i'}_R}, \eta_{i_R {i'}_R})} 
$$ 
{To each $\hat{\g}_{RS}=\d_1^+ \G_R = \d_1^- \G_S$ we assign a version of the flat cube of Theorem \ref{coher2}, namely $T(\hat{\g}_{RS}, \Q, \Q')$ given by}
$$
\begin{cases}
\d_3^-  T(\hat{\g}_{RS}, \Q, \Q') = \hat{\tau}(\g_{i_R i_S}), \quad \d_3^+  T(\hat{\g}_{RS}, \Q, \Q') = \hat{\tau}(\g_{{i'}_R {i'}_S})\\
\d_1^-  T(\hat{\g}_{RS}, \Q, \Q') = {\tau}(\g_{i_R {i'}_R}),\quad \d_1^+  T(\hat{\g}_{RS}, \Q, \Q') = {\tau}(\g_{{i}_S {i'}_S})\\
\d_2^-  T(\hat{\g}_{RS}, \Q, \Q') = \psi(\g(0)_{i_R i_S {i'}_R {i'}_S}), \quad 
\d_2^+  T(\hat{\g}_{RS}, \Q, \Q') = \psi(\g(1)_{i_R i_S {i'}_R {i'}_S})\\
\end{cases}
$$
To each $\g_{RS}=\d_2^+ \G_R = \d_2^- \G_S$ we assign a version of the flat cube of {Theorem} \ref{coher2}, namely $T(\g_{RS}, \Q, \Q')$ given by
$$
\begin{cases}
\d_3^-  T(\g_{RS}, \Q, \Q') = {\tau}(\g_{i_R i_S}), \quad \d_3^+  T(\g_{RS}, \Q, \Q') = {\tau}(\g_{{i'}_R {i'}_S})\\
\d_2^-  T(\g_{RS}, \Q, \Q') = {\tau}(\g_{i_R {i'}_R}),\quad \d_2^+  T(\g_{RS}, \Q, \Q') = {\tau}(\g_{{i}_S {i'}_S})\\
\d_1^-  T(\g_{RS}, \Q, \Q') = \psi(\g(0)_{i_R i_S {i'}_R {i'}_S}), \quad 
\d_1^+  T(\g_{RS}, \Q, \Q') = \psi(\g(1)_{i_R i_S {i'}_R {i'}_S})\\
\end{cases}
$$
Finally, to each $p_{RSTU}= \d_2^+ \d_1^+ \G_R = \d_2^+ \d_1^- \G_S = \d_2^- \d_1^+ \G_T = \d_2^- \d_1^- \G_U $, we assign the flat cube of Definition \ref{cgerbe} (2) with open set indices $i_R, i_S, i_T, i_U$ and $ {i'}_R, {i'}_S, {i'}_T, {i'}_U$. {The result follows from the fact that the set of flat 3-cubes in $\Gc$ can be turned into a strict triple groupoid; see \ref{estg}.}
\end{Proof}

As an immediate consequence we have {the following non-trivial result}:
\begin{Corollary}
Let $\G, \, \Q, \, \Q'$ be as in {Theorem} \ref{gc2hg}. If the {open set} assignments $i_R$ and $i'_R$ agree on the rectangles along the boundary of $[0,1]^2$, then $\H(\G,Q)= \H(\G, \Q')$. 
\end{Corollary}
\begin{Proof}
{If we use condition 4 of {Definition} \ref{cgerbe} and condition 1 of {Definition} \ref{ccg2b} in equation (\ref{taunotation}) we can see that $T \circ \delta_i^\pm$ each are identity 2-cubes in {$\Gc$} for $i=1,2$. Now  compare with the homotopy addition equation (\ref{ha}).} 
\end{Proof}

{Analogously it follows:}
\begin{Corollary}\label{sph}
Let $\G, \, \Q, \, \Q'$ be as in {Theorem} \ref{gc2hg}. Suppose $\G(\d [0,1]^2)=x$, for some $x\in M$, and that the open set assignments for all rectangles along the boundary of $[0,1]^2$ are chosen to be the same, i.e. all equal to $i_x$ for $\Q$ and all equal to ${i'}_x$ for $\Q'$. Then we have:
$$\H(\G,\Q)=\left(\f_{i_x{i'}_x}(x)\right)^{-1} \tr \H(\G,\Q).$$
\end{Corollary}

\subsubsection{Invariance under (free) thin homotopy}

Let $M$ be a manifold with a local connection pair $(A,B)$. It follows from {Theorem} \ref{Main3} that the two dimensional holonomy $\stackrel{(A,B)}{\mathcal H}(\G)$, where $\G\colon [0,1]^2 \to M$  is a smooth path, is invariant under thin homotopy. Now suppose that $M$ is equipped with a cubical $\Gc$-2-bundle connection. In this subsection we will study how $\H(\G)$ varies under thin homotopy.  We will consider a slightly more general definition of thin homotopy (a generality that is needed to define Wilson spheres).

\begin{Definition}
Two smooth maps $\G,\G'\colon [0,1]^2 \to M$ are said to be freely thin homotopic if there exists a smooth {map} $J\colon [0,1]^2 \times [0,1] \to M$ such that  $\Rank (\D_v J)\leq 2, \forall v \in [0,1]^3$, and such that $\d_3^-J=\G$ and $\d_3^+J=\G'$.
\end{Definition}
Note that $J$ is, in general, not a rank-2 homotopy since it does not satisfy the conditions 1 and 2 of its definition; see \ref{2-Tracks}.

\begin{Theorem}[Invariance under free thin homotopy]\label{Invthin}
Consider a free thin homotopy $J: [0,1]^3\rightarrow M$ with $\d_3^- J = \G$ and $\d_3^+ J = \G'$. Let $\Q$ denote a subdivision of 
$[0,1]^3$ into rectangular solids $\{Q_R\}_{R\in \Re}$, using partitions of the three $[0,1]$ factors, together with an assignment for each $R\in \Re$ of $i_R \in \I$ such that  $J(Q_R) \subset  U_{i_R}$. Such subdivisions exist because of the {Lebesgue Covering Lemma. Then} $\Q$ naturally induces subdivisions and open set assignments on each face of $[0,1]^3$, denoted $\d_i^\pm \Q$, $i=1,2,3$. 

Then the holonomies $\H(\G, \d_3^- \Q)$ and $\H(\G', \d_3^+ \Q)$, with respect to a fixed cubical {$\Gc$-2-}bundle with connection $\B$, are related by the homotopy addition equation {(\ref{ha})} for $T(J,\Q)$, where:
$$
\d_i^\pm T(J,\Q) = \H(\d_i^\pm J, \d_i^\pm \Q), \, i=1,2,3.
$$
\end{Theorem}
\begin{Proof}
{The proof is very similar {to} the proof of Theorem \ref{gc2hg}.}
By analogy with the definition of holonomy, we reparametrize $J$ to introduce additional 3-paths for each face separating the rectangular solids, for each edge separating these faces and for each point separating these edges. The additional 3-paths are constant in one, two or all three of the directions (horizontal, vertical, upwards). {The cube} $T(J,\Q)$ is the composition of flat cubes of various types which, for the most part, we have already encountered in the proof of Theorem \ref{gc2hg}, or are analogous versions of these obtained by rotation. The remaining flat cubes are of the type appearing in Theorem \ref{Main3}, corresponding to $J_R$, the restriction of $J$ to $Q_R$, reparametrized to be a 3-path, with the local connection pair $(A_{i_R}, B_{i_R})$, for each $R\in \Re$. Note that the curvature 3-form vanishes, since $J$ is thin.
\end{Proof}

The following analogue of {Corollary} \ref{sph} holds.
\begin{Corollary}\label{shpp}
Under the conditions of {Theorem} \ref{Invthin}, suppose $J$ is such that $J(\d [0,1]^2\times \{t\})=q(t)$, for some smooth map $q\colon [0,1] \to M$, with $q(0)=x$ and $q(1)=x'$. Suppose also that the open set assignments for the rectangular solids along $\d [0,1]^2 \times [0,1]$ only depend on the upwards direction, i.e. they are given by fixing $\d_1^- \d_2^- \Q$. Then
$$
\H(\G,  \d_3^+ \Q)=\left( \H(q, \d_1^- \d_2^- \Q)\right)^{-1} \tr \H(\G,  \d_3^- \Q),
$$
where $\H(q, \d_1^- \d_2^- \Q)$ is defined in Definition \ref{locholpath}.
\end{Corollary}

\subsubsection{Dihedral symmetry for the holonomy of general squares}

{Suppose that} $\B$ is a dihedral cubical $\Gc$-2-bundle over $(M,\U)$, with a dihedral cubical connection (see definitions \ref{Dih} and \ref{ccdg2b}). Let $(\G,\Q)$ be as in Definition \ref{lochol}, and let $r$ be some element of the dihedral group $D_4$ of the square. Then we define $\Q^r$ to be the subdivision of $[0,1]^2$ with open set assignments induced on $\G \circ r^{-1}$ by $\Q$. 
\begin{Theorem}\label{dshgs}
 We have:
$$\H(\G\circ r^{-1}, \Q^r )=r(\H(\G, \Q)).$$
\end{Theorem}
\begin{Proof}
This follows from theorems \ref{NAFT} and  \ref{dsgt} and the definition of a dihedral cubical {$\Gc$-2-bundle} with a dihedral connection; definitions \ref{Dih} and \ref{ccdg2b}. Note that the action of $r$ in $\D^2(\Gc)$ is  a double-groupoid morphism (see \ref{esdg}), so that it is enough to check the equation for all the 2-paths appearing in the definition of $\H(\G, \Q)$ and the corresponding 2-cubes of $\Dc^2(\Gc)$ - see Definition \ref{lochol}.
\end{Proof}

\subsubsection{{Dependence of the surface holonomy on the cubical $\Gc$-2-bundle with connection equivalence class}}\label{EQUIV}

Let $\B$ be a cubical $\Gc$-2-bundle with connection over $(M,\U)$, and recall from subsection \ref{Subdiv} the cubical $\Gc$-2-bundle with connection $\B_\V$ obtained from $\B$ and a subdivision $\V$ of the cover $\U$. Consider the holonomy $\stackrel{\B}{\H}(\G,\Q)$ of {Definition} \ref{lochol}. Let $\Q_\V$ denote the same subdivision of $[0,1]^2$ into rectangles $\{Q_R\}_{R\in \Re}$ as $\Q$, with assignments $R\mapsto a_R$ such that $\G(Q_R) \subset  V_{a_R}$, where $a_R\in  S_{i_R}$ (using the notation at the end of subsection \ref{Subdiv}). Then it is clear from Definition \ref{lochol} {and Proposition \ref{iusp}} that we have:
$$
\stackrel{\B_\V}{\H}(\G,\Q_\V) =
\stackrel{\B}{\H}(\G,\Q).
$$
Thus we will only consider equivalences of cubical $\Gc$-2-bundles with connection with respect to a fixed cover $\U$ of $M$.

Suppose that $\B$ and $\B'$ are equivalent cubical $\Gc$-2-bundles with connection, with the equivalence given by the triple 
$(\Phi_i, {\mathcal E_i}, \Psi_{ij})$ of subsection \ref{2bce}. Note that condition (1) of the equivalence, in view of  equation (\ref{abinv}), may be rewritten as the following equations:
$$
 (A'_i,B'_i)= (A_i,B_i) \triangleleft (\Phi_i,{\mathcal E_i}) \quad \quad  (A'_j,B'_j)=(A'_i,B'_i)\triangleleft (\f'_{ij},\eta'_{ij})
$$
$$
 (A_j,B_j)=(A_i,B_i)\triangleleft (\f_{ij},\eta_{ij}) \quad \quad  (A'_j,B'_j)= (A_j,B_j) \triangleleft (\Phi_j,{\mathcal E_j}) 
$$

We now proceed analogously to equation (\ref{taunotation}). 
Let $\g$ be a 1-path, and let $\Q$, be a subdivision of $[0,1]$ into subintervals $\{q_r\}_{r=1,\dots, s}$, with an  assignment $r\mapsto i_r\in \I$, such that $\g(q_r)\subset U_{i_r}$. Let $\g_r$ denote the restriction of $\g$ to $q_r$, rescaled and reparametrized to be a 1-path,  and denote the points separating the images of $\g_r$ by $x_r$.  We define:
\begin{equation}
\label{snotation}
 \stackrel{(\B,\B')}{s}(\g,\Q)  \doteq 
{\tau}^{(\Phi_{i_1},\mathcal E_{i_1})}_{A_{i_1}} (\g_{1}) 
(\Psi,\Phi)_{i_1i_2}(x_1) 
{\tau}^{(\Phi_{i_2},\mathcal E_{i_2})}_{A_{i_2}} (\g_{2}) 
(\Psi,\Phi)_{i_2i_3}(x_2) 
 \ldots
{\tau}^{(\Phi_{i_s},\mathcal E_{i_s})}_{A_{i_s}} (\g_{s}) .
\end{equation}
Then {the proof of Theorem} \ref{gc2hg} can be reformulated to give the dependence of the holonomy on changing $\B$ within the same equivalence class.

\begin{Theorem}[Behaviour under cubical $\Gc$-2-bundle equivalences]\label{gc2hg2} 
Let $\B$ and $\B'$ be  equivalent cubical $\Gc$-2-bundles with connection, with the equivalence given by the triple 
$(\Phi_i, {\mathcal E_i}, \Psi_{ij})$. Let $\G\colon [0,1]^2\to M$ be a smooth map and suppose $\Q$ is a  subdivision of
 $[0,1]^2$ into rectangles $\Q = \{Q_R\}_{R\in \Re}$, together with assignments $R\mapsto i_R$  such that $\G(Q_R) \subset  U_{i_R}$. Then the holonomies of $(\G,\Q)$ with respect to $\B$ and $\B'$ are related by  
 the homotopy addition equation {(\ref{ha})} for $T\in D^3$, where $T$ is given by:
$$T\circ \delta^-_3 =\stackrel{\B}{\H}(\G,\Q)
 \textrm{ and } T\circ \delta^+_3 = \stackrel{\B'}{\H}(\G,\Q)
$$
and 
$$T \circ \delta^\pm_i=
 \stackrel{(\B,\B')}{s}( \d^\pm_i(\G), \d^\pm_i \Q   ).
$$
\end{Theorem}

We have the following analogue of Corollary \ref{sph}.
\begin{Corollary}\label{sph2}
Given the conditions of Theorem \ref{gc2hg2}, suppose $\G(\d
[0,1]^2)=x$, for some $x\in M$, and that the open set 
assignments for all rectangles along the boundary of $[0,1]^2$ are
chosen to be the same, say $i_x$. Then 
 $$
\stackrel{\B'}{\H}(\G, \Q)= \left(\Phi_{i_x}(x)\right)^{-1} \tr \stackrel{\B}{\H}(\G, \Q).
$$
\end{Corollary}

\subsection{Two types of Wilson surfaces}

Let $\B$ be a cubical $\Gc$-2-bundle with connection over
$(M,\U)$. Let  $\G\colon [0,1]^2\to M$ be a 2-path such that  $\G(\d
[0,1]^2)=x$ for some $x\in M$. Thus $\G$ factors through a map $f:S^2
\rightarrow M$. We say that $\G$ and $\G'$ are equivalent if the
corresponding maps $f$ and $f'$ from $S^2$ to $M$ are related by $f' =
f \circ g$ where $g$ is an orientation-preserving {diffeomorphism} of
$S^2$. 

Let $\Q$ be a subdivision of $[0,1]^2$ into rectangles
$\{Q_R\}_{R\in \Re}$ with open set assignments $R \mapsto i_R$ such
that $\G(Q_R) \subset  U_{i_R}$, and suppose that these assignments
are the same, say $i_x$, for all rectangles along the boundary of 
$[0,1]^2$. 

\begin{Definition}
With $\B,\, \G$ and $\Q$ as above, we define the Wilson sphere
functional to be 
$$
{\mathcal W}_\B(\G,\Q)= \stackrel{\B}{\H}(\G, \Q)  \in \ker \d\subset E.
$$
\end{Definition}

\begin{Theorem}\label{WP}
Up to acting by elements of $G$, the Wilson sphere functional
${\mathcal W}_\B(\G,\Q)$ is independent of the choice of $\Q$, the
  choice of $\G$ within the same equivalence class, and the choice of
  $\B$ within the same equivalence class. 
For $\B$ a dihedral bundle with dihedral connection and $r\in D_4$ an
orientation {reversing} element, we have, following the notation of
Theorem \ref{dshgs}, 
$$
{\mathcal W}_\B(\G\circ r^{-1},\Q^r) =
\left( {\mathcal W}_\B(\G,\Q) \right)^{-1}.
$$
\end{Theorem}
\begin{Proof}
The statement for $\Q$ follows from subsection \ref{subdivide} and Corollary
\ref{sph}. Since the mapping class group of $S^2$ is $\{\pm 1\}$, when
$\G$ and $\G'$ are equivalent, then they are isotopic. Thus there exists a thin free homotopy
$J:[0,1]^3\rightarrow M$ of the type appearing in Corollary
\ref{shpp} ($J$ is thin since {it factors through a smooth family of diffeomorphisms of $S^2$)}, and satisfying $\d_3^- J = \G$ and $\d_3^+ J = \G'$. Thus the
statement for $\G$ follows from Corollary \ref{shpp}. The statement
for $\B$ follows from Corollary \ref{sph2}. The final statement, when
the bundle and connection are {dihedral}, is an immediate consequence
of Theorem \ref{dshgs}.
\end{Proof}

If the image of $\G$ is an embedded sphere $\Sigma$ in $M$, then any two orientation-preserving parametrizations of $\Sigma$ are equivalent. In this case we may state the result as follows:
\begin{Theorem}[Embedded Wilson Spheres]
The holonomy of an oriented embedded sphere $\Sigma$ does not depend on the chosen parametrization of $\Sigma$ up to acting by elements of $G$. We denote it by ${\mathcal W}_\B(\Sigma)$. 
\end{Theorem}
 This may have  applications in 2-knot theory, {c.f. \cite{W,CR}}.

With  $\B$ as before, suppose now that $\G$ is such that $\d_u \G =
\d_d \G$ and $\d_l \G = \d_r \G$. Then the 2-path $\G$ factors through
a map $f$ from the torus $T^2$ to $M$.  We say that $\G$ and $\G'$ are equivalent if the
corresponding maps $f$ and $f'$ from $T^2$ to $M$ are related by $f' =
f \circ g$ where $g$ is an automorphism of
$T^2$ which is isotopic to the identity (note that the mapping class
group of the torus is $\GL(2,\Z)$). 

Let $\Q$ be a subdivision of $[0,1]^2$ into rectangles
$\{Q_R\}_{R\in \Re}$ with open set assignments $R \mapsto i_R$ such
that $\G(Q_R) \subset  U_{i_R}$, and suppose that these assignments
are such that they match along the upper and lower boundary of
$[0,1]^2$, and along the left and right boundary of $[0,1]^2$, i.e. $\d_u \Q =
\d_d \Q$ and $\d_l \Q = \d_r \Q$.

\begin{Definition}
With $\B,\, \G$ and $\Q$ as above, we define the Wilson torus
functional to be 
$$
{\mathcal W}_\B(\G,\Q)= \stackrel{\B}{\H}(\G, \Q)  \in  \d^{-1}(G^{(1)})\subset E,
$$
where $G^{(1)}$ is the commutator subgroup of $G$.
\end{Definition}

Note that the value of the Wilson torus functional indeed belongs to
${\d^{-1}\big (G^{(1)}\big)}$, since
$$
\d (\stackrel{\B}{\H}(\G, \Q)) = [ 
\stackrel{\B}{\H}(\d_d\G, \d_d\Q)),
\stackrel{\B}{\H}(\d_r\G, \d_r\Q))
].
$$

Analogous arguments to the proof of Theorem \ref{WP}, now using Theorem \ref{gc2hg}, Theorem \ref{Invthin} and Theorem
\ref{gc2hg2}, give:
\begin{Theorem}
The Wilson torus functional ${\mathcal W}_\B(\G,\Q)$ is independent of the choice of $\Q$, the
  choice of $\G$ within the same equivalence class , and the choice of
  $\B$ within the same equivalence class, up to changes of the form 
 of the following simultaneous horizontal and vertical conjugation:
$${\mathcal W}_\B(\G,\Q)\mapsto {\begin{array}{ccc}\ulcorner &  e_2^{-\vm} &\urcorner \\ e_1 & {\mathcal W}_\B(\G,\Q) & e_1^{-\h} \\ 
\llcorner & e_2 & \lrcorner \end{array}}.$$
\end{Theorem}

\begin{Remark}
If the image of $\G$ is an embedded torus $\Sigma$ in $M$, then unlike in the case of the sphere, the holonomy of $\Sigma$ will in general depend on the mapping class of $\G$ and not just on the oriented embedded surface itself. {This is a consequence of the fact that the mapping class group of the torus is ${\rm GL}(2,\Z)$ rather than $\{\pm 1\}$, which is the case of the sphere.}
\end{Remark}

\section*{Acknowledgements}
{This work was  partially supported by the {\em Programa Operacional  Ci\^{e}ncia e Inova\c{c}\~{a}o 2010},   and the 
{\em Projecto 3599: Promover a Produ\c{c}\~{a}o Cientifica, o Desenvolvimento Tecnol\'{o}gico e a Constitui\c{c}\~{a}o de Redes Tem\'{a}ticas}
financed by the  {\em Funda\c{c}\~{a}o para a Ci\^{e}ncia e a Tecnologia} (FCT) and  cofinanced by the European Community fund FEDER, in part through the FCT research project Quantum Topology (POCI/MAT/60352/2004 and PPCDT/MAT/60352/2004).}
{We would like to thank Jim Stasheff and {Tim Porter} for useful comments.}

\end{document}